\documentclass[mathpazo]{cicp}

\newcommand*\diff[1]{\mathop{}\!{\mathrm{d}#1}} %% d in integrand

\newcommand{\td}[2]{\frac{\diff #1}{\diff #2}}
\newcommand{\pd}[2]{\frac{{\partial}#1}{{\partial}#2}}

\newcommand{\mb}[1]{\mathbf{#1}}
\newcommand{\mbb}[1]{\mathbb{#1}}
\newcommand{\mc}[1]{\mathcal{#1}}
\newcommand{\nor}[1]{\left\| #1 \right\|}

\newcommand{\LRp}[1]{\left( #1 \right)}
\newcommand{\LRs}[1]{\left[ #1 \right]}

\newcommand{\LRb}[1]{\left| #1 \right|}
\newcommand{\LRc}[1]{\left\{ #1 \right\}}
\newcommand{\LRlim}[1]{\left.{ #1 }\right|}

\newcommand{\grad}{\nabla}
\renewcommand{\div}{\grad \cdot}

\usepackage{subfig}
\usepackage{hhline}
\usepackage{amsmath}
\usepackage{graphicx}
\usepackage{amsfonts}
\usepackage{hyperref}
\usepackage{color}
\usepackage{array}
\usepackage{algorithm}
\usepackage[noend]{algpseudocode}
\usepackage{algorithmicx}

\begin{document}

\title{Variations on Hermite methods for wave propagation}

%\author{Arturo Vargas\thanks{Department of Computational and Applied Mathematics, Rice University, 6100 Main Street, Houston, TX 77005 ({\tt av29@rice.edu}).}
%        \and Jesse Chan\thanks{Department of Mathematics, Virginia Tech, Blacksburg, VA 24060} \and Thomas Hagstrom\thanks{Department of Mathematics, Southern Methodist University, Dallas, TX} \and T. Warburton$^{\ddagger}$}

%\author{First Author         \and
%        Second Author %etc.
%}
%
%%\authorrunning{Short form of author list} % if too long for running head
%
%\institute{F. Author \at
%              first address \\
%              Tel.: +123-45-678910\\
%              Fax: +123-45-678910\\
%              \email{fauthor@example.com}           %  \\
%%             \emph{Present address:} of F. Author  %  if needed

\author{Arturo Vargas\affil{1}, Jesse Chan\affil{2}\corrauth, Thomas Hagstrom\affil{3}, and T. Warburton\affil{2}}
\address{\affilnum{1}\ {Department of Computational and Applied Mathematics, Rice University.}\\
\affilnum{2}{Department of Mathematics, Virginia Tech}\\
\affilnum{3}{Department of Mathematics, Southern Methodist University}
}
\email{{\tt jlchan@vt.edu}}

%\tableofcontents

\begin{abstract}
Hermite methods, as introduced by Goodrich et al.\ in \cite{goodrich2006hermite}, combine Hermite interpolation and staggered (dual) grids to produce stable high order accurate schemes for the solution of hyperbolic PDEs. We introduce three variations of this Hermite method which do not involve time evolution on dual grids.  Computational evidence is presented regarding stability, high order convergence, and dispersion/dissipation properties for each new method.  Hermite methods may also be coupled to discontinuous Galerkin (DG) methods for additional geometric flexibility \cite{chen2014hybrid}.  An example illustrates the simplification of this coupling of this coupling for the Hermite methods.  
\end{abstract}

\maketitle
%\begin{abstract}
%\textcolor{red}{Todo: abstract}.
%\end{abstract}

\section{Introduction}

%For the time-explicit solution of hyperbolic systems of equations, high order finite element (FEM) and discontinuous Galerkin (DG) methods have the advantage of both rapid convergence and decreased numerical dissipation \cite{karniadakis1999spectral, deville2002high, hesthaven2007nodal} compared to low order methods, especially for problems in high frequency wave propagation \cite{melenk2011wavenumber}.  Additionally, high order methods yield a computational structure well-suited to modern computing architectures \cite{klockner2009nodal, markidis2015openacc, remacle2015gpu}.  However, for typical hyperbolic problems, the stable timestep for high order FEM and DG methods decreases as 
%\[
%dt < C \frac{h}{N^2}
%\]
%where $C$ is some constant depending on the problem, $h$ is the mesh size and $N$ is the order of approximation.  This $O(h/N^2)$ behavior of the timestep restriction limits the efficiency of time-explicit methods at higher orders of approximation \cite{klockner2009nodal, modave2015nodal}, though approaches have been proposed which decrease this dependence to $O(h/N)$  \cite{warburton2008taming}.  

The computational simulation of wave propagation is central to geophysical applications, such as seismic imaging and exploration, the modeling of seismic waves induced by earthquakes, and problems in structural acoustics.  However, the numerical modeling of intermediate frequency waves is known to be challenging for many standard low-order methods, requiring a large number of points per wavelength to adequately resolve oscillatory behavior.  Additionally, the simulation of propagating waves using low order methods is typically subject to significant non-physical (numerical) dissipation and dispersion.  High order methods have the advantage of both rapid convergence and decreased numerical dissipation \cite{deville2002high, hesthaven2007nodal} compared to low order methods, especially for problems in intermediate frequency wave propagation \cite{fornberg1987pseudospectral, melenk2011wavenumber}.   High order methods also tend to have a high number of operations per data access, yielding a computational structure well-suited to modern computing architectures \cite{klockner2009nodal, medina2014high, markidis2015openacc}.  

Hermite methods, as introduced by Goodrich et al.\ in \cite{goodrich2006hermite}, are high order methods for wave propagation which represent the solution using a piecewise polynomial basis by collocating the solution and its derivatives on a structured grid.  Solution and derivative information at grid nodes is then used to reconstruct and evolve the solution in time on a staggered (dual) grid.  Hermite methods are provably stable and high order accurate for hyperbolic equations, including problems with varying coefficients.  

Furthermore, though the reconstruction step requires the access of non-local data at neighbor nodes, the computation of derivatives then depends only locally on the reconstructed data at each node.  This is advantageous for high order or multi-stage timestepping methods compared to finite difference methods, where neighboring data must be accessed each time derivatives are approximated.  This structure has also been noted to be well-suited for parallel implementations on modern architectures \cite{appelo2011hermite, dye2015performance}.  Hermite schemes, which were initially introduced for Cartesian domains, have also been coupled with discontinuous Galerkin (DG) schemes for numerical simulations on complex geometries \cite{chen2014hybrid}.  They have also been applied to problems in aeroacoustics \cite{appelo2011hermite}, electromagnetics \cite{chen2014hybrid}, and fluid dynamics \cite{hagstrom2012simulation}.  

%By employing a reconstruction process on a staggered grid, such methods result in a timestep restriction which is independent of the order of approximation, depending only on physical parameters and grid spacing.  These reconstructions require a constant stencil of $2^d$ nodes in $d$-dimensions with increasing order.  Since each node contains high order data, the total amount of data accessed is not decreased compared to equivalent finite difference methods; however, the compact stencil results in a more structured communication pattern than high order finite differences.  Since the stencil is only used in the reconstruction operation, the approximation of derivatives does not require non-local data.  

The Hermite schemes of Goodrich et al.\ \cite{goodrich2006hermite} are one instance of a broader family of methods involving collocation of the solution and its dervatives.  Other methods in this family include shape-preserving methods \cite{fischer1978convective, williamson1989two, rasch1990shape} and jet schemes \cite{nave2010gradient, seibold2011jet, chidyagwai2011comparative}, which use Hermite interpolation in conjunction with semi-Lagrangian techniques to solve advective problems.  These differ from the Hermite schemes discussed here in terms of the characteristic time evolution procedure; however, the analysis and stability of both Hermite and jet schemes both rely primarily on properties of Hermite interpolation under high order Sobolev seminorms.  

Sections~\ref{sec:time} and \ref{sec:framework} present a generalized view of Hermite methods, and motivate new one-step Hermite schemes based on variations in the reconstruction procedure.  These procedures also aim to simplify the implementation and coupling of Hermite and DG schemes \cite{chen2014hybrid}.  Section~\ref{sec:num1D} presents numerical experiments which confirm the high order convergence and stability of each method for the advection equation in one dimension.  Section~\ref{sec:2D} extends each method to two space dimensions and includes numerical results for the two-dimensional advection and acoustic wave equations.  

\subsection{Time evolution}
\label{sec:time}

In this section, we introduce one-dimensional Hermite schemes for the approximation of an evolving solution and its derivatives at a collection of points over an interval $[a,b)\subset \mbb{R}$.  Each Hermite scheme presented has a timestep restriction based only on the domain of dependence for hyperbolic partial differential equations.  For simplicity of presentation, we illustrate this using the 1D periodic scalar advection equation 
\begin{align*}
\pd{u}{t}{} &= c\pd{u}{x}{}\\
u\left(x,t_0\right) &= \mb{u}_0(x)\\
u(a,t) &= u(b,t),
\end{align*}
where  $c$ is a constant advection speed and ${u}_0(x)$ is a smooth initial condition.  Hermite methods may be extended in a straightforward manner to non-uniform grids and more general systems of equations with variable coefficients \cite{goodrich2006hermite}, though these details are omitted for brevity.  

We define first a primary grid $\Omega$ as a collection of $K$ equispaced points 
\[
\Omega = \LRc{x_m: \quad x_m = a + m h_x, \quad m = 0,\ldots, K-1 },
\]
where $K$ is the number of grid points on the interval $(a,b)$ and $h_x = (b-a)/K$ denotes the spacing between the nodes.  For periodic domains, we assume that $x_{m+K} = x_m$.  

Next, we introduce the interpolation length scale $h$ (distinct from the grid spacing $h_x$).  We assume that a smooth function $u(x)$ is well-approximated over some interval $(L_m,R_m)$ with size $h = R_m - L_m$ by a degree 
%$2N+1$ 
$\tilde{N}$
expansion around some point $\tilde{x}_m$.  This expansion takes the form
\[
%u(x) \approx u_m(x) = \sum_{j=0}^{2N+1} \mb{u}_j \LRp{\frac{x-x_m}{h_x}}^j,
u(x) \approx \tilde{u}_m(x) = \sum_{j=0}^{\tilde{N}} \tilde{\mb{u}}_j \LRp{\frac{x-\tilde{x}_m}{h_x}}^j,
\]
where $h_x$ is some spatial length scale, and is typically taken to be some grid spacing in practice.  The vector $\mb{u}$ contains Hermite degrees of freedom, which are scaled spatial derivatives at $\tilde{x}_m$
\[
%\mb{u}_{j} = \LRlim{\frac{h_x^j}{j!}\frac{d^ju_m(x)}{dx^j}}_{x_m},\quad j = 0, \ldots, {2N+1}.
\tilde{\mb{u}}_{j} = \LRlim{\frac{h_x^j}{j!}\frac{d^j\tilde{u}_m(x)}{dx^j}}_{\tilde{x}_m},\quad j = 0, \ldots, \tilde{N}.
\]
For convenience, we express the advection operator applied to $u_m(x)$ as an expansion around $x_m$
\[
%a\pd{u_m}{x} = \sum_{j=0}^{2N+1} \mb{w}_j \LRp{\frac{x-x_m}{h_x}}^j, \qquad {\mb{w}}_j = \LRlim{\frac{h_x^j}{j!}\frac{\partial^j }{\partial x^j}\LRp{a\pd{u_m}{x}}}_{x_m}.
c\pd{\tilde{u}_m(x)}{x} = \sum_{j=0}^{\tilde{N}} {\mb{w}}_j \LRp{\frac{x-\tilde{x}_m}{h_x}}^j, \qquad {\mb{w}}_j = \LRlim{\frac{h_x^j}{j!}\frac{\partial^j }{\partial x^j}\LRp{a\pd{\tilde{u}_m(x)}{x}}}_{\tilde{x}_m}.
\]
The coefficients ${\mb{w}}_j$ are related to ${\mb{u}}_j$ through the derivative matrix $\mb{D}$
\[
%\mb{q}_i = \sum_{j=0}^{2N+1}a\mb{D}_{ij}\mb{u}_j, 
{\mb{w}}=c\mb{D}\mb{u}, 
\qquad 
\mb{D}_{ij} = \begin{cases} 
\frac{i+1}{h_x}&, \quad j = i+1\\
0&, \quad \text{otherwise},
\end{cases}
\qquad 0 \leq i,j \leq \tilde{N}.
\]
This yields a semi-discrete system for the degrees of freedom $\mb{u}_j$
\[
\td{{\tilde{\mb{u}}(t)}}{t} = -c \mb{D} {\tilde{\mb{u}}} (t).
\]

We approximate $u(x,t_n + dt)$ for some timestep $dt > 0$ by solving this semi-discrete system.  This is achieved in \cite{goodrich2006hermite} using a temporal Taylor series: assuming an expansion centered around time $t_n$, 
\[
%u_m(x,t) = \sum_{j=0}^{2N+1} \sum_{k=0}^{{2N+1}-j} \mb{U}_{jk} \LRp{\frac{x-x_m}{h_x}}^j\LRp{\frac{t-t_n}{dt}}^k,
\tilde{u}_m(x,t) = \sum_{j=0}^{\tilde{N}} \sum_{k=0}^{\tilde{N}-j} \tilde{\mb{U}}_{jk} \LRp{\frac{x-\tilde{x}_m}{h_x}}^j\LRp{\frac{t-t_n}{dt}}^k,
\]
and using the Cauchy-Kowalevski relation, the coefficients $\mb{U}_{jk}$ may be shown to satisfy
%and differentiating the advection equation $r$ times in $x$ and $s$ times in $t$ yields 
%\begin{equation}
%\pdp{}{x}{r}\pdp{u}{t}{s} = a\pdp{}{x}{r+1}\pdp{u}{t}{s-1}.
%\label{eq:ck}
%\end{equation}
%Plugging $u_m(x,t)$ into the above and evaluating at $(x_m, t_n)$, we have that 
%\[
%\mb{U}_{j0} = \mb{u}_j, \qquad \LRlim{\pdp{}{x}{j}\pdp{u}{t}{k}}_{x_m,t_n} = \mb{U}_{jk}\frac{j!k!}{h^j dt^k}, \qquad 0 < j \leq {2N+1}-k.
%\]
%Combining the above with (\ref{eq:ck}) implies that the coefficients $\mb{U}_{jk}$ may be expressed in terms of known coefficients, 
\[
%\mb{U}_{jk} = \frac{1}{k dt} a\frac{(j+1)}{ h_x } \mb{U}_{j+1,k-1} \qquad j = 0, \ldots , {2N+1} - k.
\tilde{\mb{U}}_{jk} = \frac{c}{k dt} \frac{(j+1)}{ h_x } \tilde{\mb{U}}_{j+1,k-1} \qquad j = 0, \ldots , \tilde{N} - k.
\]
or more succinctly using matrix-vector notation
\[
\tilde{\mb{U}}_{(\cdot, k)} = \frac{c}{k dt}  \mb{D}\tilde{\mb{U}}_{(\cdot, k-1)}.
\]
where $\tilde{\mb{U}}_{(\cdot,k)}$ refers to the $k$th column of the array $\tilde{\mb{U}}$.  Inserting these coefficients into the temporal Taylor expansion and evaluating at time $t_n + dt$ yields an update for the solution $u_m(x,t_n)$.  In practice, this is computed using Algorithm~\ref{alg:time} \cite[p109]{canuto2006spectral}, which may be generalized to linear autonomous systems in multiple dimensions.  Alternatively, more standard time integration techniques may be used to advance the solution forward in time.  
\begin{algorithm}
\begin{algorithmic}[1]
\Procedure{Temporal Taylor series evaluation }{}
\State $\mb{w} = \tilde{\mb{u}}^n$
%\For{$\ell = 2N+1,\ldots, 2$}
\For{$\ell = \tilde{N},\tilde{N}-1,\ldots, 0$}
\State $\mb{w} = \mb{w} + \frac{dt}{1 + \ell} (-c \mb{D}){\tilde{\mb{u}}^n }$
\EndFor
\State $\tilde{\mb{u}}^{n+1} =  {\mb{w}}$
\EndProcedure
\end{algorithmic}
\caption{Time evolution procedure for 1D scalar advection.}
\label{alg:time}
\end{algorithm}
This evolution of the solution from $t_n$ to $t_{n+1} = t_n + dt$ may be represented by the application of some update matrix $\mb{T}^{dt}$ to the degrees of freedom at time $t_n$, such that
\[
%\mb{u}^{n+1} = \mb{T}^{dt} \mb{u}^n, \qquad u_m(x,t_n) = \sum_{j=0}^{2N+1}\mb{u}^n_j \LRp{\frac{x-x_m}{h_x}}^j.
\tilde{\mb{u}}^{n+1} = \mb{T}^{dt} \tilde{\mb{u}}^n, \qquad \tilde{u}_m(x,t_n) = \sum_{j=0}^{\tilde{N}}\tilde{\mb{u}}^n_j \LRp{\frac{x-\tilde{x}_m}{h_x}}^j.
\]

%\begin{equation}
%\mb{U}_{jk} =\frac{1}{kdt} \frac{a(j+1)}{ h_x } \mb{U}_{j+1,k-1} \quad j = 0, \ldots , {2N+1} - k.
%\end{equation}

%Iterating this gives
%\begin{equation}
%\mb{U}_{jk} = \LRp{\frac{a}{h dt}}^k \mb{U}_{j+k} \frac{(j+1)\ldots(j+k)}{k!}.
%\label{eq:cku}
%\end{equation}
%Suppose that at some time $t_n$, the spatial expansion of $u(x,t_n)$ around $x_m$ is represented by coefficients $\mb{u}_j^n$, and that for $t_{n+1}= t_n + dt$, $u(x, t_{n+1})$ is represented by coefficients $\mb{u}_j^{n+1}$.  
%Then, the Taylor series may be evolved in time through evaluation of the Taylor series and its derivatives at $x = x_m$ and $t_n + dt$ and use of (\ref{eq:cku}) yields the update formula
%\[
%\mb{u}_{i}^{n+1} = \sum_{j=0}^{2N+1} \mb{T}^{dt}_{ij} \mb{u}_j^n, \qquad 
%\mb{T}^{dt}_{ij} = \begin{cases}
%0, &j < i\\
%\LRp{\frac{a}{h_xdt}}^{j-i} \frac{1}{(j-i)!}\prod\limits_{s=0}^{i-1}(j-s), & j \geq i.  
%\end{cases}
%\]
%The matrix $\mb{T}^{dt}$ evolves coefficients $\mb{u}^n$ at timestep $t_n$ forward in time by $dt$.   The exact form of the matrix $\mb{T}^{dt}$ varies between problems and discretizations, though for linear problems, it is always possible to construct such a matrix or compute its the action on a vector using the temporal Taylor series and (\ref{eq:ckc}).

As noted in \cite{goodrich2006hermite}, for the scalar advection equation, the Taylor expansion in time yields an exact evolution of the approximation to $\mb{u}_0$, so long as the domain of dependence at $t_0 + dt$ lies within the interval $(L_m ,R_m)$.  This translates into a degree-independent timestep restriction 
\[
cdt < \frac{h}{2}, \qquad h = \min\{{\tilde{x}_m-L_m}, {R_m - \tilde{x}_m}\}.
\]
If $L_m$ and $R_m$ are chosen symmetrically around $x_m$, $h$ the interpolation interval reduces to $(L_m, R_m) = (x_m-h,x_m + h)$.  

The same domain of dependence argument motivates timestep restrictions for ``tent-pitching'' space-time finite element methods, which use a coupled discretization in both time and space.  The stability of the space-time formulation results in a similar causal timestep restriction \cite{falk1999explicit}.  

%Alternatively, time-varying coefficients may be assumed 
%\[
%u(x,t) = \sum_{j=0}^{2N+1} \mb{u}_j(t) \LRp{\frac{x-x_m}{h_x}}^j.
%\]
%Substituting this expression into the original differential equation, taking spatial derivatives of both sides and evaluating at $x_m$ yields ordinary differential equations for each coefficient
%\[
%\td{{\mb{u}_j}(t)}{t} = a\mb{u}_{j+1}(t), \quad j = 0,\ldots, {2N}, \qquad \td{{\mb{u}_{2N+1}}(t)}{t} = 0,  
%\]
%which may be solved to yield a solution at some time $t = t_n + dt$.  

\subsection{Local Hermite interpolation}

The above evolution procedure hinges on a degree $\tilde{N}$ polynomial representation of some smooth function $u(x)$ in an interval $(L_m,R_m)$ which accurately approximates the solution and its derivatives at the point $\tilde{x}_m$.  This is addressed in \cite{goodrich2006hermite} using degree $N$ Hermite interpolation to produce a degree $\tilde{N} = 2N+1$ reconstruction at specific points $\tilde{x}_m$.  

The degree $N$ Hermite interpolant of $u(x)$ over $(L_m, R_m)$ (which we refer to as $\tilde{u}_m(x)$) is constructed by specifying $(N+1)$ scaled spatial derivatives at the left and right endpoints of the interval $(L_m,R_m)$
\[
\mb{u}^{L}_{j} = \left.\frac{h_x^j}{j!}\frac{d^ju}{dx^j}\right|_{L_m}, \qquad \mb{u}^R_{j} = \left.\frac{h_x^j}{j!}\frac{d^ju}{dx^j}\right|_{R_m}, \qquad j = 0, \ldots, N.
\]
The resulting expansions around $L_m, R_m$%$x_{m}-h$ and $x_{m}+h$
\begin{align*}
u^L_{m}(x) = \sum_{j=0}^N \mb{u}^L_j\LRp{\frac{ x - L_m}{h_x}}^j, \qquad u^R_{m}(x) = \sum_{j=0}^N \mb{u}^R_j\LRp{\frac{ x - R_m}{h_x}}^j
\end{align*}
coincide with the value and first $N$ derivatives of $u(x)$ at each point.  

The Hermite interpolant $\tilde{u}_m(x)$ is a polynomial of order $2N+1$ over $(L_m, R_m)$, which is defined by interpolating the $(N+1)$ solution and derivative values at each endpoint.  $\tilde{u}_m(x)$ is represented using the following expansion
\[
\tilde{u}_m(x) = \sum_{j=0}^{2N+1} \tilde{\mb{u}}_j\LRp{\frac{ x - \tilde{x}_{m}}{h_x}}^j
\]
where $\tilde{\mb{u}}_j$ are values of the solution and $2N+1$ scaled derivative values at $\tilde{x}_m$.  These coefficients may be determined by solving the interpolation problem
\[
\left.\frac{\partial^i \tilde{u}_m}{\partial x^i}\right|_{L_m} = \left.\frac{\partial^i u_m^L}{\partial x^i}\right|_{L_m}, \qquad 
\left.\frac{\partial^i \tilde{u}_m}{\partial x^i}\right|_{R_m}= \left.\frac{\partial^i u_m^R}{\partial x^i}\right|_{R_m}, \qquad i = 0,\ldots, 2N+1.
\]
This results in the system
\begin{equation}
\LRs{\begin{array}{c}
\mb{C}^L \\
\mb{C}^R
\end{array}}\tilde{\mb{u}} = \LRs{\begin{array}{c}
\mb{u}^L \\
\mb{u}^R
\end{array}}, 
\label{eq:interp}
\end{equation}
where the constraint matrices $\mb{C}^L, \mb{C}^R \in \mbb{R}^{(N+1)\times(2N+2)}$ enforce conditions on $u_m(x)$ and its $2N+1$ derivatives at the points $L_m, R_m$%$x_m-h$ and $x_m+h$
\begin{align*}
\mb{C}^L_{mn} &= 
\begin{cases}
\LRp{\frac{L_m-\tilde{x}_m}{h_x}}^{n-m}\frac{1}{m!}\prod\limits_{s=0}^{m-1}(n-s), & n\geq m\\
0, &n < m
\end{cases}\\
\mb{C}^R_{mn} &= 
\begin{cases}
\LRp{\frac{R_m-\tilde{x}_m}{h_x}}^{n-m}\frac{1}{m!}\prod\limits_{s=0}^{m-1}(n-s), & n\geq m\\
0, &n < m.
\end{cases}
\end{align*}
%In addition to the order of approximation, both $\mb{C}^L$ and $\mb{C}^R$ depend only the ratio ${h}/{h_x}$ of the interpolation interval width to the spatial length scale, such that specifying $h$ is enough to determine these matrices in 1D.  
For convenience, we define $\mb{H}$
\[
\mb{H} = \LRs{\begin{array}{c}
\mb{C}^L \\
\mb{C}^R
\end{array}}^{-1}
\]
as the interpolation matrix which maps solution and scaled derivative data at the endpoints of the interval $(L_m, R_m)$ to a degree $2N+1$ expansion at the node $\tilde{x}_m$.  Variations in the construction of $\mb{H}$ result in different interpolation/reconstruction schemes.  

For the remainder of this work, we refer to the expansion of $\tilde{u}_m(x)$ as a Hermite reconstruction at $\tilde{x}_m$.  This expansion may then be evolved in time a distance of $dt$ using the evolution matrix $\mb{T}^{dt}$.  %The interpolation points $(L_m, R_m)$ may be distinct from $\tilde{x}_m$, the point at which the solution is evolved forward in time.  

\section{Hermite methods in one dimension}
\label{sec:framework}
Having introduced evolution and interpolation procedures, one-dimensional Hermite methods may now be specified.  Given some grid $\Omega$, Hermite approximation spaces are associated with each grid.  The degree $N$ space $\mc{P}^N_{\Omega}$ is defined to be
\[
\mc{P}^N_{\Omega} = \LRc{\bigoplus_{m=0}^{K-1} \mc{P}^{2N+1}[x_m,x_{m+1})} \cap C^N[a,b].  
\]
which consists of piecewise polynomials of degree $2N+1$ with $N$ globally continuous derivatives at each node $x_m \in \Omega$.  Each $u(x)\in \mc{P}^N_{\Omega}$ is defined by $(N+1)$ pieces of Hermite interpolation data at each node $x_m$
\[
\mb{U}_{mj} = \frac{h_x^j}{j!}\left.\frac{\partial^j u}{\partial x^j}\right|_{x_m}, \qquad j = 0,\ldots, N+1, \qquad m = 0,\ldots,K-1.
\]
For brevity, we refer to $\mb{U}_{(m,\cdot)}$, the vector of derivatives at the node $x_m$, as $\mb{U}_{m}$.  

Finally, in addition to the evolution and interpolation operators $\mb{T}^{dt}$ and $\mb{H}$, respectively, we introduce the restriction operator $\mb{R} \in \mbb{R}^{(N+1)\times(2N+1)}$ such that $\mb{R}_{ij} = \delta_{ij}$ and multiplication of a vector by $\mb{R}$ extracts the first $N+1$ entries of that vector.  

Hermite methods march forward in time by combining three procedures over one or more stages:
\begin{enumerate}
\item \textbf{Interpolation} using points in the primary grid to produce Hermite reconstructions of higher degree $\tilde{N}$ centered around points $\tilde{x}_m$.  
\item \textbf{Evolution} of the higher degree Hermite reconstructions at points $\tilde{x}_m$ forward in time to $\tilde{t} = t_n + dt$.  
\item \textbf{Restriction} of the solution by truncating the degree of the polynomial expansion at $\tilde{x}_m$.  
\end{enumerate}
We emphasize that for each Hermite method, the stable timestep restriction is independent of the degree of approximation $N$.  This is due to the fact that, by evolving the solution in time using a temporal Taylor series, a single update step of a Hermite method may be interpreted as the composition of the \emph{exact} evolution of piecewise polynomial data with a \emph{projection} in a seminorm which is preserved by the solution \cite{goodrich2006hermite}.  In particular, Hermite interpolation in one space dimension results in the projection onto piecewise degree-$2N+1$ polynomial in the $H^{N+1}$ seminorm.  Thus, as described in Section~\ref{sec:time}, the timestep restriction is determined only by the domain of dependence of the equation and the interval of the Hermite reconstruction.  

We refer to the collection of of reconstruction points $\tilde{x}_m$ as an auxiliary grid $\tilde{\Omega}$ on which the solution is evolved in time.  Different schemes use differing combinations of interpolation and evolution procedures, which are summarized in Table~\ref{table:diagram}.  Different Hermite schemes may also vary parameters of the interpolation process.  For example, the Hermite schemes of Goodrich et al.\ \cite{goodrich2006hermite} (referred to henceforth in this paper as Dual Hermite schemes) produce a Hermite reconstruction centered between two grid points, which are represented over an auxiliary grid consisting of midpoints of the primary grid.  These reconstructions are then evolved forward in time and truncated.  The resulting auxiliary grid data may then be used to compute Hermite reconstructions at the original primal points, which are then evolved in time and truncated to complete a single timestep.  

We introduce here the Virtual Hermite method, which is equivalent to the Dual Hermite method with a timestep of size zero on the auxiliary grid.  As a result, the interpolations to and from the auxiliary grid may be combined, resulting in a step which interpolates and reconstructs on the same primary grid.  

The Central and Upwind Hermite methods aim to interpolate and reconstruct on the same primary grid through redefinitions of the interpolation operator.  The Central scheme expands the interval of interpolation, producing a Hermite reconstruction at a point $x_m$ information at neighboring nodes $x_{m-1}, x_{m+1}$, while the Upwind Hermite scheme uses data from a single neighbor to produce a directional reconstruction.  

\begin{table}[!h]
\centering                                                          
\begin{tabular}{|c||c|c|c|c|@{}m{0pt}@{}}
\hline
& \multicolumn{2}{c|}{Stage 1} &  \multicolumn{2}{c|}{Stage 2}  \\ 
\hline
 & Interpolate & Evolve & Interpolate & Evolve \\ 
\hhline{|=||=|=|=|=|}
Dual & ${\mb{U}^{n}_{{m}},\; \mb{U}^{n}_{{m+1}}}\rightarrow \tilde{\mb{U}}^n_{{m+1/2}}$ & $[t_n, t_n + dt)$ & $\mb{U}^{n+1}_{{m-1/2}},\; {\mb{U}}^{n}_{m+1/2}\rightarrow \tilde{\mb{U}}^n_{{m}}$ & $[t_n, t_n + dt)$ &\\[10pt]
\hline                                                              
Virtual & ${\mb{U}^{n}_{{m}},\; \mb{U}^{n}_{{m+1}}}\rightarrow \tilde{\mb{U}}^n_{{m+1/2}}$ & $[t_n, t_n)$ &  $\mb{U}^{n+1}_{{m-1/2}},\; {\mb{U}}^{n}_{m+1/2}\rightarrow \tilde{\mb{U}}^n_{{m}}$ & $[t_n, t_n + dt)$&\\[10pt]
\hline                                                              
Central & ${\mb{U}^{n}_{{m-1}},\; \mb{U}^{n}_{{m+1}}}\rightarrow \tilde{\mb{U}}^n_{{m}}$ & $[t_n, t_n + dt)$ & &&\\[10pt]
\hline                                                              
Upwind &${\mb{U}^{n}_{{m-1}},\; \mb{U}^{n}_{{m}}}\rightarrow \tilde{\mb{U}}^n_{{m}}$& $[t_n, t_n + dt)$ & &&\\[10pt]
\hline                                                              
\end{tabular}                                                       
\caption{Overview of interpolation and evolution operations over different stages for various Hermite method, where $\mb{U}^n_{m}$ refers to the Hermite solution at the $n$th timestep and $m$th grid point.  Note that the stable time step restriction results in a different $dt$ for each method.}
\label{table:diagram}
\end{table}
%\textcolor{red}{Discuss general structure of each Hermite scheme.  Introduce secondary grid.  }

%Hermite schemes rely on the use of interpolation operators to transfer information between grids.  For two grids $\Omega_1$ and $\Omega_2$, an interpolation procedure is defined by specifying left and right interpolation points $L_i, R_i \in \Omega_1$ which define the interpolation interval for each node $x_i \in \Omega_2$.  
%At each point $x_i$, a Hermite expansion is constructed using data from $L_i, R_i$
%\[
%\mb{u}^C = \mb{H}
%\LRs{\begin{array}{c}
%\mb{U}_{(L_i,\cdot)}\\
%\mb{U}_{(R_i,\cdot)}
%\end{array}},
%\]
%where $\mb{U}_{i,\cdot}$ refers to the $i$th column of the array $\mb{U}$.  
%
%$L_i, R_i$ are typically specified by the specific Hermite method.  By the domain of dependence arguments in Section~\ref{sec:time}, the length of the interpolation interval $h$ will also determine the maximum timestep ${dt}$ for the evolution operator $\mc{T}_{dt}$.
%
%\textcolor{red}{CONVERT EVERYTHING TO MATRICES.}
%
%Each Hermite method of degree $N$ begins with a primary grid $\Omega$, and may also utilize one or more auxiliary grids $\Omega_A$.  The primary grid $\Omega$ contains equispaced points
%\[
%\Omega = \LRc{x_i: \quad x_i = a + i h_x, \quad i = 0,\ldots, K },
%\]
%where $K$ is the number of grid points on the bi-unit interval and $h_x = (b-a)/K$ denotes the spacing between the nodes (which is distinct from $h$).  For a periodic domain, we assume the grid is also periodic, such that $x_{i+K} = x_i$.  

While the implementation of the time evolution operator changes from problem to problem, the interpolation procedure remains the same irregardless of the equation solved.  For this reason, we focus on the description of the interpolation procedure for three specific Hermite methods in the following Sections.  Common to each Hermite scheme is a timestep restriction which is independent of the polynomial degree $N$.   

For each one-dimensional Hermite method discussed, the timestep restriction is presented only for the scalar advection equation.  An extension to $k\times k$ systems of hyperbolic equations 
\[
\pd{u}{t} = A\pd{u}{x},\qquad u\in \mbb{R}^k, \qquad A\in \mbb{R}^{k\times k}
\]
results in timestep restrictions which are identical to the scalar advection equation, except that the wavespeed $c$ is replaced by the spectral radius $\rho(A)$.  A similar analysis extends these results to variable coefficient problems, and both are described in more detail in \cite{goodrich2006hermite}.  

In the following sections, we introduce the Dual, Virtual, Central, and Upwind Hermite methods in more detail.  

\subsection{Dual Hermite method}
\label{sec:dual}

Hermite methods were originally introduced by Goodrich, Hagstrom, and Lorenz in \cite{goodrich2006hermite}, using ``primal'' and ``dual'' grids to facilitate Hermite reconstructions at specific points.  We refer to this specific Hermite method as the \textit{Dual} Hermite method.  The Dual Hermite method introduces the auxiliary grid $\tilde{\Omega}$, which is taken to be a dual or co-volume grid 
\[
\tilde{\Omega} = \LRc{\tilde{x}_{m+1/2} = a + \LRp{m + {1}/{2}} h_x, \quad m = 0,\ldots, K-1}
\]
such that the nodes of $\tilde{\Omega}$ are staggered a distance of $h_x/2$ between the nodes of $\Omega$.  For periodic grids, the dual grid also satisfies $\tilde{x}_{m +1/2+ K} = \tilde{x}_{m+1/2}$.  We associate also an approximation space $\mc{P}^N_{\tilde{\Omega}}$ to the dual grid
\[
\mc{P}^N_{\tilde{\Omega}} = \LRc{\bigoplus_{m=1}^{K-1} \mc{P}^{2N+1}[\tilde{x}_{m-1/2}, \tilde{x}_{m+1/2})} \cap C^N[a,b].
\]
Each $q(x)\in \mc{P}^N_{\tilde{\Omega}}$ is defined by degrees of freedom $\mb{Q}_{mj}$
\[
\mb{Q}_{m+1/2,j} = \frac{h_x^j}{j!}\LRlim{\frac{\partial^j q}{\partial x^j}}_{\tilde{x}_{m+1/2}}, \qquad j = 0,\ldots, N+1, \qquad m = 0,\ldots,K-1.
\]
In the Dual Hermite method, the interpolation length scale is $h_x/2$, such that the timestep restriction on both primary and dual grids is 
\[
cdt < h_x/2.
\]  
The utility of the dual grid comes in considering the Hermite reconstruction at $x_m \in \Omega$, which interpolates at the points $L_m = x_m-h_x/2$ and $R_m = x_m + h_x/2$, such that
\[
(L_m, R_m) = \LRp{x_m - \frac{h_x}{2}, x_m + \frac{h_x}{2}} = \LRp{\tilde{x}_{m-1/2}, \tilde{x}_{m+1/2}}.  
\]
In other words, nodal data at points on the primary grid is used to produce a Hermite reconstruction $\tilde{u}_m(x)$ at each point $\tilde{x}_m$ on the dual grid. This also defines the interpolation length scale $h = h_x/2$, implying a time step restriction of $cdt < h_x/2$.  

Denoting the vector of nodal data at a point $x_{m \pm 1/2}$ as $\mb{Q}_{m\pm 1/2}$, the interpolation procedure is 
\[
\tilde{\mb{U}}_{m} = \mb{H}\LRs{\begin{array}{c}
\mb{Q}_{m-1/2}\\ \mb{Q}_{m+1/2}
\end{array}}
\]
where $\tilde{\mb{U}}_{m}$ represents a degree $2N+1$ expansion at a point $x_m$ on the primary grid.  

\begin{figure}[h!]
\centering
\subfloat[Primary to dual grid]{\includegraphics[width=.4\textwidth]{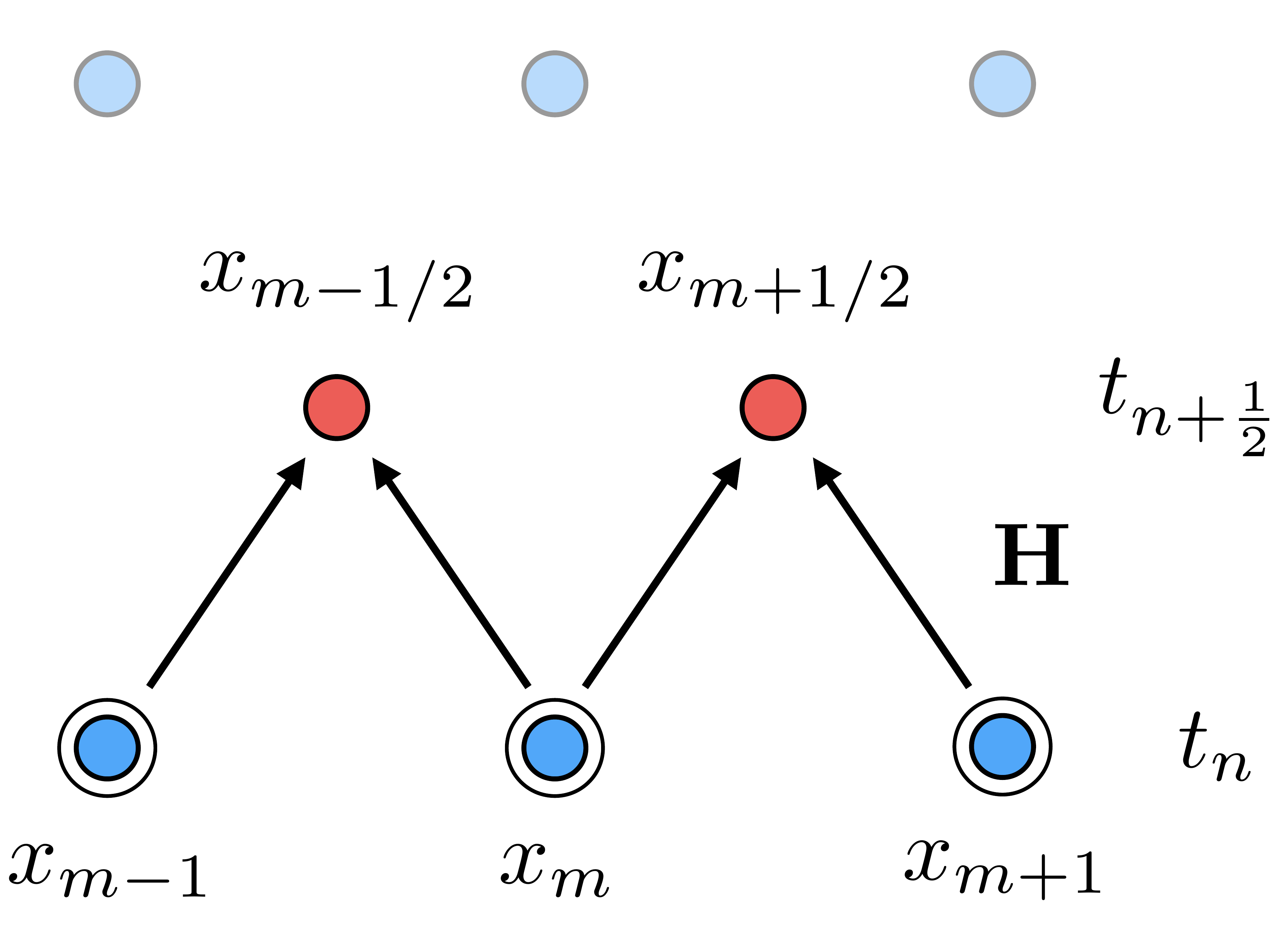}}
\hspace{4em}
\subfloat[Dual to primary grid]{\includegraphics[width=.4\textwidth]{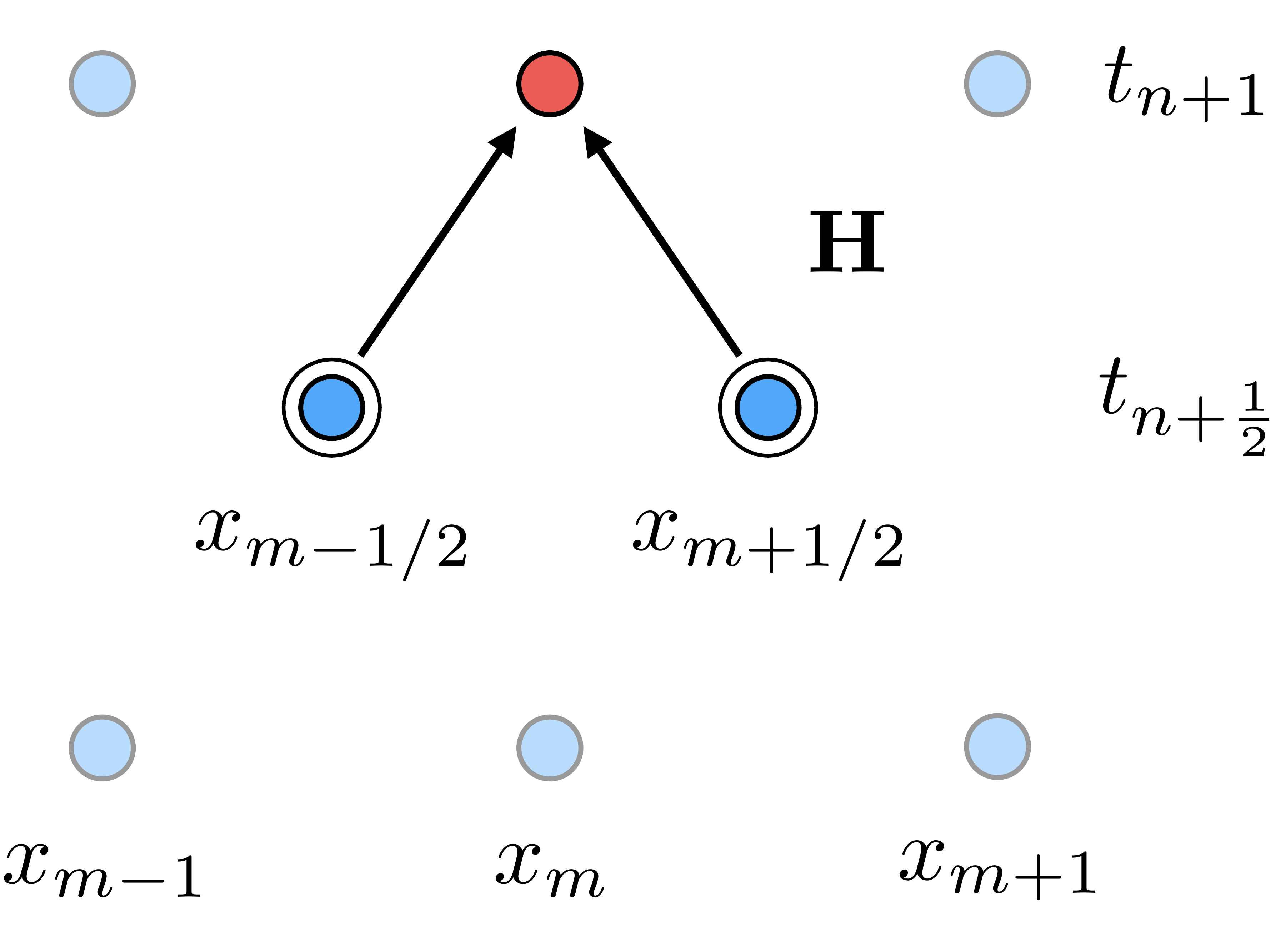}}
\caption{Interpolation procedures to and from primary and dual grids for the Dual Hermite method.  Nodes which contribute data to the reconstruction are circled. }
\label{fig:dgH}
\end{figure}

Suppose $u^n(x) \in P^N_{\Omega}$ is the solution at timestep $n$ with degrees of freedom $\mb{U}^n_{mj}$.  The Dual Hermite method interpolates the primary grid solution at $x_m$ and $x_{m+1}$ to the dual grid and evolves it in time.  The degree $2N+1$ solution is then truncated, producing $v^{n+1/2}_{m+1/2}(x) \in P^{N}_{\tilde{\Omega}}$.  Denoting degrees of freedom for $v^{n+1/2}_{m+1/2}(x)$ as $\mb{Q}^{n+1/2}_{m+1/2}$, this step of a Hermite method may be expressed as 
\begin{align*}
\mb{Q}^{n+1/2}_{m+1/2} = \mb{R}\mb{T}^{dt}\mb{H}\LRs{\begin{array}{c}
\mb{U}^{n}_{m-1}\\ \mb{U}^{n}_{m+1}
\end{array}}, \quad m = 0,\ldots, K-1.
\end{align*}
In order to update the solution on the primary grid, the process is repeated, except that data from the dual grid is transferred to the primary grid before being evolved in time
\begin{align*}
\mb{U}^{n+1}_{m} = \mb{R}\mb{T}^{dt}\mb{H}\LRs{\begin{array}{c}
\mb{Q}^{n+1/2}_{m-1/2}\\ \mb{Q}^{n+1/2}_{m+1/2}
\end{array}}, \quad m = 0,\ldots, K-1.
\end{align*}
The complete process illustrated in Figure~\ref{fig:dgH}.  

\subsection{Virtual Hermite method}
\label{sec:virtual}

The Virtual Hermite method is motivated by the fact that two timesteps in the Dual Hermite method may be collapsed into a single update step on the primary grid involving nodal data at points $x_{m-1}, x_m, x_{m+1}$, bypassing explicit time evolution on the dual grid.  This requires the formation of the time evolution operator explicitly, which unfortunately depends on physical and discretization parameters, and may vary between timesteps for nonlinear problems.  We propose the Virtual Hermite method to avoid the explicit construction of $\mb{T}^{dt}$ for the dual grid.  

The Virtual Hermite method is identical to the Dual Hermite method except for time-evolution on the dual grid. In the Dual Hermite method, the solution on both the primary and dual grids is evolved with timestep $dt$.  The Virtual Hermite method skips one evolution step, taking a timestep of $dt = 0$ on the dual grid.  The evolution operator $\mb{T}^{dt}$ on the dual grid then becomes the identity matrix, and the two steps of the Dual Hermite method may be collapsed into a single update step 
\[
\mb{U}^{n+1}_{m} = \mb{F}
\begin{bmatrix}
\mb{U}^{n}_{m-1}\\
\mb{U}^{n}_{m}\\
\mb{U}^{n}_{m+1}
\end{bmatrix}, \qquad 
\mb{F} = \mb{R}\mb{T}^{dt}\mb{H}
\begin{bmatrix}
\boxed{\,\mb{R}\mb{H}\;\;\;\;} \;\;\;\mb{0}\\
\mb{0}\;\;\; \boxed{\;\;\;\;\mb{R}\mb{H}\,}
\end{bmatrix}
\]
where $\mb{0} \in \mbb{R}^{(N+1)\times (N+1)}$ and $\mb{R}\mb{H}\in \mbb{R}^{(N+1) \times (2N+2)}$.  

The matrix $\mb{F}$ resembles the co-volume filter analyzed in \cite{warburton2008taming}, which projects the solution on a primary grid to and from a staggered dual grid, suppressing spurious gradients of the solution on the primal grid.  The Virtual Hermite method uses a similar procedure, where $\mb{F}$ maps Hermite data of degree $N$ to Hermite data of degree $2N+1$ by transferring to and from the dual grid.  We expect Virtual Hermite solutions to resemble filtered Dual Hermite solutions, which is supported by numerical experiments in Section~\ref{sec:spectra}.  

\begin{figure}[h!]
\centering
\subfloat[Virtual dual nodes]{\includegraphics[width=.4\textwidth]{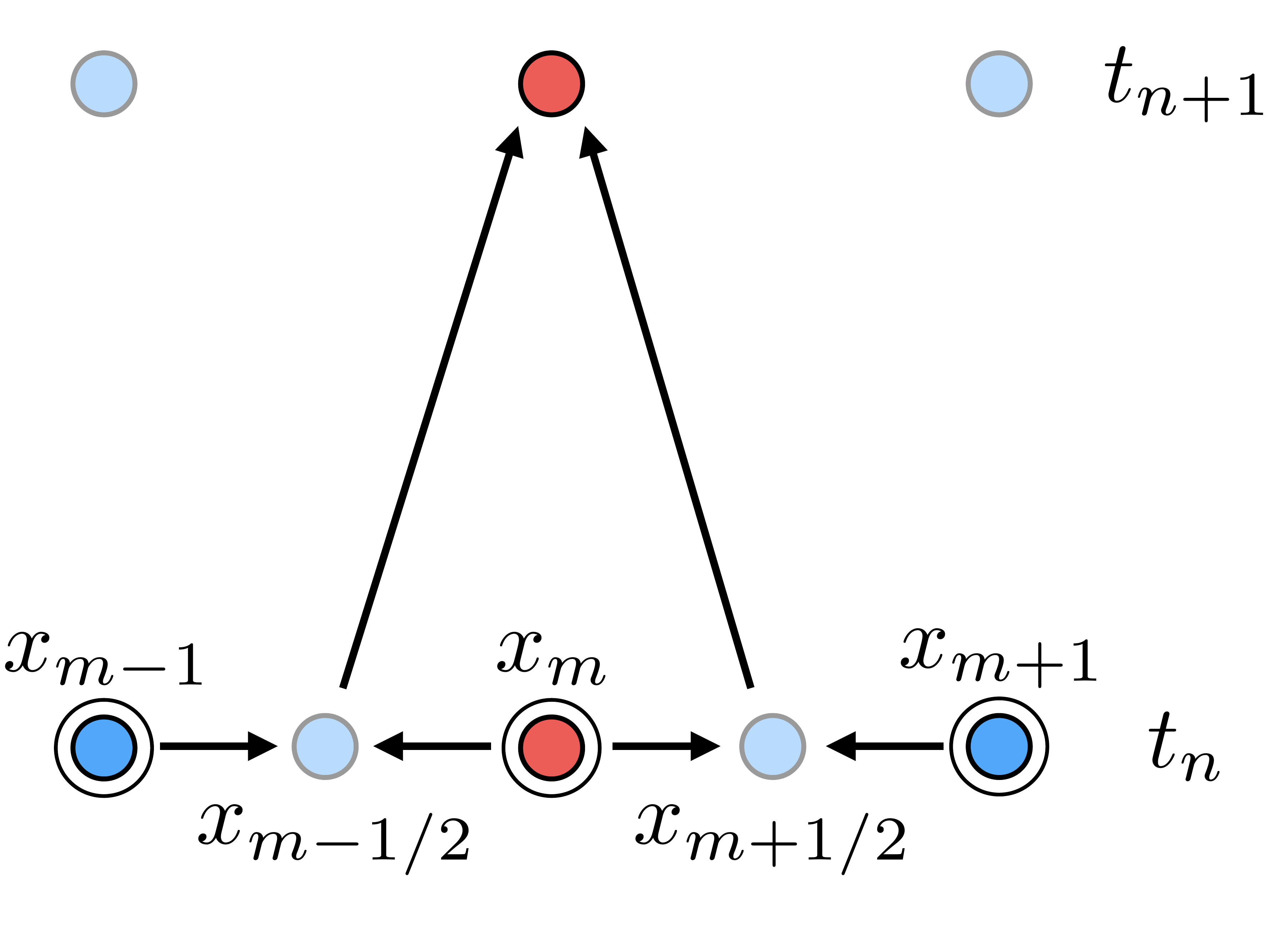}}
\hspace{4em}
\subfloat[Collapsed dual nodes]{\includegraphics[width=.4\textwidth]{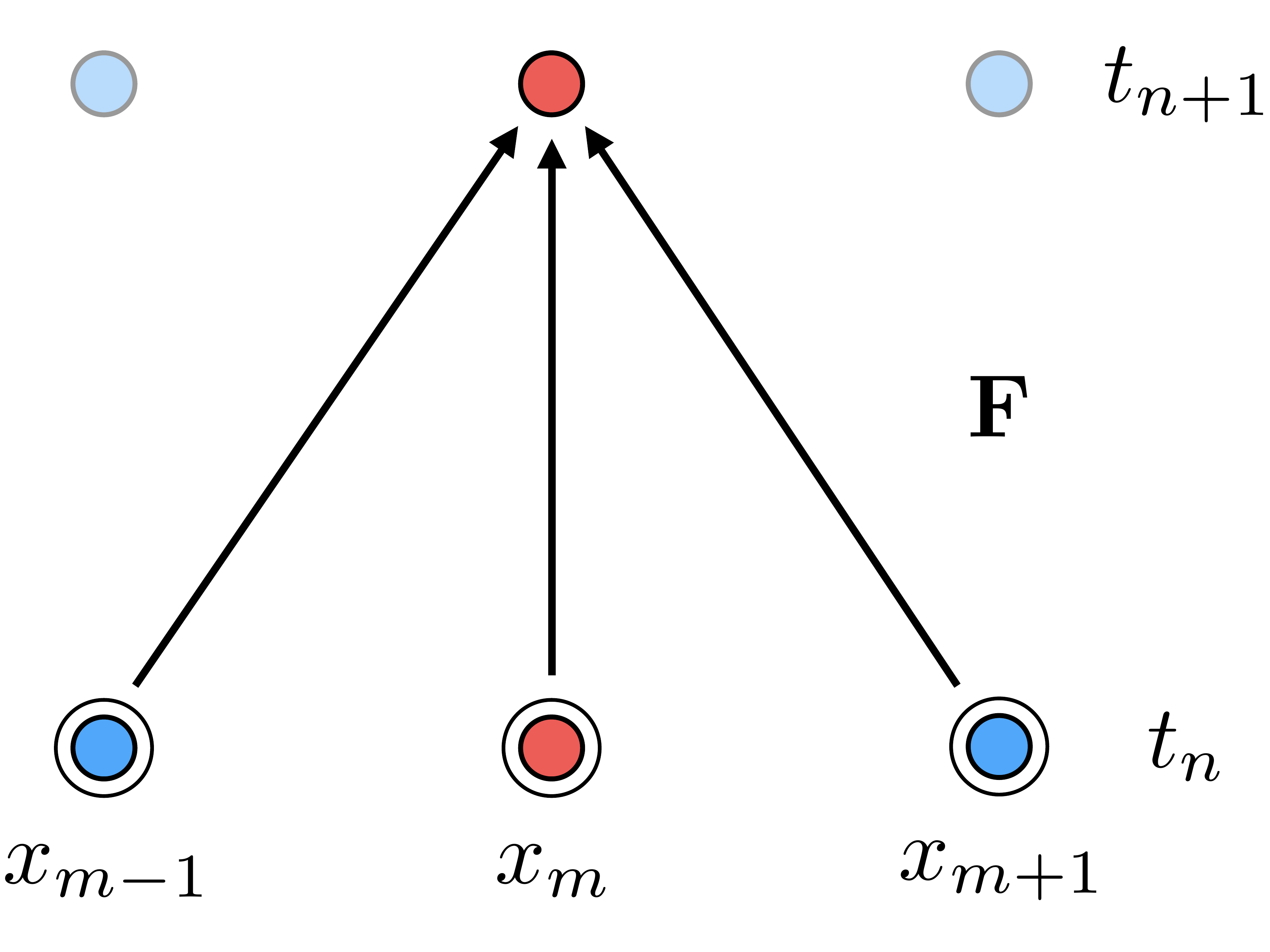}}
\caption{Virtual Hermite interpolation procedure, which transfers to and from an auxiliary grid (left) or as a single reconstruction step involving a three-node stencil (right).  Nodes which contribute data to the reconstruction are circled.}
\label{fig:virtual}
\end{figure}

The Virtual Hermite method obeys the timestep restriction $c dt < h_x/ c$.   This is the same restriction observed for a single step of the Dual Hermite method, since the interpolation length scale $h = h_x/2$ is the same in both cases.  However, by forming $\mb{F}$ and using a single update step on the primary grid, the Virtual Hermite method eliminates the need to explicitly store and evolve dual grid solutions.  We refer to $\mb{F}$ as an operator with a 3-node stencil since degree $N$ data from the three nodes $x_{m-1}, x_m, x_{m+1}$ is required to produce degree $2N+1$ data at $x_m$, as shown in Figure~\ref{fig:virtual}.  %$\mb{F}$ may also be viewed as a map from the primary grid to itself, which reconstructs the solution at each primary grid point prior to time evolution.  %This map from the primary grid to itself motivates the Central Hermite scheme, which defines the interpolation operator in a similar manner.  

We note that it is also possible to produce higher degree reconstructions using the same 3-node stencil.  One approach is to directly using the degree $2N+1$ reconstructions at the two dual grid nodes to produce a final reconstruction of degree $4N+3$ at $x_m$, instead of truncating the dual grid reconstructions.  Another option is to directly interpolating the $3(N+1)$ solution and derivative values at each node to produce a degree $3N+2$ reconstruction at $x_m$.  In both cases, the higher degree reconstruction may then be evolved in time using a higher order scheme; however, we do not observe significantly improved convergence rates under such a procedure, and for some values of $N$, higher degree reconstructions results in an unstable scheme.  Section~\ref{sec:spectra} describes an alternative way to determine a higher degree reconstruction based on optimization of discrete dispersion and dissipation relations \cite{tam1993dispersion}.  

%\textcolor{red}{mention stability problems of each method, and maybe optimization of Hermite stencils.}

\subsection{Central Hermite method}

While the Virtual Hermite method removes the need to update the solution in time on the dual grid, the Central Hermite method sidesteps the use of an dual grid altogether by defining the Hermite reconstruction at $x_m$ through interpolation at neighboring points
\[
(L_m, R_m) = (x_m - h_x, x_m + h_x) = (x_{m-1},x_{m+1}).
\]
As a result, the interpolation length scale is $h = h_x$ and the timestep restriction is 
\[
cdt < h_x.  
\]
%\[
%\LRp{x_m - h_x, x_m + h_x} = \LRp{x_{m-1},x_{m+1}}.  
%\]
\begin{figure}[h!]
\centering
\includegraphics[width=.4\textwidth]{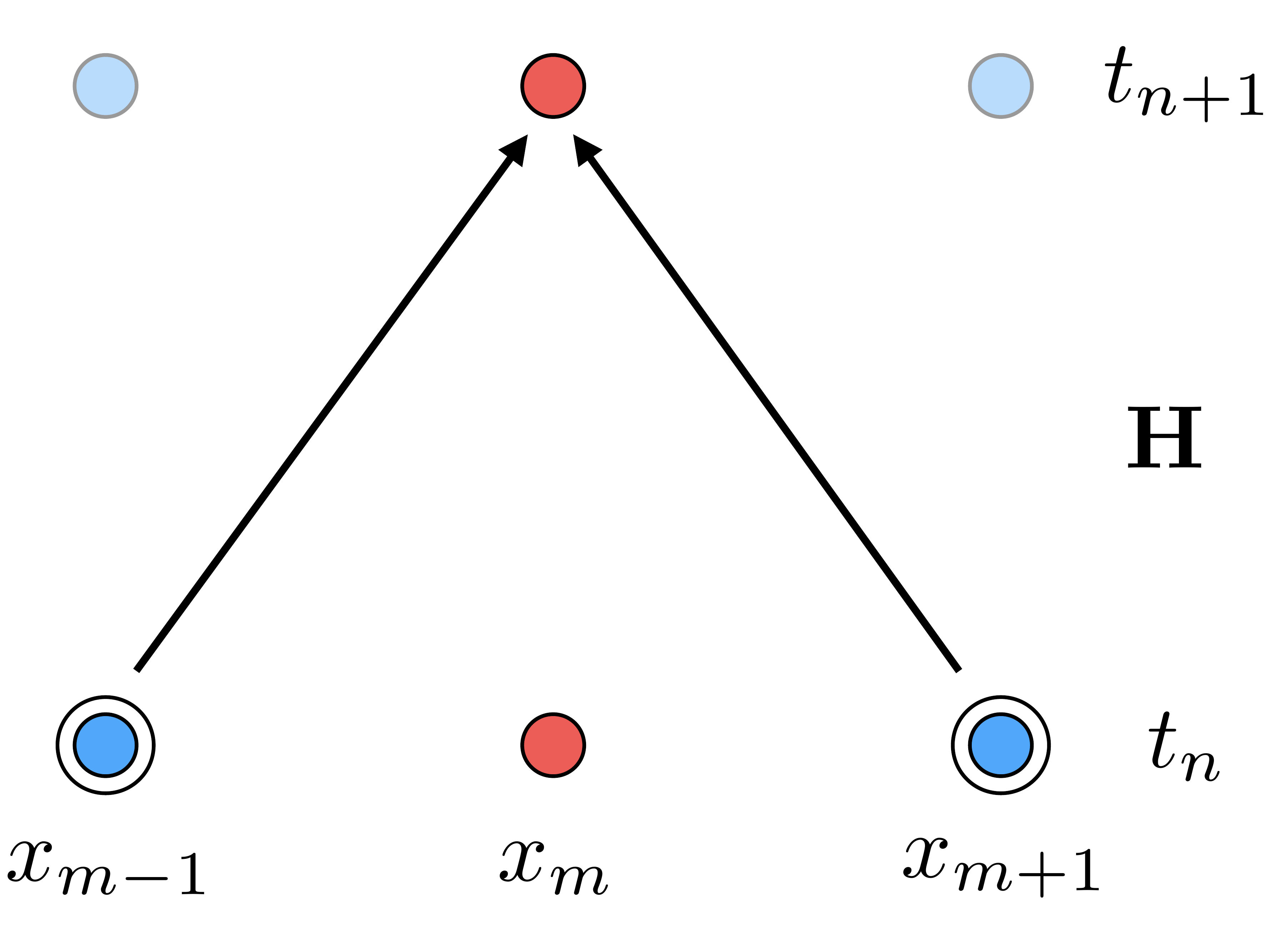}
\caption{Central Hermite interpolation procedure, where a Hermite polynomial is constructed centered at $x_i$ using nodal information from $x_{i+1}$ and $x_{i-1}$.  Nodes which contribute data to the reconstruction are circled.}
\label{fig:central}
\end{figure}
Then, nodal data at $x_m$ is constructed using data from nodes at $x_{m-1}, x_{m+1}$.  
In other words, Hermite interpolation at a node on the primary grid is performed using its two neighbors and then evolved in time, resulting in an update step
\[
\mb{U}_{m} = \mb{R}\mb{T}^{dt}\mb{H}\LRs{\begin{array}{c}
\mb{U}_{m-1}\\ \mb{U}_{m+1}
\end{array}}.
\]
The Central Hermite method may thus also be interpreted as two decoupled Dual Hermite methods on grids of size $2h_x$.  The timestep restriction and numerical results in Section~\ref{sec:num1D} also confirm this interpretation.  

We note that by increasing the size of interpolation interval $h$, the resulting timestep restriction increases independently of the grid spacing $h_x$.  However, doing so also decreases the quality of the interpolation procedure, and Section~\ref{sec:spectra}  describes deleterious effects on the error and spectra of the resulting method.  
 
The Central Hermite interpolation operator results in the 2-node stencil of Figure~\ref{fig:central}, since nodal information from $x_{m-1}$ and $x_{m+1}$ is required to construct information at a $x_m$.  Computationally, a smaller stencil results in fewer memory accesses for the reconstruction.  While the difference between the Central Hermite and Virtual Hermite stencils in one space dimension is small, the difference becomes more pronounced in multiple dimensions.  For a degree $N$ Hermite method in $d$ dimensions, each node in the stencil requires $(N+1)^d$ accesses, and a Central Hermite stencil contains $4$ nodes in 2D, and $8$ nodes in 3D, while the Virtual Hermite stencil contains $9$ nodes in 2D and $27$ nodes in 3D.  

\subsection{Upwind Hermite methods} 

Each Hermite method presented has utilized a centered stencil, where solution values and derivatives are interpolated in a symmetric fashion around the reconstruction point.  The Upwind Hermite method constructs instead a directional or one-sided Hermite reconstruction.  This concept was used in \cite{goodrich2006hermite} to enforce boundary conditions, though the use of such reconstructions may also take advantage of the directional nature of hyperbolic equations \cite{courant1952solution, lesaint1974finite}.  

\begin{figure}[h!]
\centering
\subfloat[One-dimensional reconstruction]{\includegraphics[width=.4\textwidth]{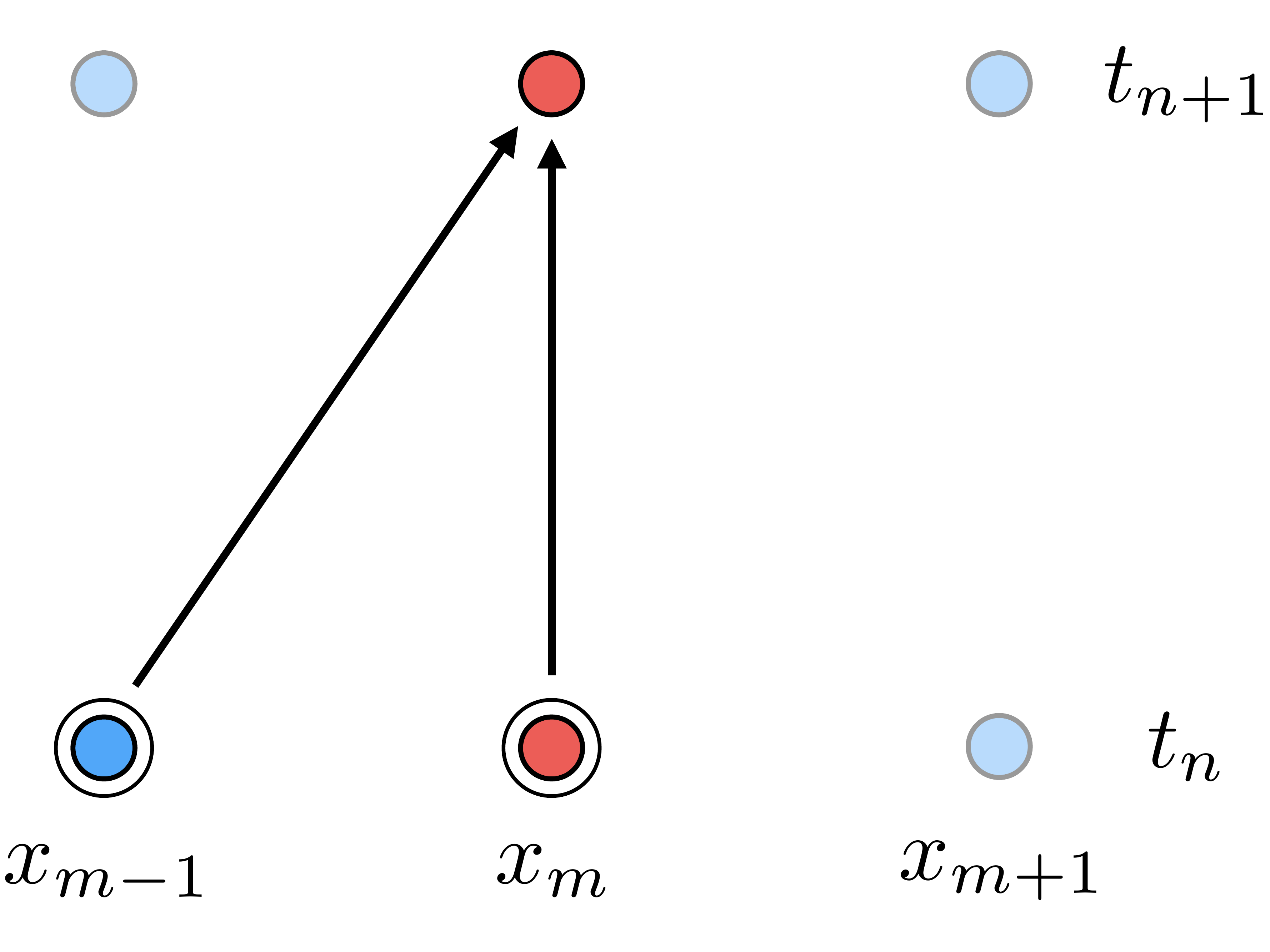}}
\caption{Interpolation procedure in one space dimension, where a Hermite polynomial is constructed centered at $x_m$ using additional nodal information from and $x_{m-1}$.  Nodes which contribute data to the reconstruction are circled.}
\label{fig:upwind}
\end{figure}

As shown in Figure~\ref{fig:upwind}, an Upwind Hermite reconstruction $\tilde{u}_m(x)$ may be defined at $x_m$ by interpolating solution and derivative values at the endpoints of the interval $(x_{m-1}, x_{m})$
\[
\left.\frac{\partial^i \tilde{u}_m}{\partial x^i}\right|_{x_{m-1}} = \left.\frac{\partial^i u_{m-1}}{\partial x^i}\right|_{x_{m-1}}, \qquad 
\left.\frac{\partial^i \tilde{u}_m}{\partial x^i}\right|_{x_m}= \left.\frac{\partial^i u_m}{\partial x^i}\right|_{x_m}, \qquad i = 0,\ldots, 2N+1.
\]
The solution to this problem results in a degree $2N+1$ reconstruction at $x_m$.  However, since the solution and derivative values at $x_m$ are interpolated, the first $(N+1)$ coefficients at $x_m$ remain unchanged, and only the remaining $(N+1)$ coefficients need to be computed.  This may be done by multiplying the derivative values at $x_{m-1}, x_m$ by the last $(N+1)$ rows of the interpolation matrix $\mb{H}$.  A downwind interpolation operator may be defined in a similar manner using information at $(x_m, x_{m+1})$ to produce a reconstruction at $x_m$.  

For the advection equation specifically, the domain of influence is biased, such that the solution at time $t_n + dt$ at the point $x_m$ depends only on $x_m - cdt$ at time $t_n$.  Redefining the width of the interpolation interval for the Upwind Hermite method as $h = x_m-x_{m-1} = h_x$, the same domain of dependence arguments used previously imply that the method is stable if $cdt < h_x$.  The timestep restriction of the Upwind Hermite method then matches that of the Central Hermite method, with the caveat that this result is specific to scalar advection equations and the sign of $c$.  For example, for $c>0$, a downwind reconstruction would be unstable due to the fact that the interval $(x_m, x_{m+1})$ does not contain the domain of dependence for $x_m$ at any time greater than $t_n$.  We also note that, by similar arguments made in \cite{goodrich2006hermite}, the time evolution of the Upwind Hermite method by temporal Taylor series is also exact.   

The Upwind Hermite method also requires special treatment when directionality is not readily apparent, such as for systems of hyperbolic equations.  In one space dimension, the procedure may be adapted to reconstruct upwind and downwind characteristic variables, similar to the approach used in WENO reconstructions \cite{qiu2002construction, ren2003characteristic}.  However, the effectiveness of the characteristic approach does not appear to extend to all systems of equations in higher dimensions, as discussed in Section~\ref{sec:2D}.

\section{Numerical experiments in 1D}
\label{sec:num1D}
To compare the performance of the new Hermite methods, we examine convergence rates and qualitative behavior for the Virtual. Central, and Upwind Hermite methods.  Numerical results are shown for the periodic constant-coefficient scalar advection equation on the interval $[-1,1]$, using the Taylor expansion discussed in Section~\ref{sec:time} to evolve in time.  

We introduce also a CFL constant $C > 0$ such that $dt = Ch/c$, where $h$ is the size of the interpolation interval (for Virtual Hermite methods, $h = h_x/2$, while for Central and Upwind Hermite methods, $h = h_x$).  $C\approx 1$ sets the timestep as large as possible based on the timestep restriction for each method, while $C \ll 1$ results in more timesteps than necessary as implied by stability.  Since a filter-like step (the Hermite reconstruction) is applied at each timestep, small values of the CFL constant $C$ (i.e.\ smaller timesteps than strictly necessary) may result in stronger filtering than necessary and larger errors.  

\subsection{Convergence rates}

\begin{figure}[!h]
\centering
\subfloat[$C = .1$]{\includegraphics[width=.475\textwidth]{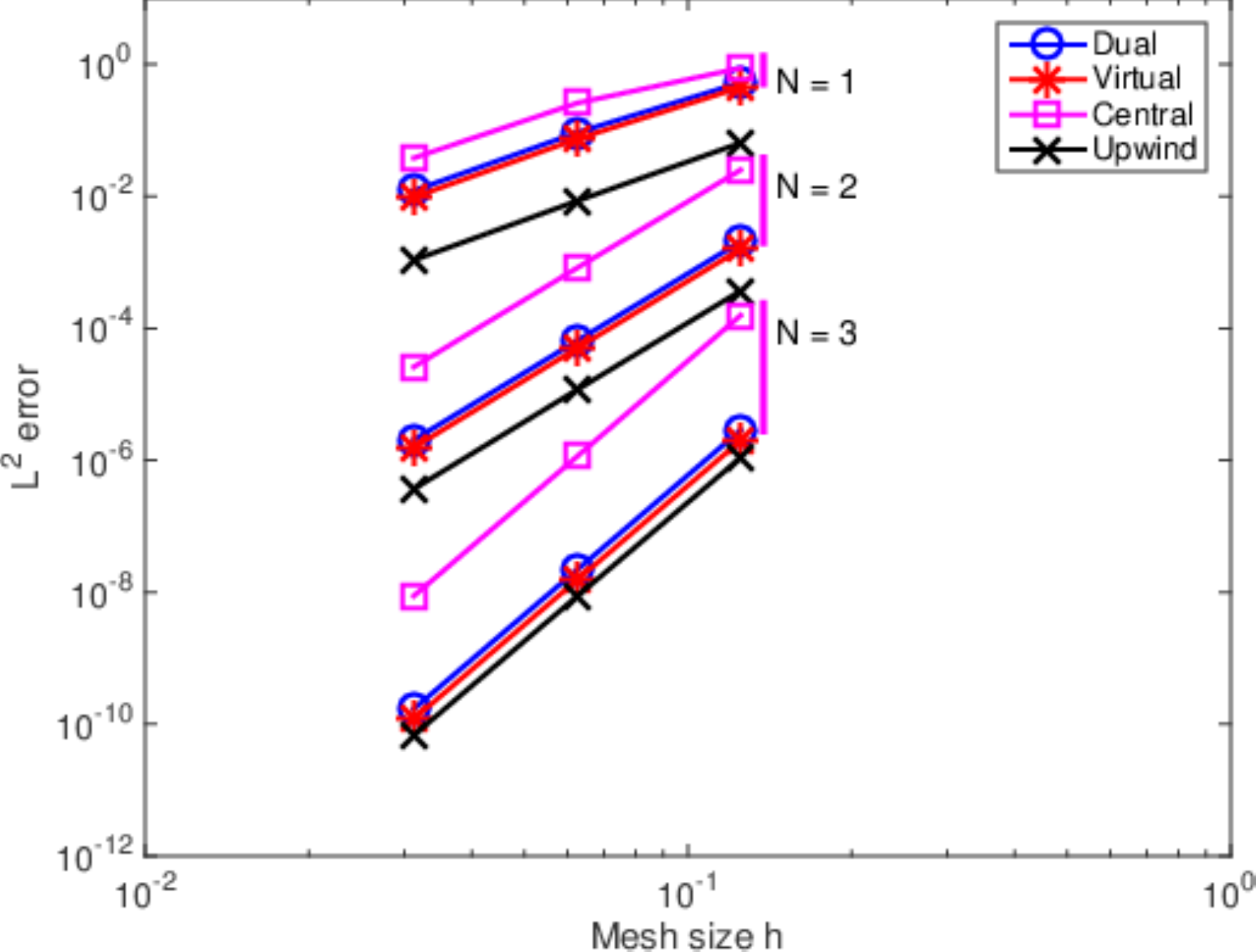}}
\subfloat[$C = .5$]{\includegraphics[width=.475\textwidth]{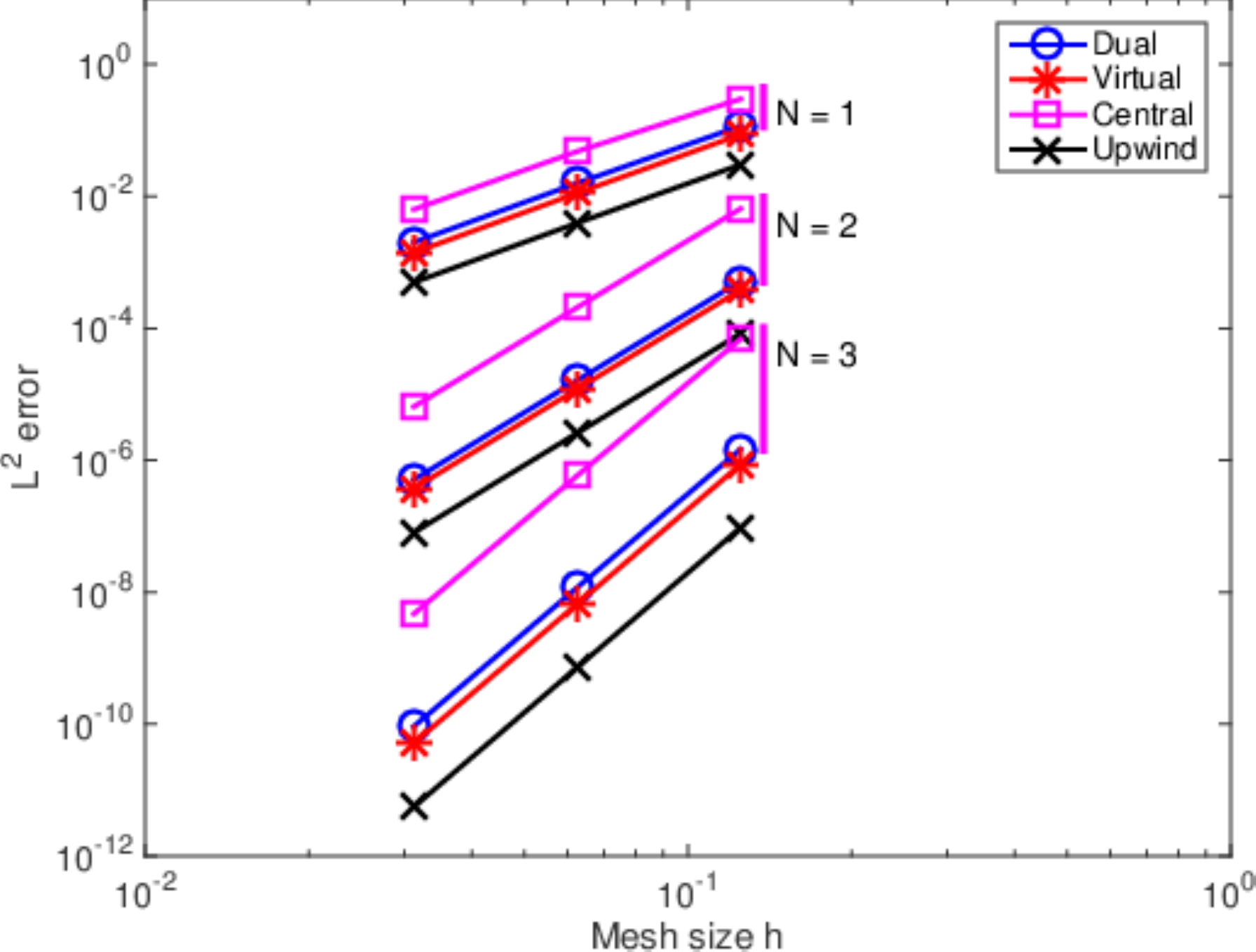}}\\
\subfloat[$C = .9$]{\includegraphics[width=.475\textwidth]{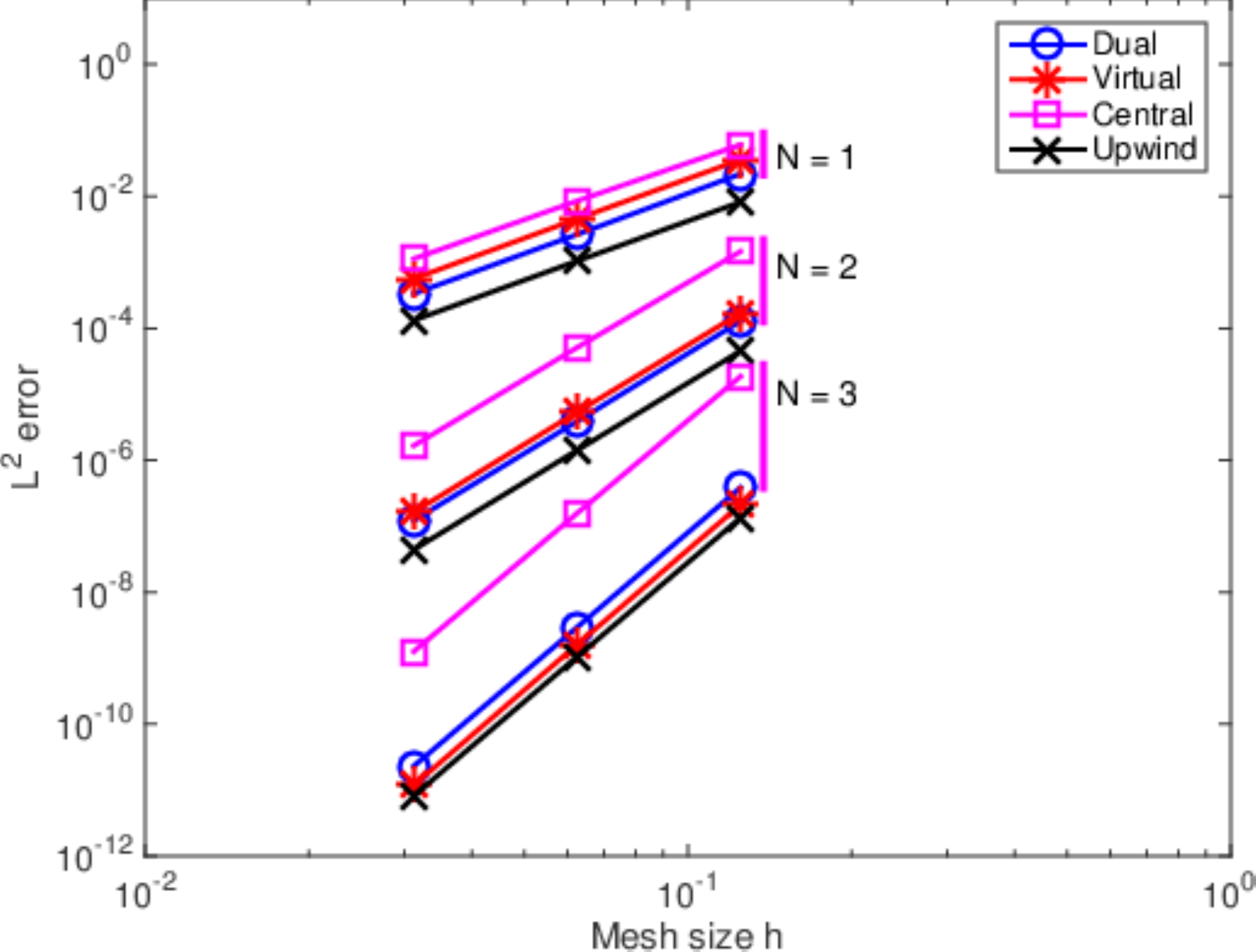}}
\caption{Convergence of $L^2$ errors for Hermite methods for the one-dimensional advection equation with a smooth sine solution.}
\label{fig:rates1}
\end{figure}

We report convergence rates for the one-dimensional scalar advection equation with speed $a=1$ and solution 
\[
u(x,t) = \sin(\pi (x-t)).
\]
We vary the CFL constant between $C = .1$, $C=.5$, and $C=.9$ and calculate $L^2$ errors at time $T=10$.  For both methods, the error is smaller the closer $C$ is to $1$ as shown in Figure~\ref{fig:rates1}.  

Convergence rates are reported in Table~\ref{table:rates1}, and except for the lowest order case $N=1$ and smallest CFL constant $C=.1$, optimal $O(h^{2N+1})$ rates of convergence were observed for all methods.  Additionally, at higher orders of approximation, errors for the Virtual and Upwind Hermite method are very similar to those of the Dual Hermite method.  For the Central Hermite method, the error is roughly a factor of $2^N$ greater than that of the other methods.  

\begin{table}[h!]
\centering                                                          
\begin{tabular}{|c||c|c|c|c|c|c|c|c|c|}
\hline
& \multicolumn{3}{c|}{$C=.1$} &  \multicolumn{3}{c|}{$C=.5$}  &  \multicolumn{3}{c|}{$C=.9$} \\ 
\hline
$N$ & $1$ & $2$ & $3$ & $1$ & $2$ & $3$  & $1$ & $2$ & $3$ \\ 
\hhline{|=|=|=|=|=|=|=|=|=|=|}
Dual & 2.72 & 4.99 & 7.02 & 2.93 & 5.0 & 6.98 & 3.02 & 5.02 & 7.01\\        
\hline
Virtual  &   2.67 &    5.0  &  7.00 &   2.96  &  5.0 &   7.02 &   2.99 &  5.01 &  7.07\\
\hline
Central & 1.71 & 4.94 & 7.06 & 2.62 & 4.98 & 6.92 & 2.92 & 4.98 & 6.99\\        
\hline                                                              
Upwind & 2.94 &   4.99 &   6.98 &   2.96 &   5.0 &  7.02 &   3.02 &  5.03 &   7.03\\
\hline                                                              
\end{tabular}                                                       
\caption{$L^2$ rates of convergence for Hermite methods for the one-dimensional advection equation with a smooth sinusoidal solution.}
\label{table:rates1}
\end{table}              

\begin{figure}[!h]
\centering
\subfloat[Dual Hermite]{\includegraphics[width=.4\textwidth]{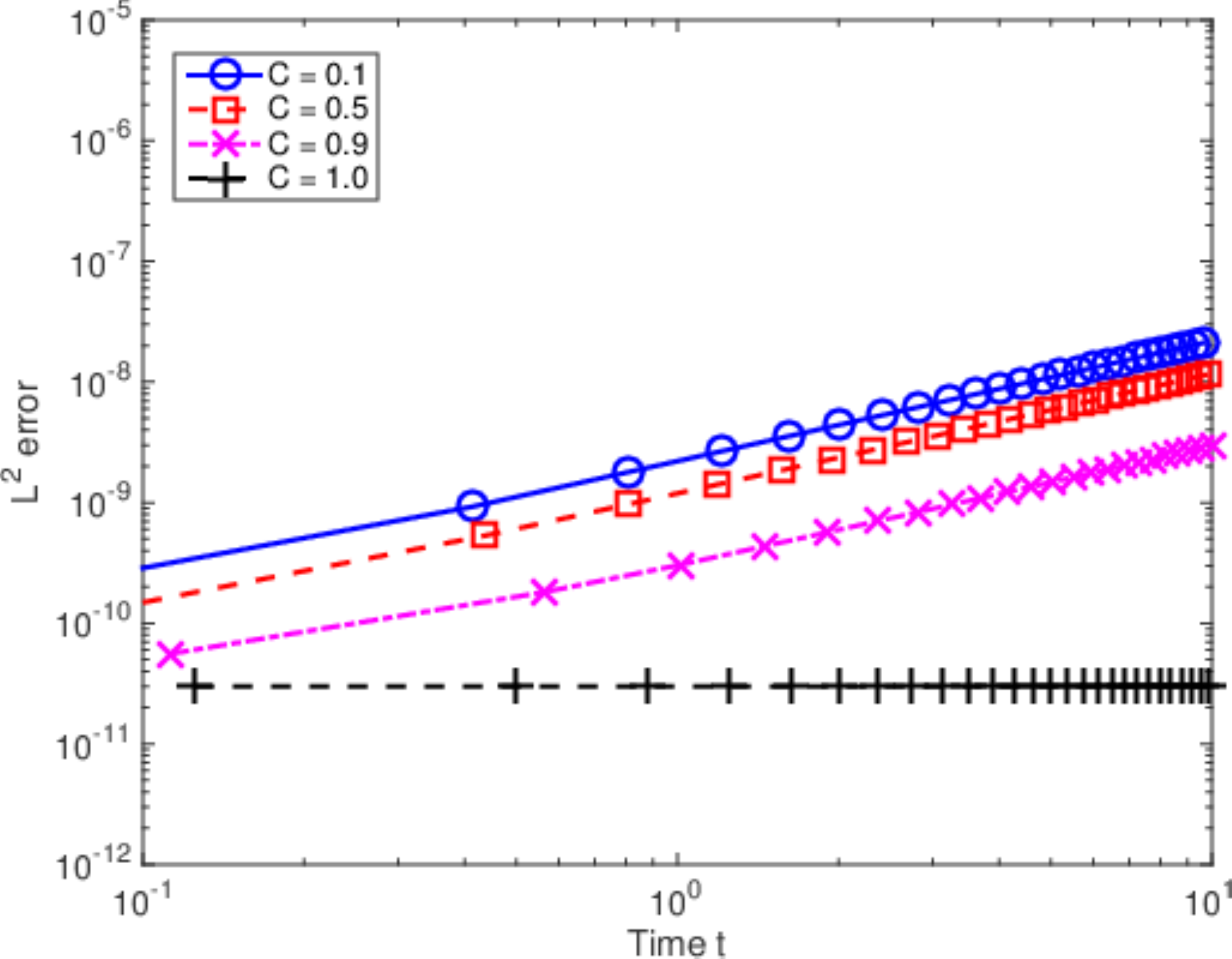}}
\hspace{2em}
\subfloat[Virtual Hermite]{\includegraphics[width=.4\textwidth]{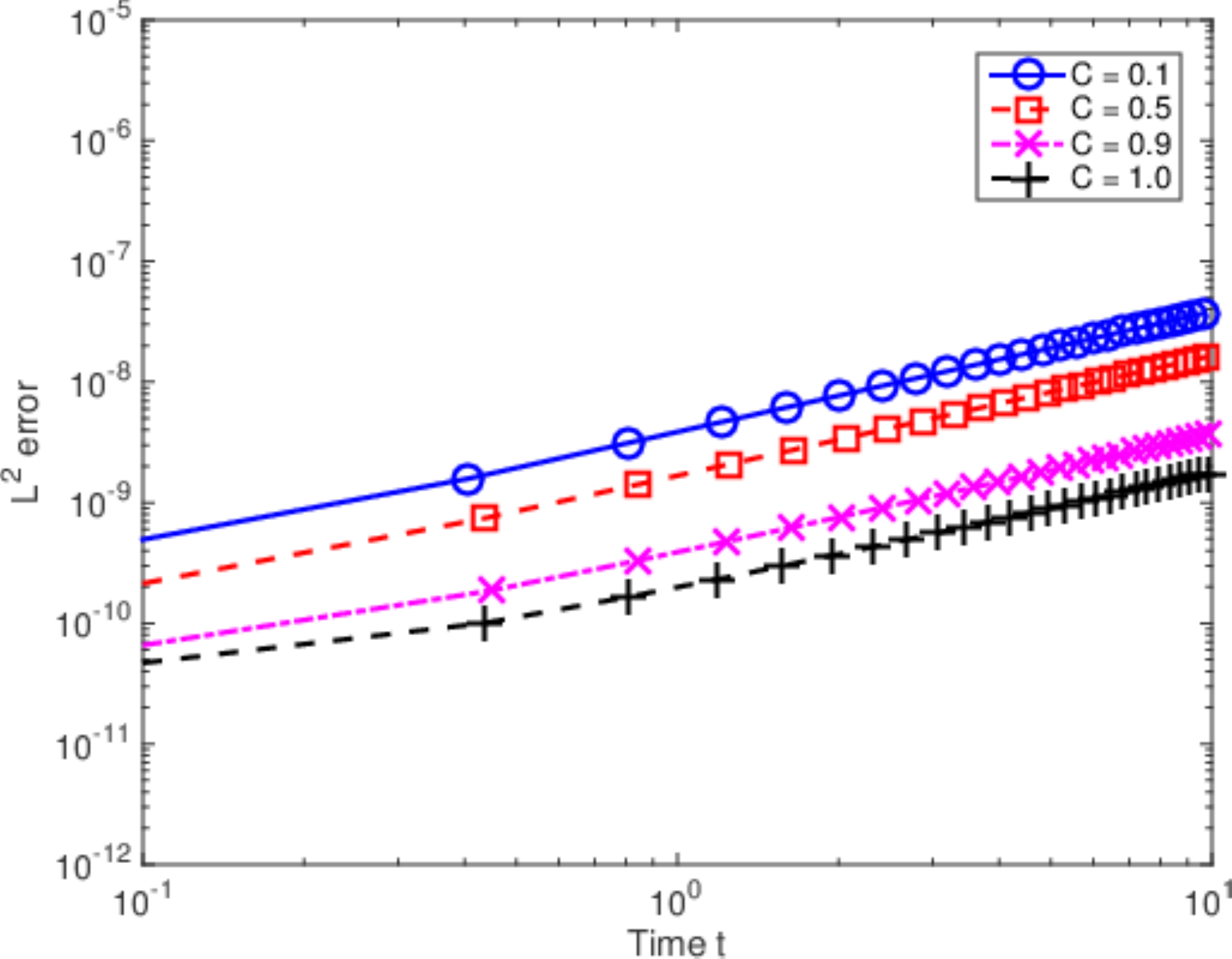}}\\
\subfloat[Central Hermite]{\includegraphics[width=.4\textwidth]{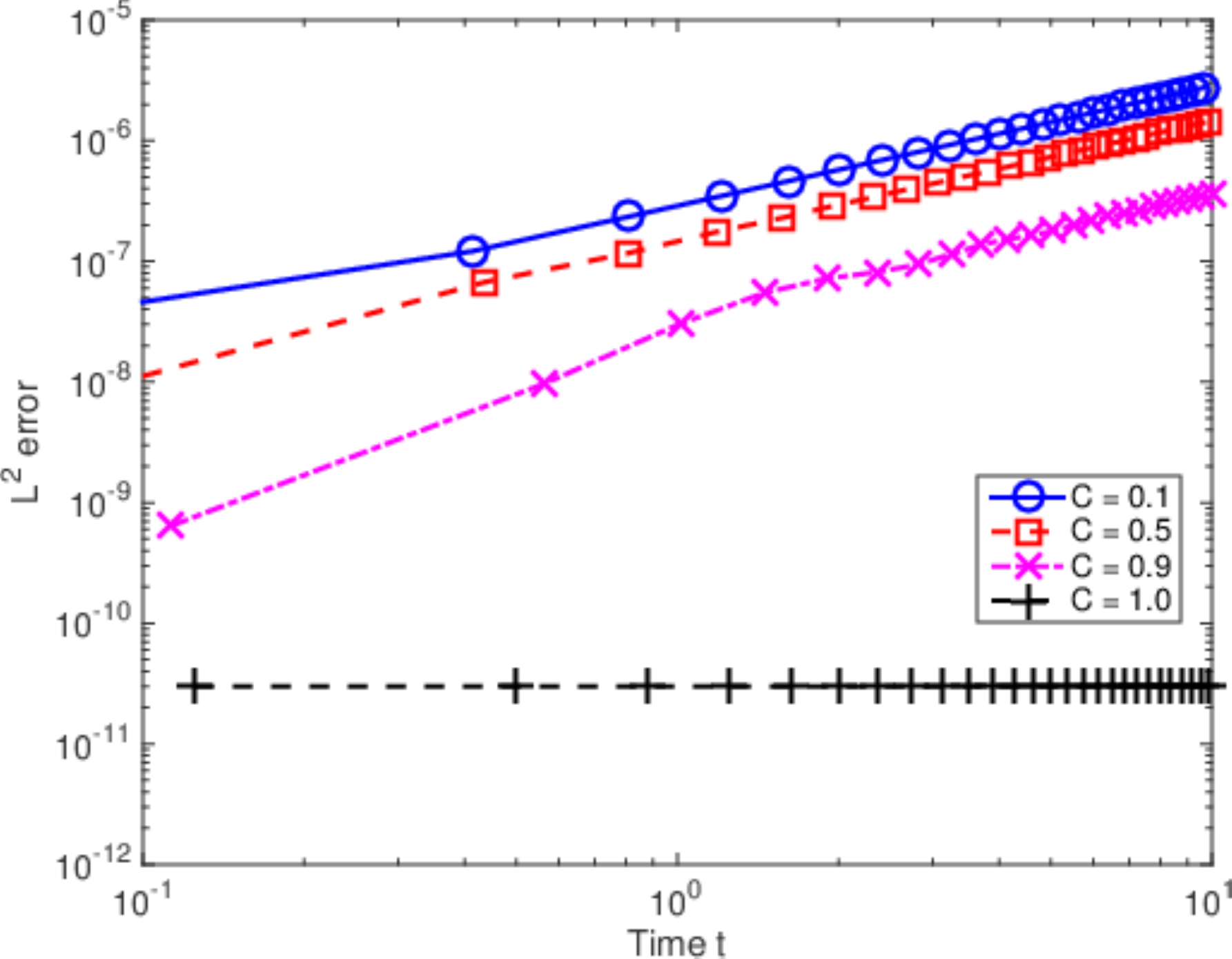}}
\hspace{2em}
\subfloat[Upwind Hermite]{\includegraphics[width=.4\textwidth]{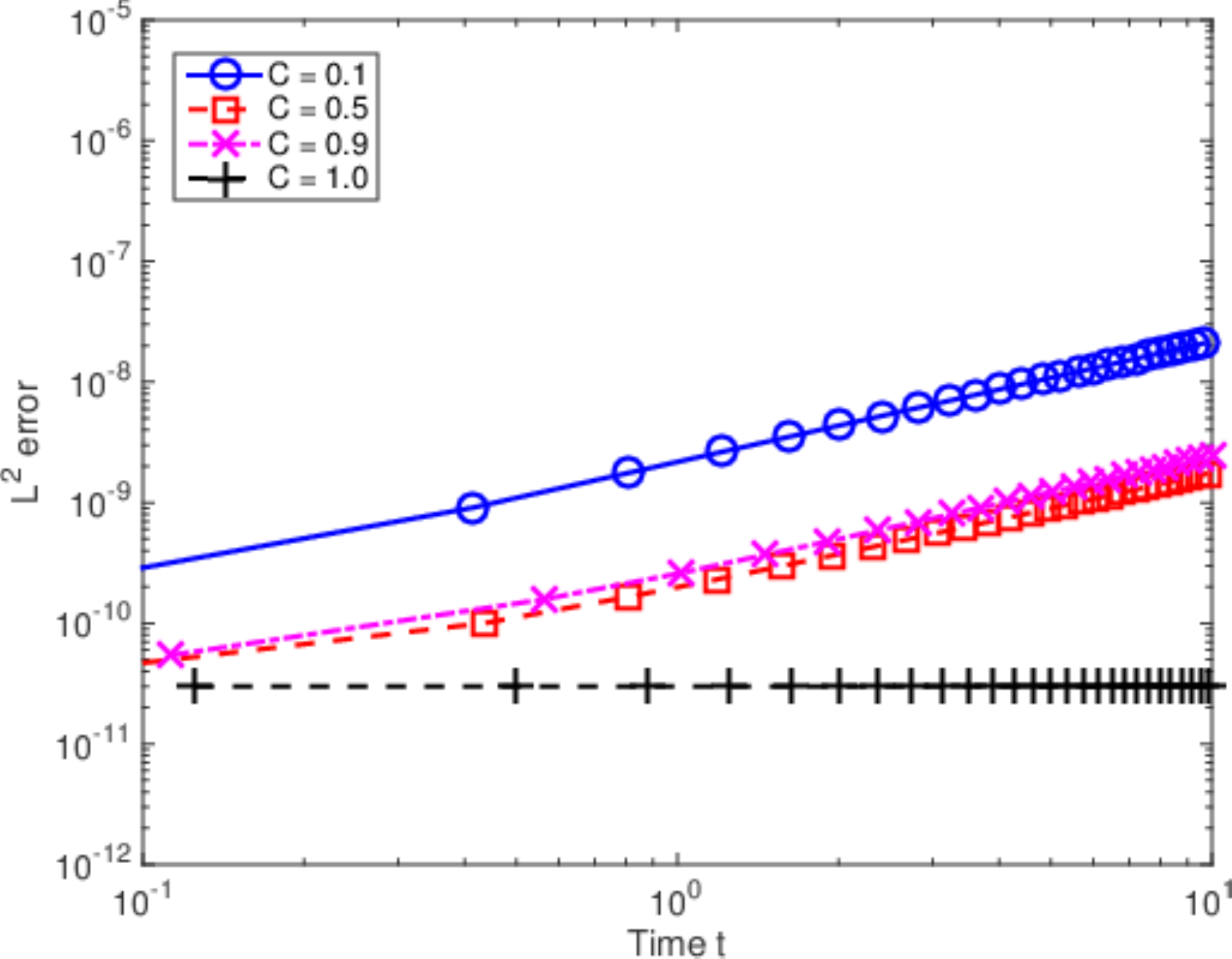}}
\caption{Growth of $L^2$ error in time for various Hermite schemes with $N=3$, $K = 16$ in one dimension.  The advection equation is solved up to time $T=10$ with a smooth sine solution, and the $L^2$ error is computed at each timestep.}
\label{fig:time_err}
\end{figure}

The growth of error in time is also examined for each Hermite method.  Error estimates for discretizations of transient hyperbolic problems typically contain two terms which characterize spatially-dependent and time-dependent errors, respectively.  Standard bounds are of the form
\[
e(T) \leq (C_1 + C_2 T) h^{r(N)}
\]
where $e(T)$ is some measure of error at time $T$, and $r(N)$ is some rate of convergence depending on the degree of approximation \cite{cockburn1998local, hesthaven2007nodal}.  We confirm the linear growth of error in time for all Hermite methods in Figure~\ref{fig:time_err}, with the exception of the case when $C=1$, which is discussed in Section~\ref{sec:spectra}.  The growth of error for the Upwind and Dual Hermite methods is very similar.  While the Central Hermite method develops larger errors than the Virtual Hermite method, the long-time rate of growth is identical for each value of $C$ for both methods.  Moreover, the Central Hermite method with $K=32$ results in time-dependent errors very similar to the Virtual Hermite method for $K=16$, indicating a strong dependence of the error on the size of the interpolation interval.  

\begin{figure}
\centering
\subfloat[Dual, $C = .1$]{\includegraphics[width=.3\textwidth]{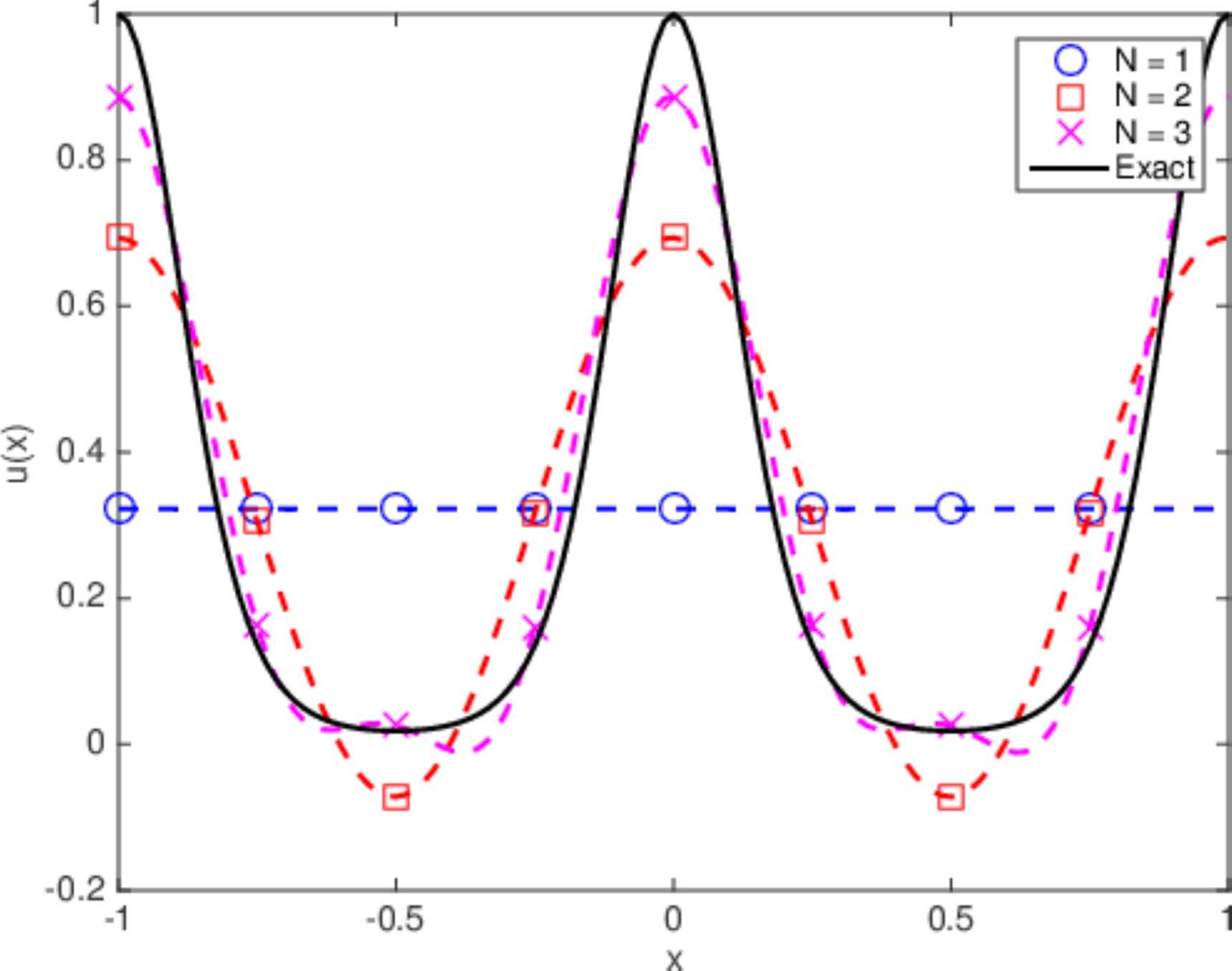}}
\hspace{1em}
\subfloat[Dual, $C = .5$]{\includegraphics[width=.3\textwidth]{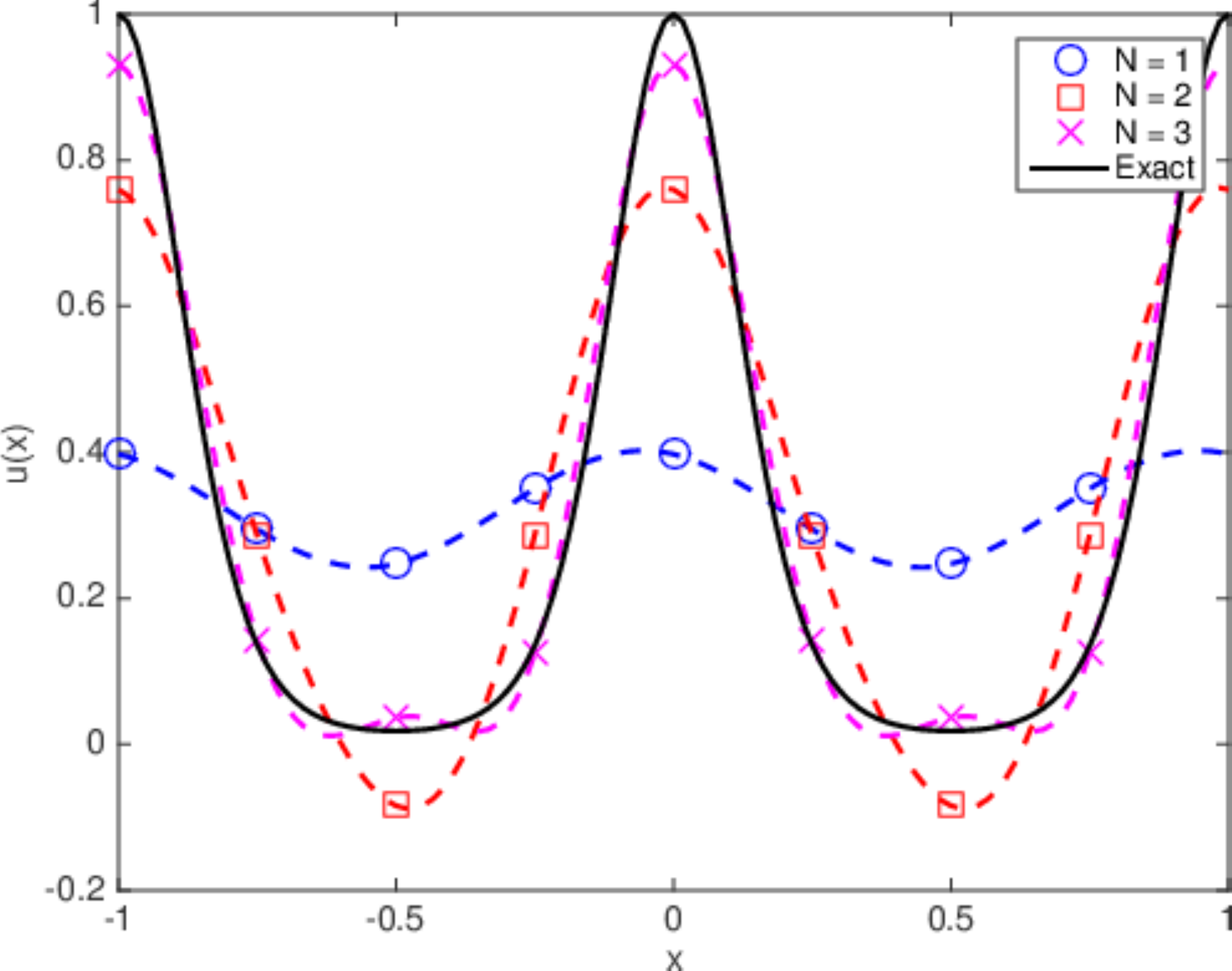}}
\hspace{1em}
\subfloat[Dual, $C = .9$]{\includegraphics[width=.3\textwidth]{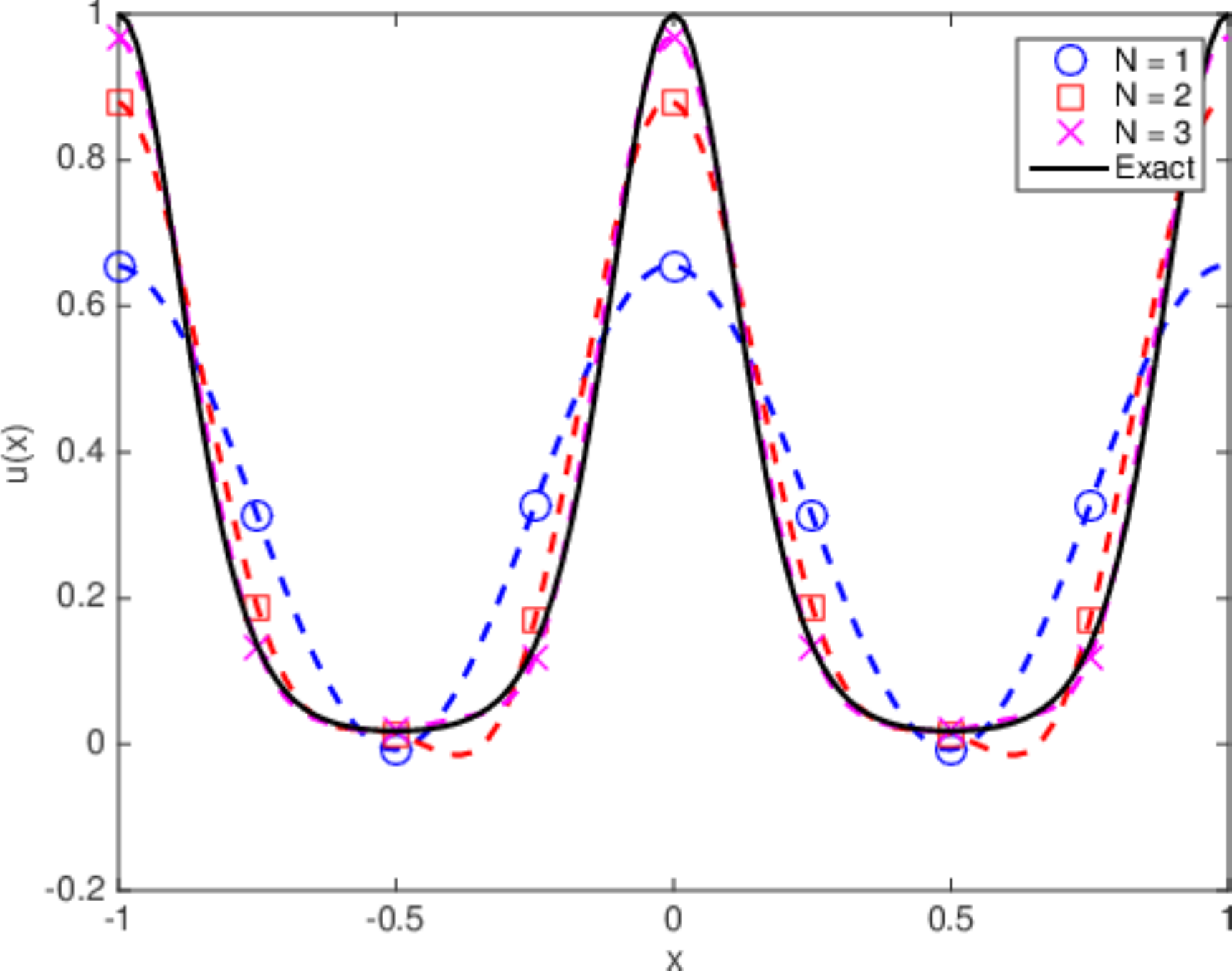}}
 
\subfloat[Virtual, $C = .1$]{\includegraphics[width=.3\textwidth]{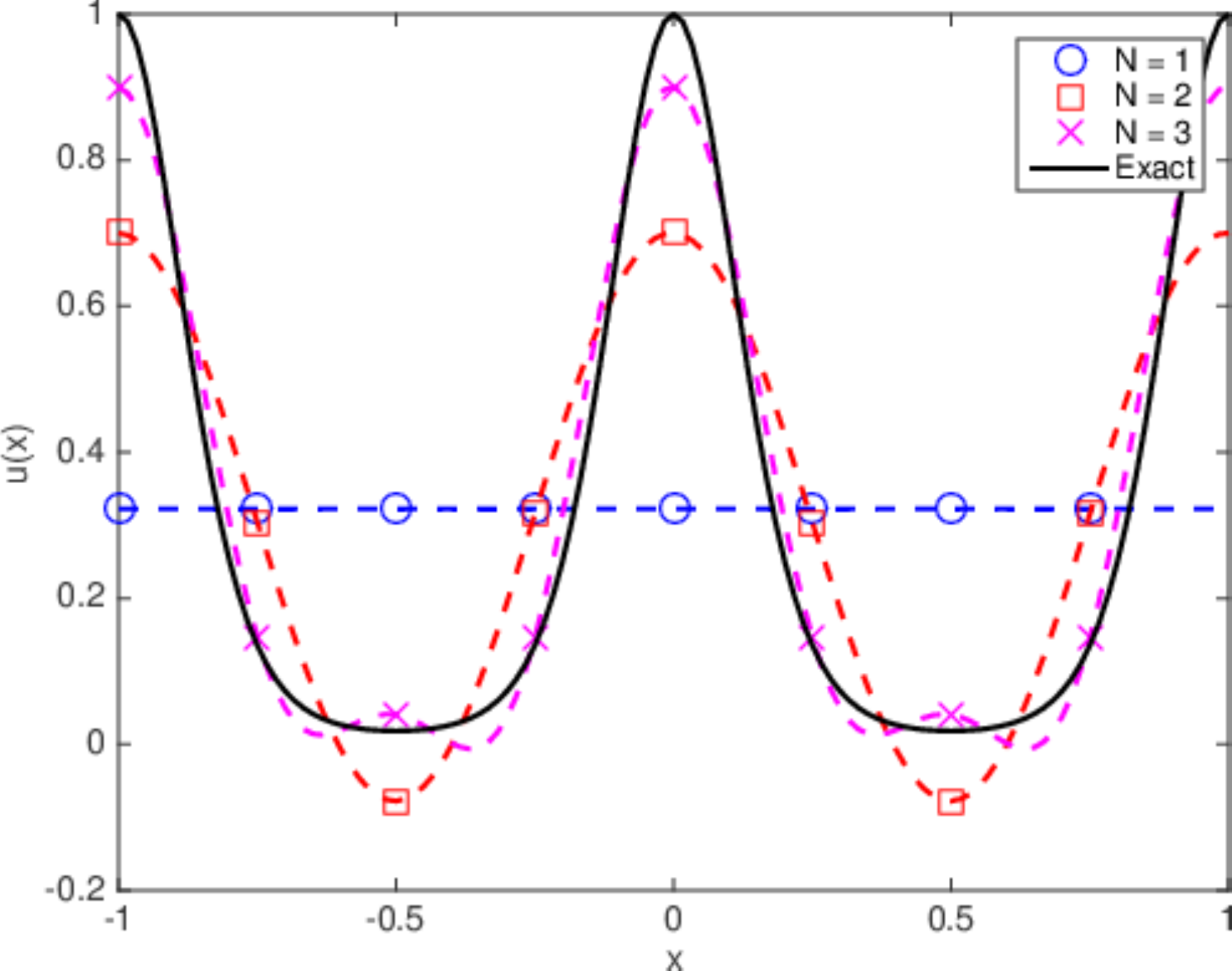}}
\hspace{1em}
\subfloat[Virtual, $C = .5$]{\includegraphics[width=.3\textwidth]{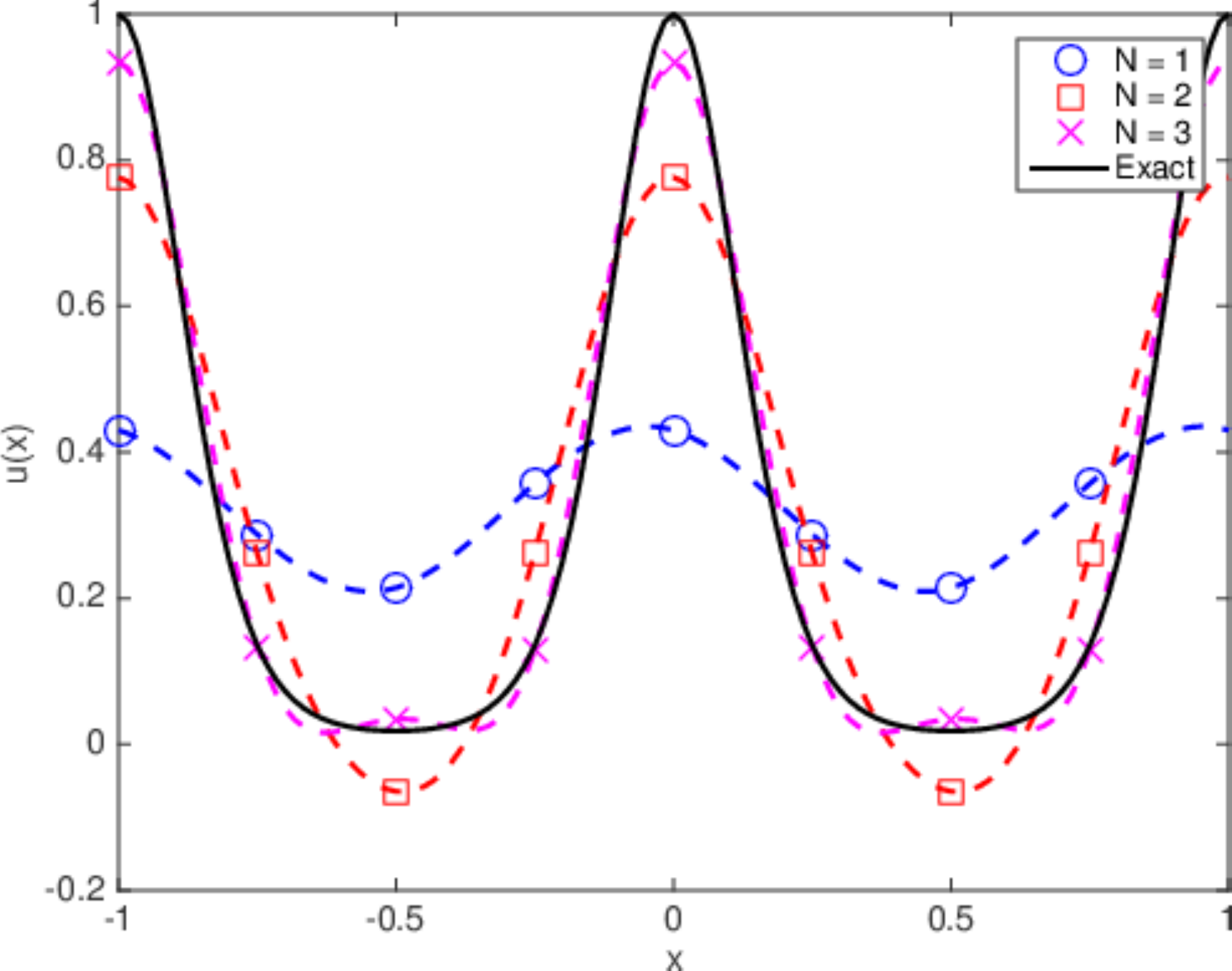}}
\hspace{1em}
\subfloat[Virtual, $C = .9$]{\includegraphics[width=.3\textwidth]{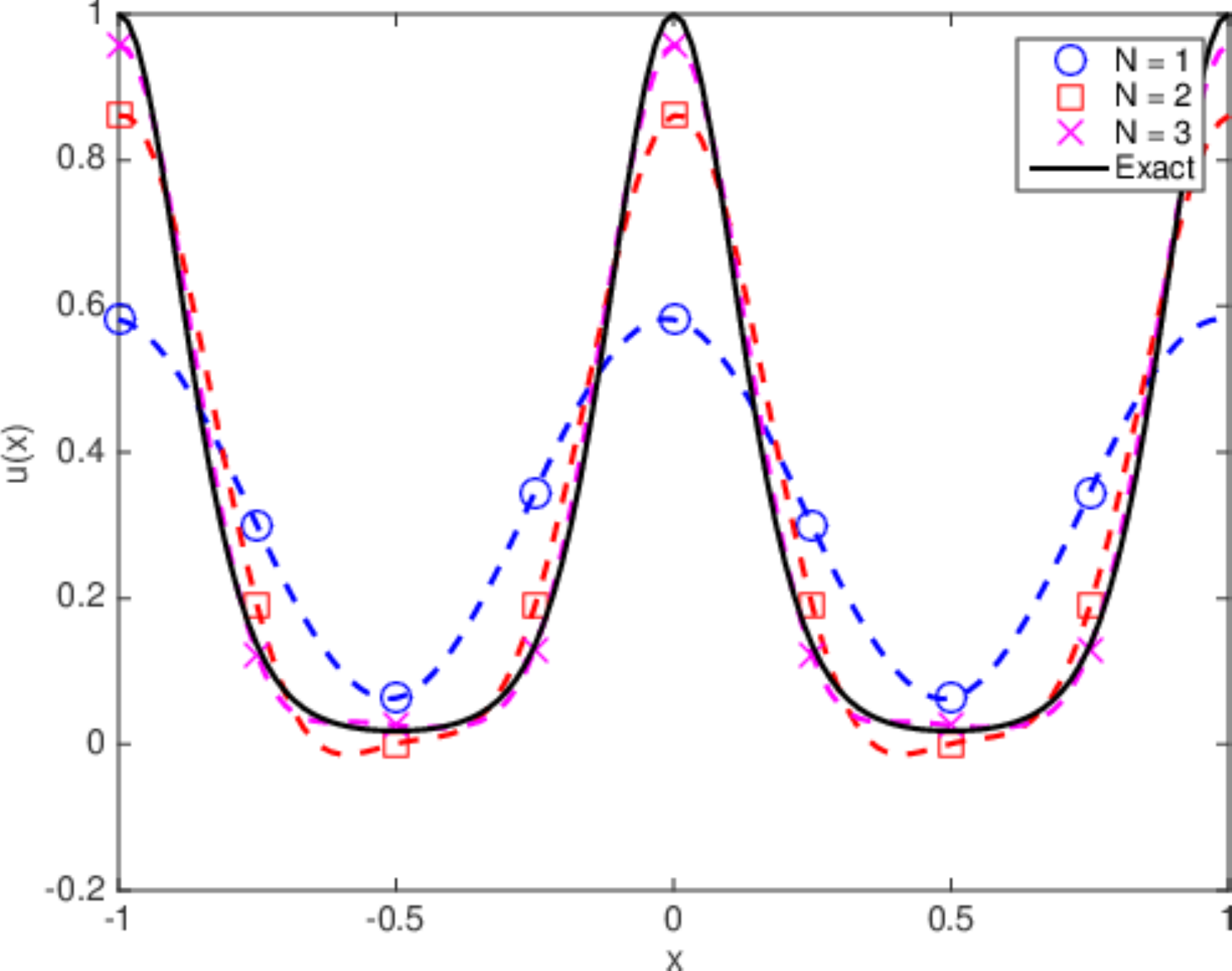}}

\subfloat[Central, $C = .1$]{\includegraphics[width=.3\textwidth]{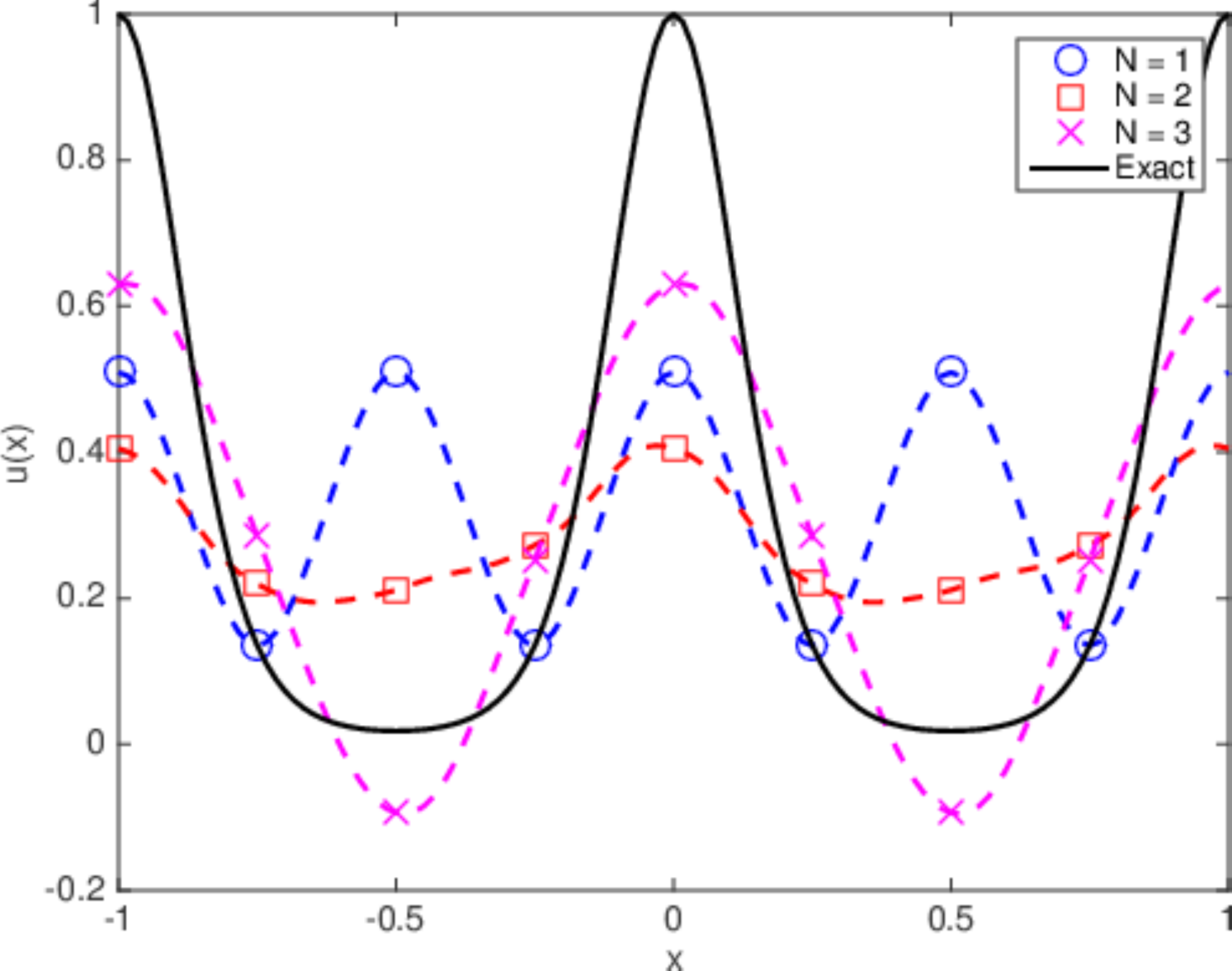}}
\hspace{1em}
\subfloat[Central, $C = .5$]{\includegraphics[width=.3\textwidth]{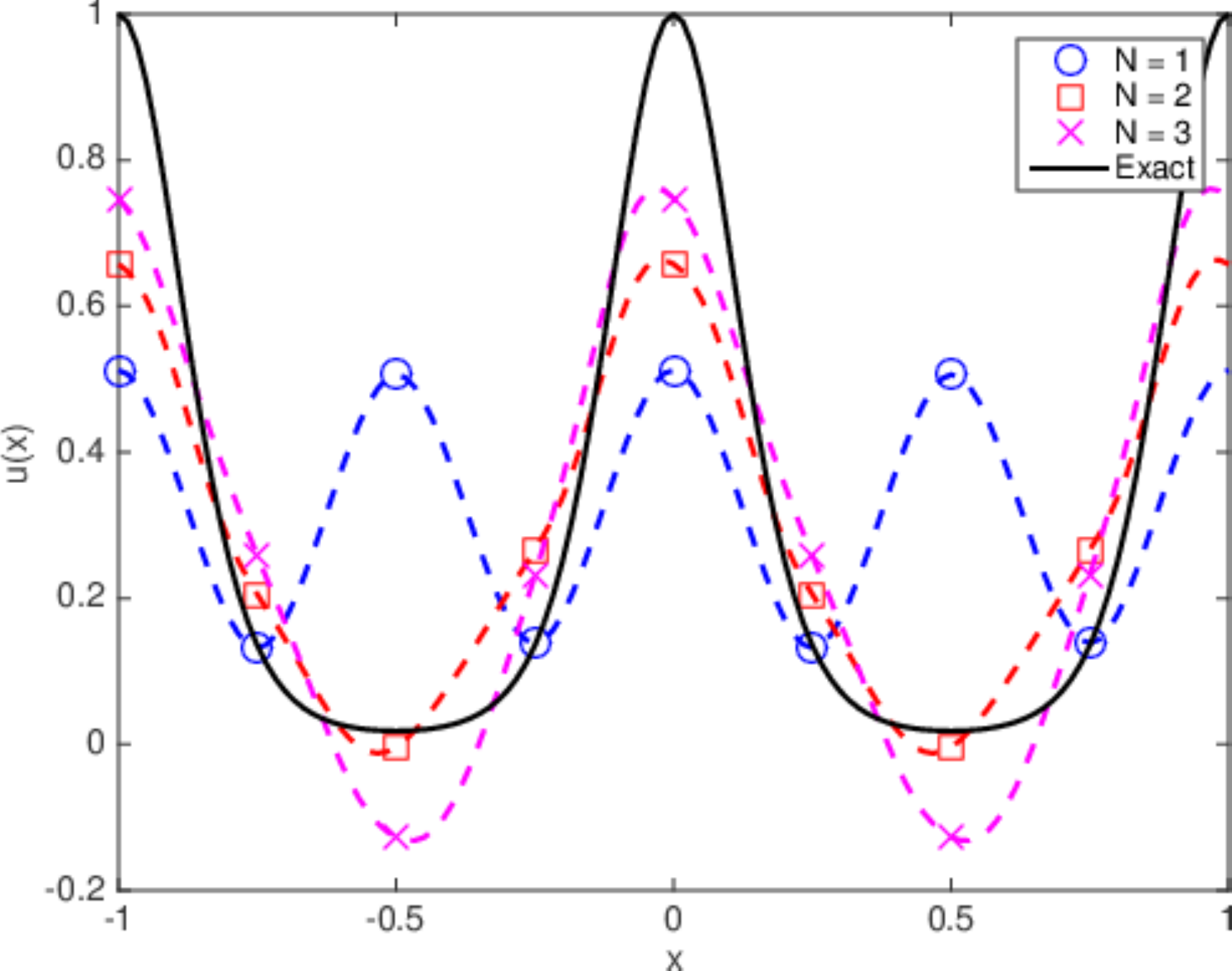}}
\hspace{1em}
\subfloat[Central, $C = .9$]{\includegraphics[width=.3\textwidth]{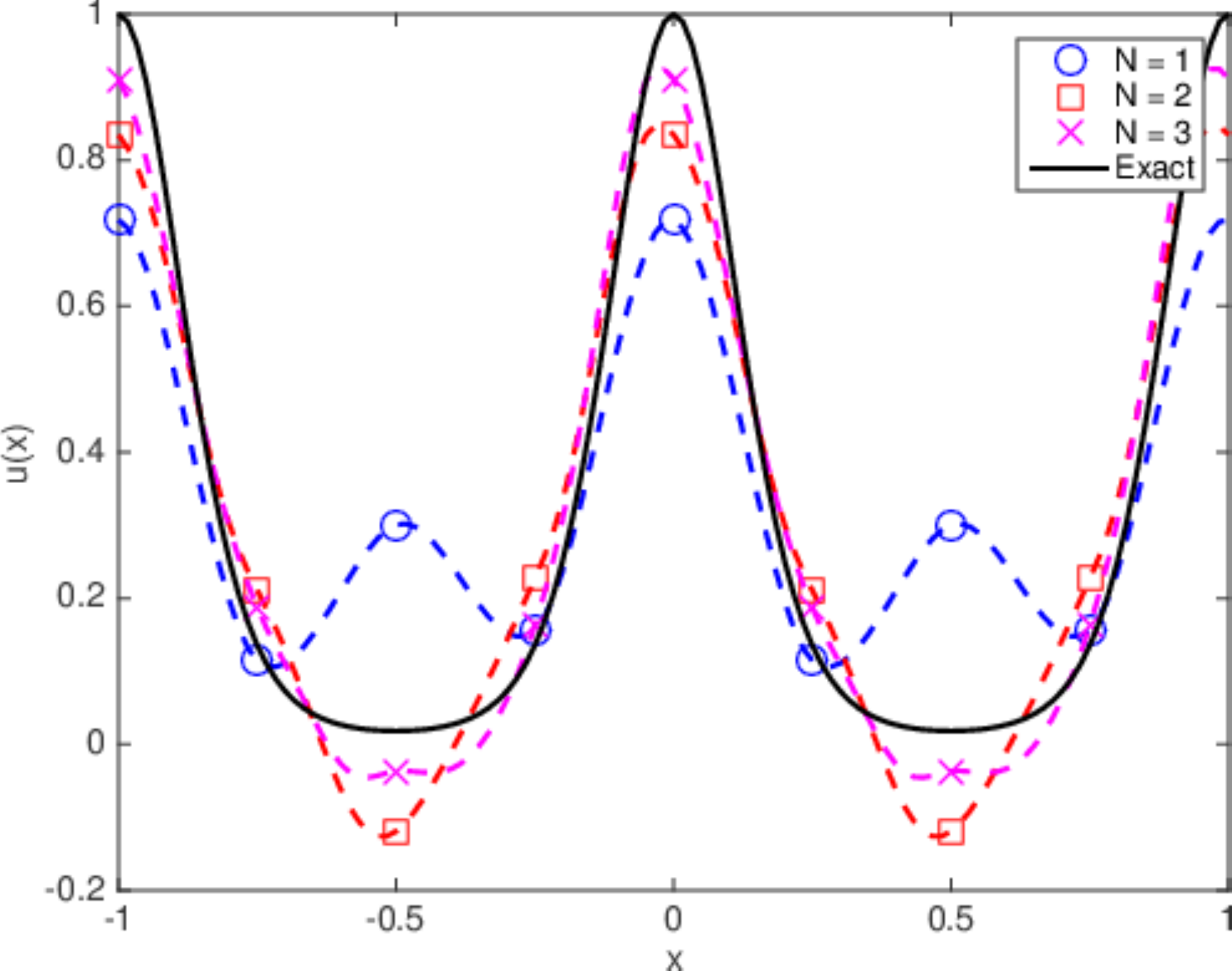}}

\subfloat[Upwind, $C = .1$]{\includegraphics[width=.3\textwidth]{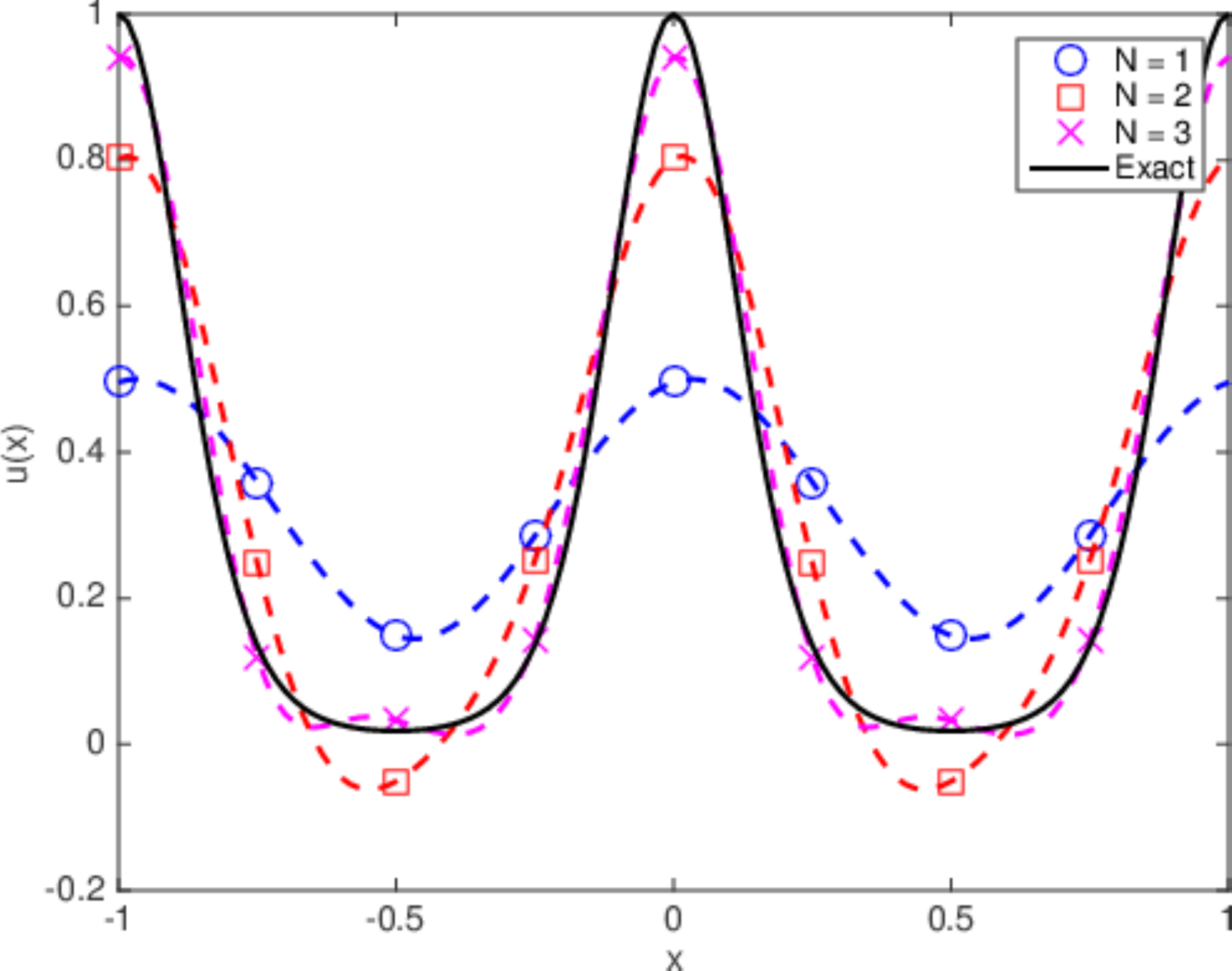}}
\hspace{1em}
\subfloat[Upwind, $C = .5$]{\includegraphics[width=.3\textwidth]{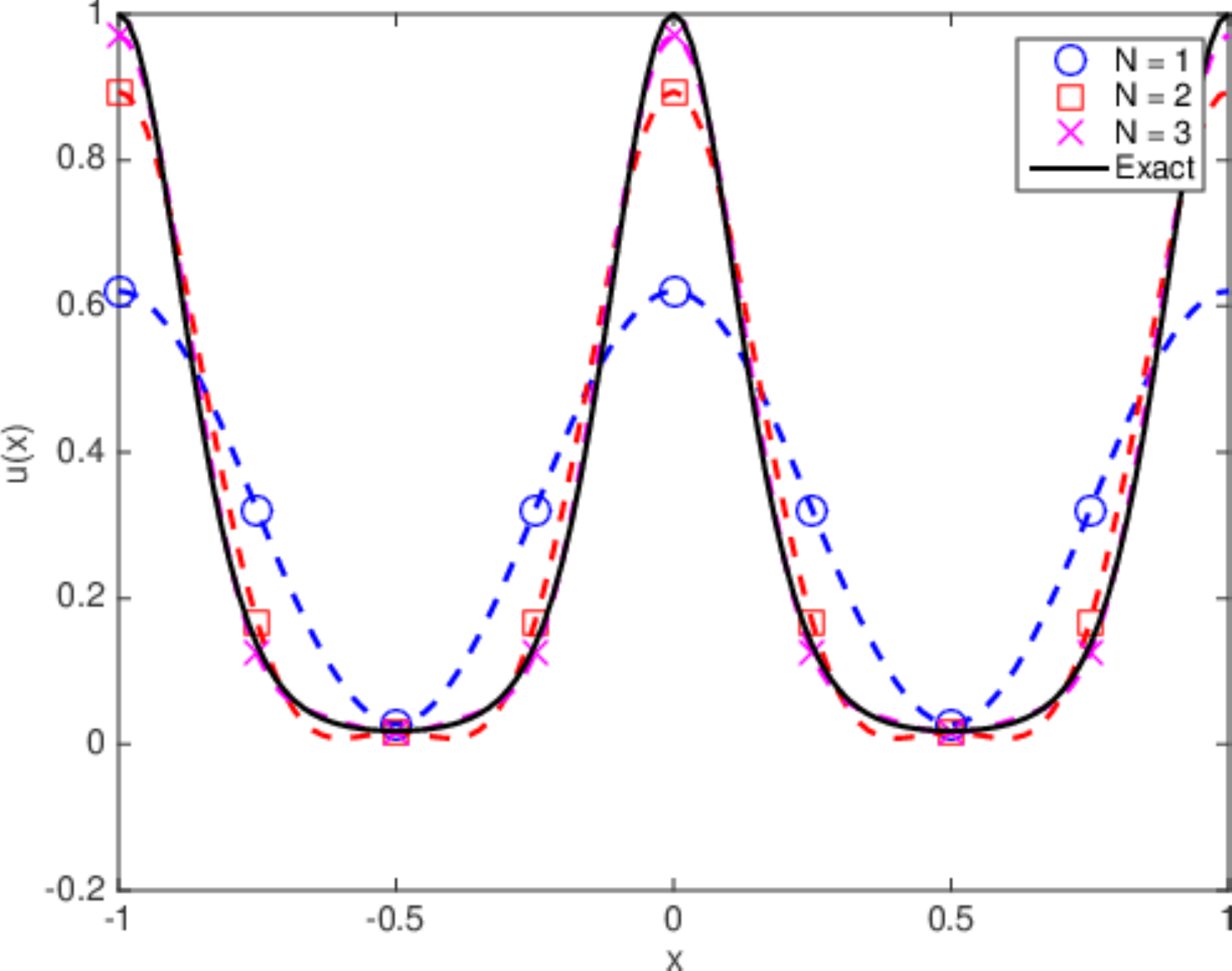}}
\hspace{1em}
\subfloat[Upwind, $C = .9$]{\includegraphics[width=.3\textwidth]{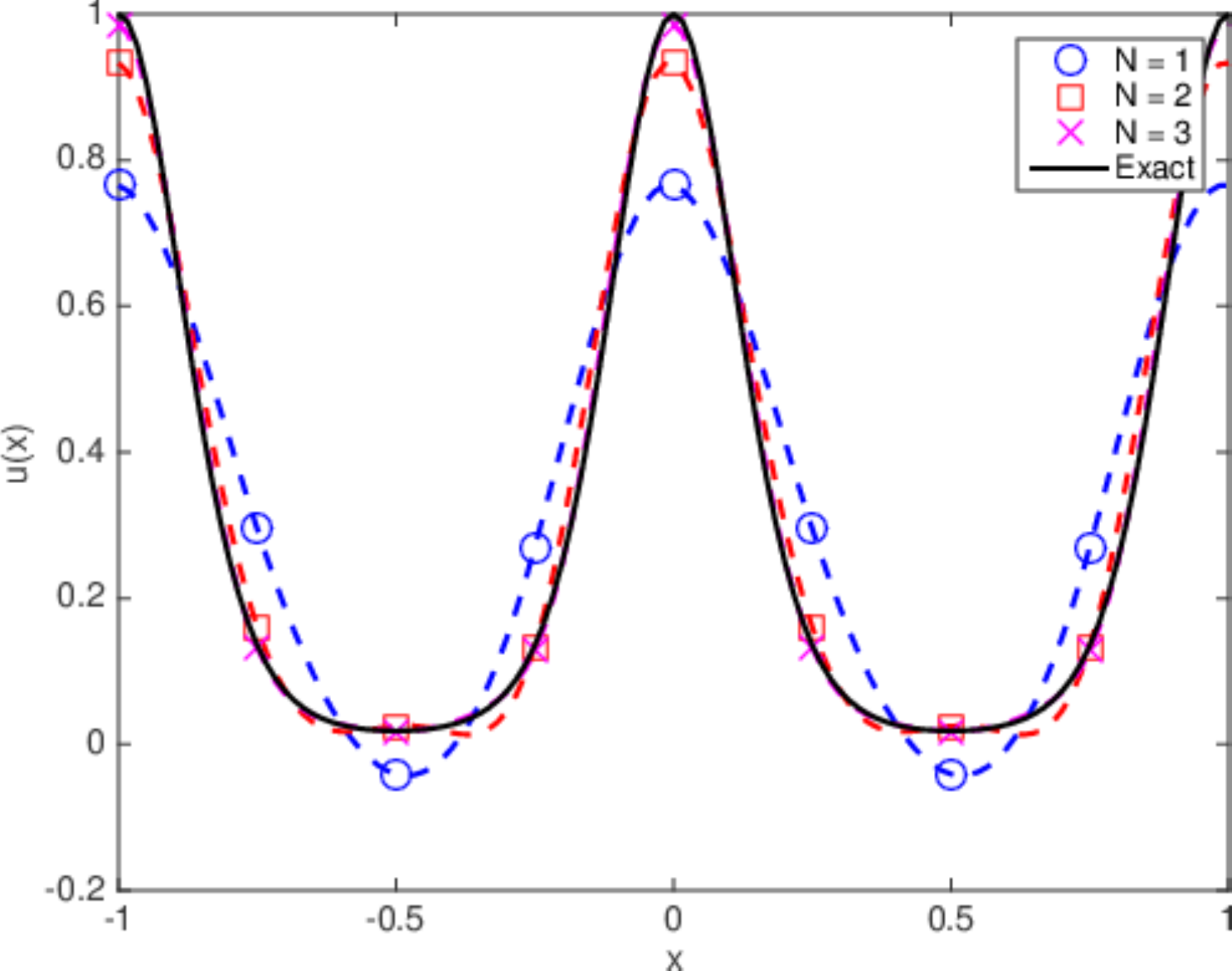}}
\caption{Advection of a periodic Gaussian pulse by Hermite schemes with various CFL constants $C$ on a grid of $K = 8$ cells.}
\label{fig:gauss1}
\end{figure}

Growth of the error in time is often described in terms of dispersive and/or dissipation mechanisms intrinsic to numerical methods \cite{roberts1966convective}.  This may be illustrated by advecting an under-resolved function over several periods; the effect of numerical diffusion will be to smooth the profile out as time increases.  Figures~\ref{fig:gauss1} shows advection of the periodic Gaussian pulse initial condition $e^{-4 \sin(\pi x)^2}$ over 5 periods for orders of approximation $N = 1,2,3$ and a grid of 8 nodes.  %grids consisting of $8$ and $16$ nodes.  

As expected, all methods display diffusive behavior at low orders of approximation and low values of $C$, which is improved as $N$ and $C$ increase.  The Dual and Upwind Hermite methods appear to be the least diffusive, though the difference between each method is small at higher orders of approximation.  On a coarse $K=8$ mesh, the Central Hermite scheme behaves particularly poorly, displaying both spurious oscillations and significant numerical diffusion for all $C$.  

%\begin{figure}
%\centering
%\subfloat[Dual, $C = .1$]{\includegraphics[width=.3\textwidth]{images/pulse_dual_K16_cflopt1.pdf}}
%\hspace{1em}
%\subfloat[Dual, $C = .5$]{\includegraphics[width=.3\textwidth]{images/pulse_dual_K16_cflopt2.pdf}}
%\hspace{1em}
%\subfloat[Dual, $C = .9$]{\includegraphics[width=.3\textwidth]{images/pulse_dual_K16_cflopt3.pdf}}
%
%\subfloat[Virtual, $C = .1$]{\includegraphics[width=.3\textwidth]{images/pulse_virtual_K16_cflopt1.pdf}}
%\hspace{1em}
%\subfloat[Virtual, $C = .5$]{\includegraphics[width=.3\textwidth]{images/pulse_virtual_K16_cflopt2.pdf}}
%\hspace{1em}
%\subfloat[Virtual, $C = .9$]{\includegraphics[width=.3\textwidth]{images/pulse_virtual_K16_cflopt3.pdf}}
%
%\subfloat[Central, $C = .1$]{\includegraphics[width=.3\textwidth]{images/pulse_central_K16_cflopt1.pdf}}
%\hspace{1em}
%\subfloat[Central, $C = .5$]{\includegraphics[width=.3\textwidth]{images/pulse_central_K16_cflopt2.pdf}}
%\hspace{1em}
%\subfloat[Central, $C = .9$]{\includegraphics[width=.3\textwidth]{images/pulse_central_K16_cflopt3.pdf}}
%
%\subfloat[Upwind, $C = .1$, $K=16$]{\includegraphics[width=.3\textwidth]{images/pulse_upwind_K16_cflopt1.pdf}}
%\hspace{1em}
%\subfloat[Upwind, $C = .5$, $K=16$]{\includegraphics[width=.3\textwidth]{images/pulse_upwind_K16_cflopt2.pdf}}
%\hspace{1em}
%\subfloat[Upwind, $C = .9$, $K=16$]{\includegraphics[width=.3\textwidth]{images/pulse_upwind_K16_cflopt3.pdf}}
%
%\caption{Advection of a periodic Gaussian pulse by Hermite schemes for a $K = 16$  grid.}
%\label{fig:gauss2}
%\end{figure}

Increasing to a finer mesh $K=16$, the Central Hermite scheme behaves comparably to the Virtual Hermite scheme.  Qualitatively, the behavior of the Central Hermite scheme for $K=16$ resembles that of the Dual Hermite scheme for $K=8$.  The $L^2$ errors for advection of a Gaussian, while not identical, are very close --- for $N=2$, the Central Hermite scheme with $K=16$ results in an $L^2$ error of $.0841903$, while the Dual Hermite scheme with $K=8$ results in an $L^2$ error of $.0839913$.  This is expected since, for Central with $K=16$ and Dual with $K=8$, the timestep restrictions and interpolation intervals are identical.  

\subsection{Spectra and dispersion/dissipation relations}
\label{sec:spectra}

Numerical experiments confirm the high order convergence of each method; however, the qualitative behavior of each method in convecting an under-resolved solution varies significantly.  We seek to further analyze this behavior by computing the spectra of the update matrix and dispersion/dissipation relations for each Hermite method.

We define the update operator $\mb{S}$ such that, for solution degrees of freedom $\mb{U}^n$ at time $t_n$, the application of $\mb{S}$ evolves the solution at time $t_n + dt$ 
\[
\mb{U}^{n+1} = \mb{S}\mb{U}^n.
\]
Both the Dual and Central Hermite method march forward by $dt = C h_x/c$ over a single timestep (the Dual Hermite method defines $dt = Ch_x/(2c)$, but takes timesteps on both primary and dual grids).  Since the Virtual Hermite method takes the timestep to be $dt = Ch_x/(2c)$, we analyze $\mb{S}^2$ for the Virtual grid (corresponding to taking two timesteps instead of one) in order to normalize how far in time the update operator evolves the solution.  

\begin{figure}
\centering
\subfloat[Dual, $C = .1$]{\includegraphics[width=.3\textwidth]{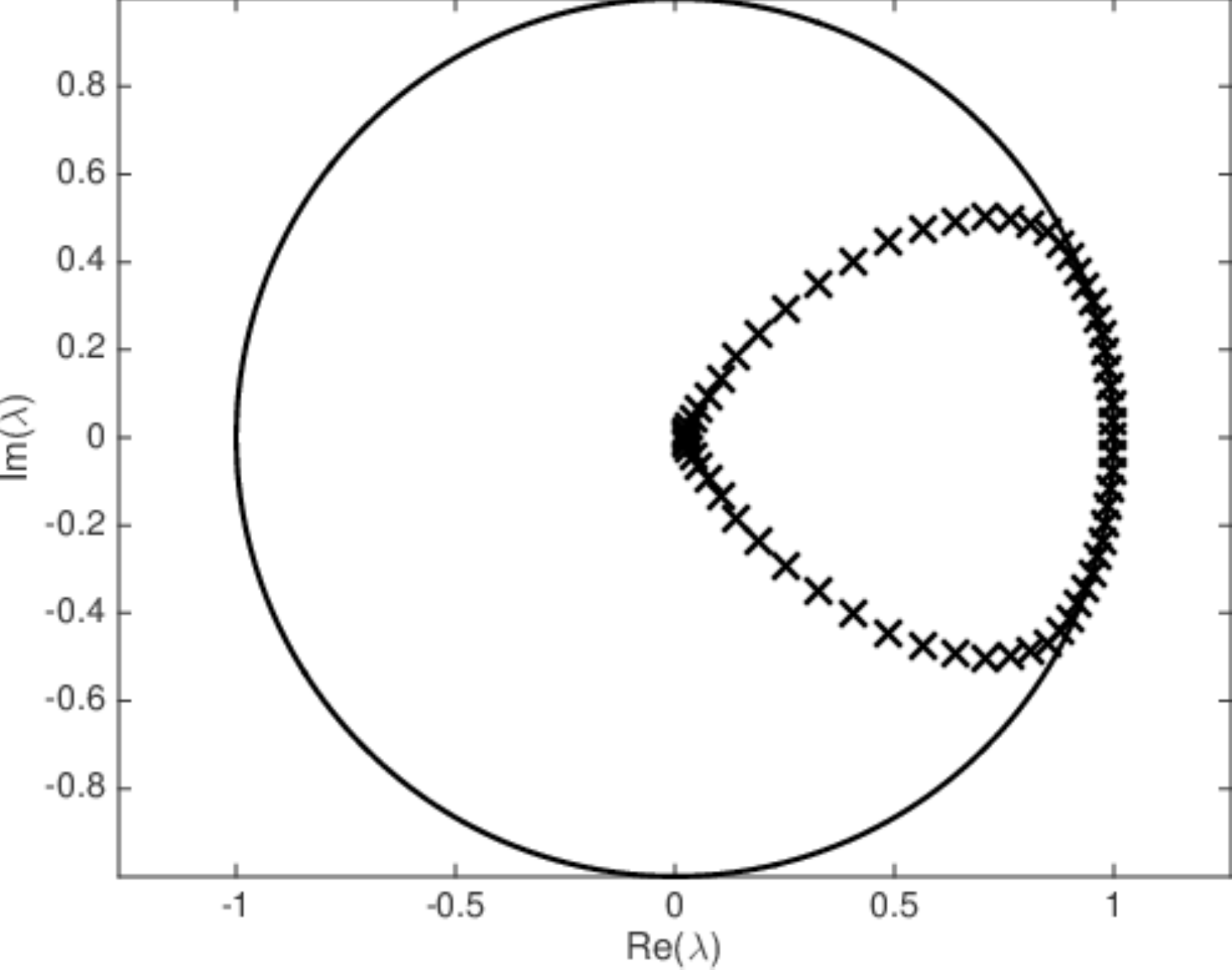}}
\hspace{1em}
\subfloat[Dual, $C = .5$]{\includegraphics[width=.3\textwidth]{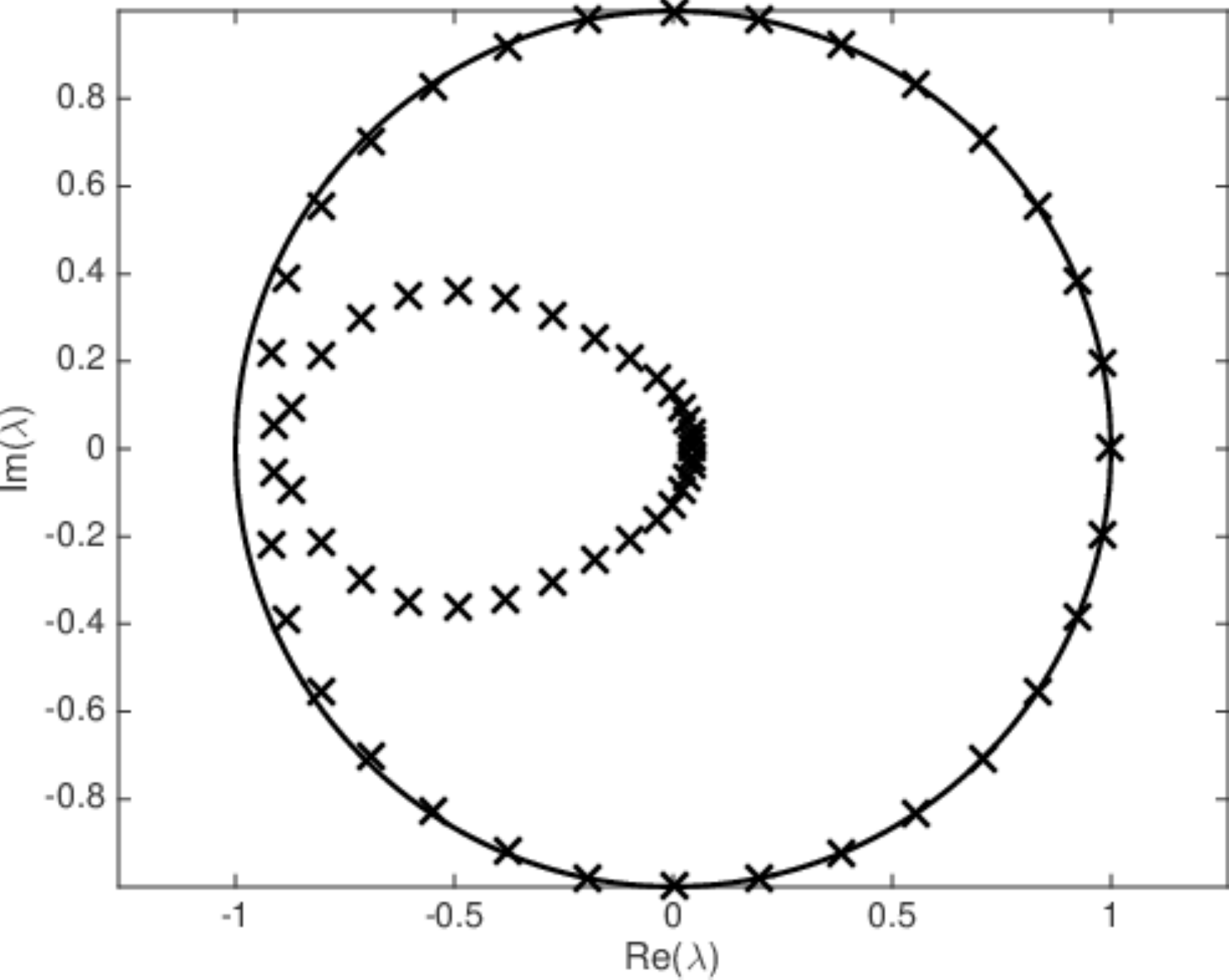}}
\hspace{1em}
\subfloat[Dual, $C = .9$]{\includegraphics[width=.3\textwidth]{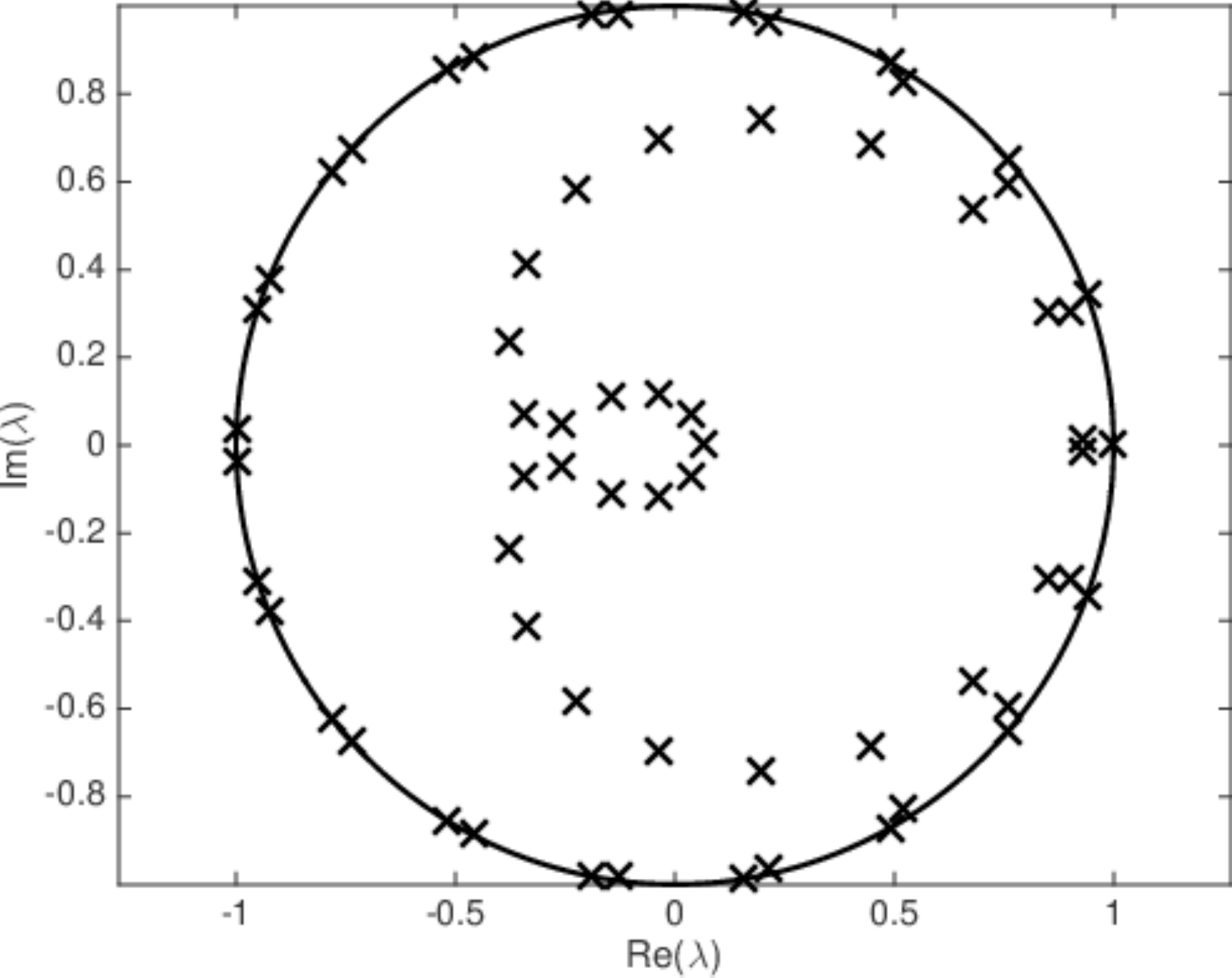}}

\subfloat[Virtual, $C = .1$]{\includegraphics[width=.3\textwidth]{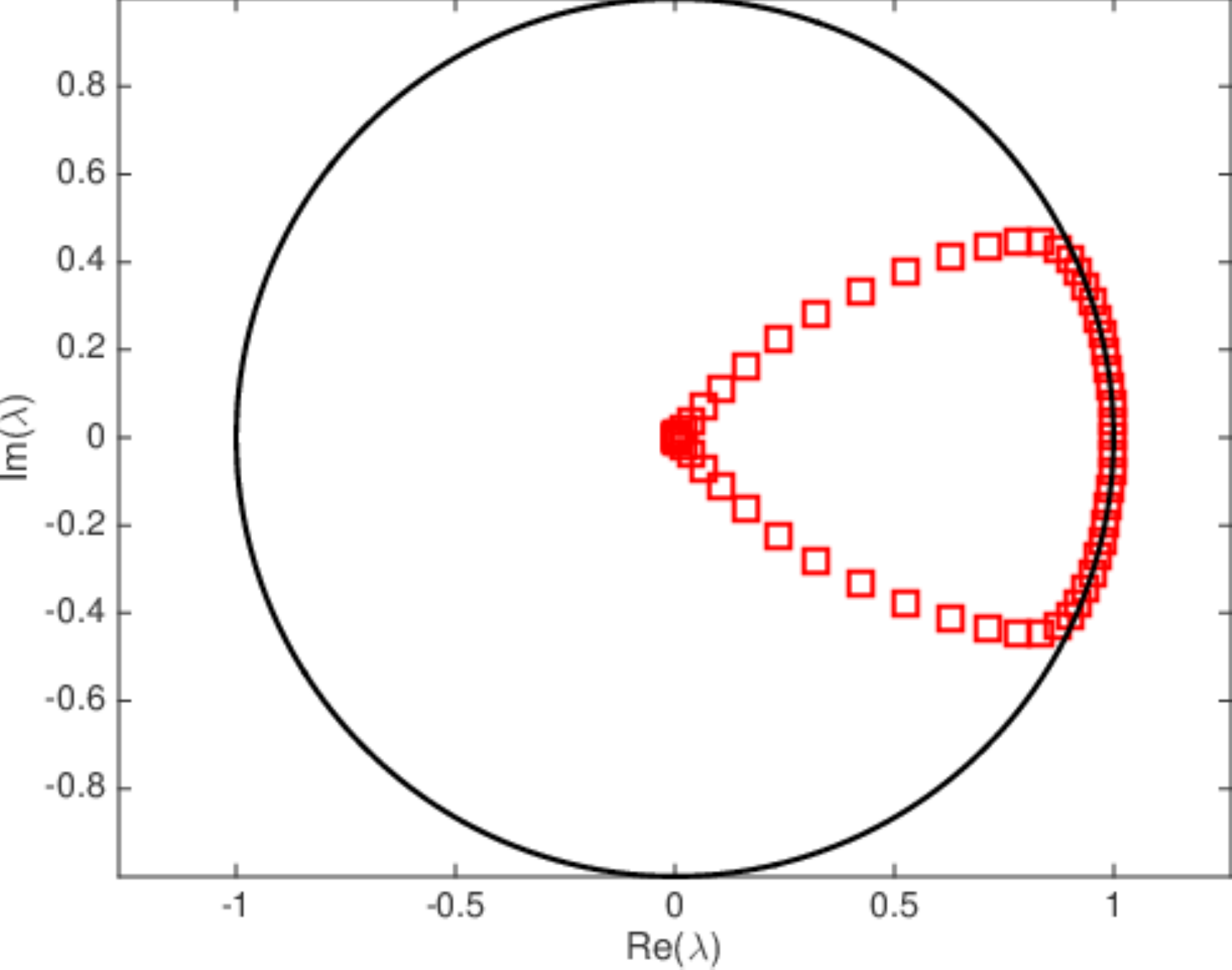}}
\hspace{1em}
\subfloat[Virtual, $C = .5$]{\includegraphics[width=.3\textwidth]{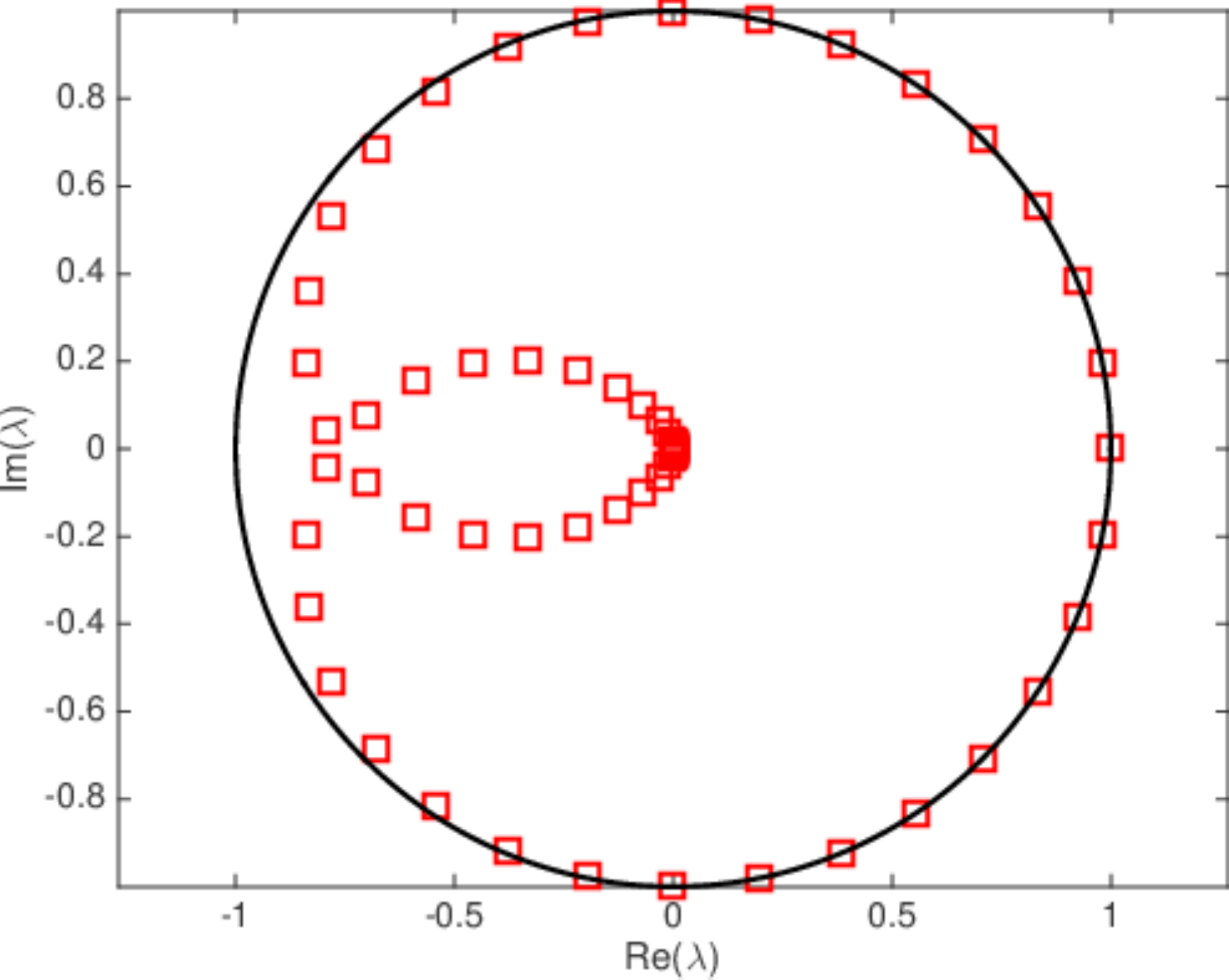}}
\hspace{1em}
\subfloat[Virtual, $C = .9$]{\includegraphics[width=.3\textwidth]{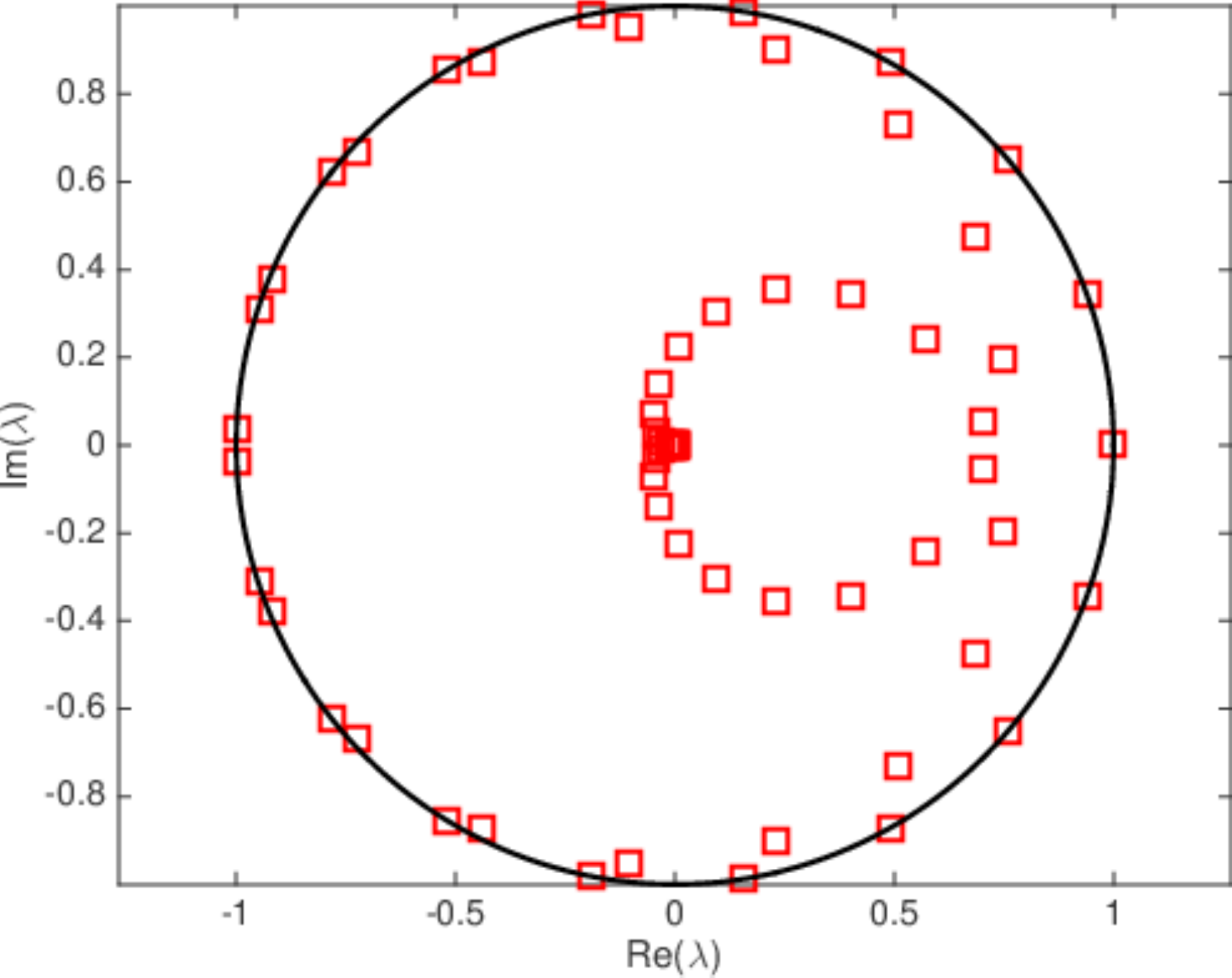}}

\subfloat[Central, $C = .1$]{\includegraphics[width=.3\textwidth]{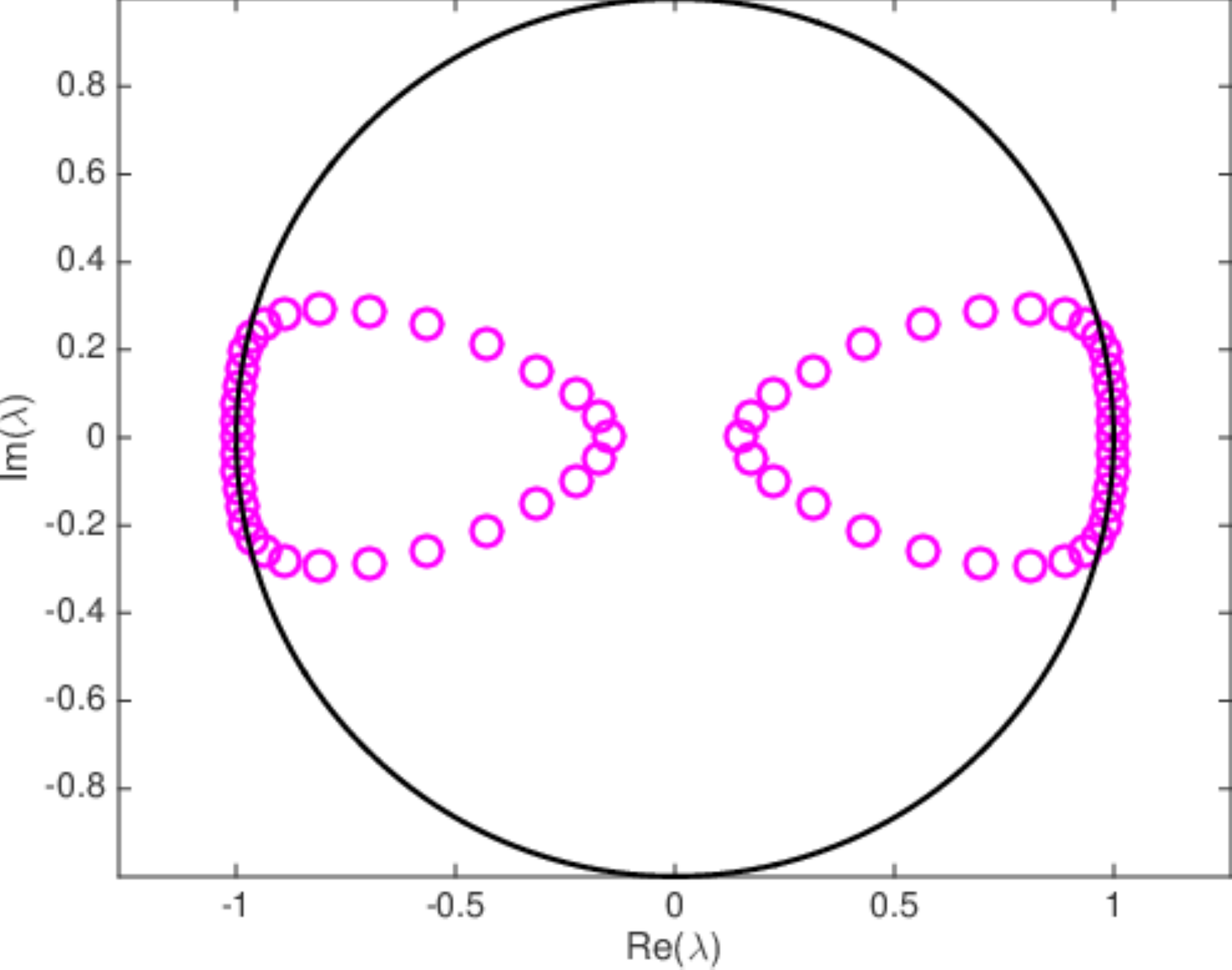}}
\hspace{1em}
\subfloat[Central, $C = .5$]{\includegraphics[width=.3\textwidth]{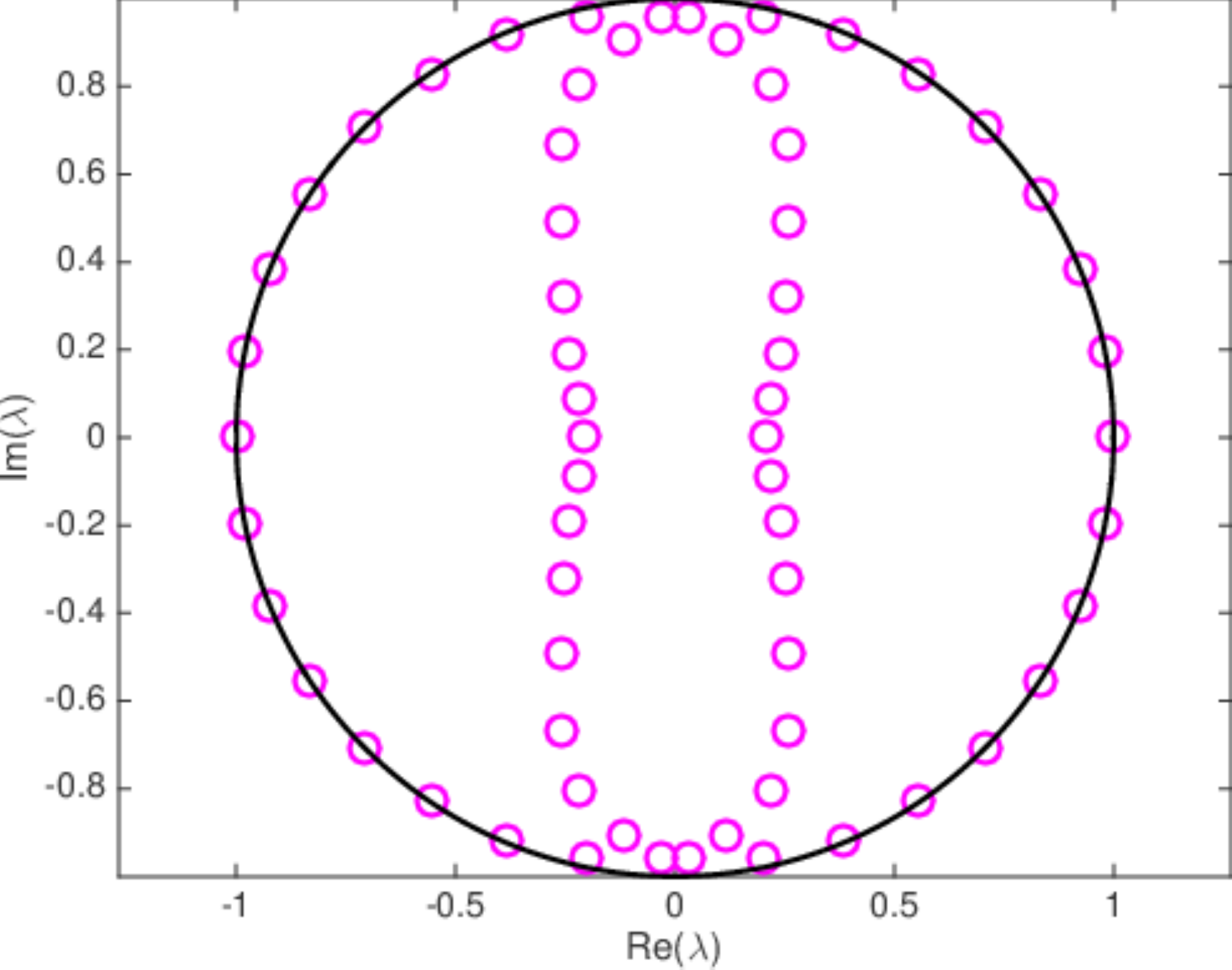}}
\hspace{1em}
\subfloat[Central, $C = .9$]{\includegraphics[width=.3\textwidth]{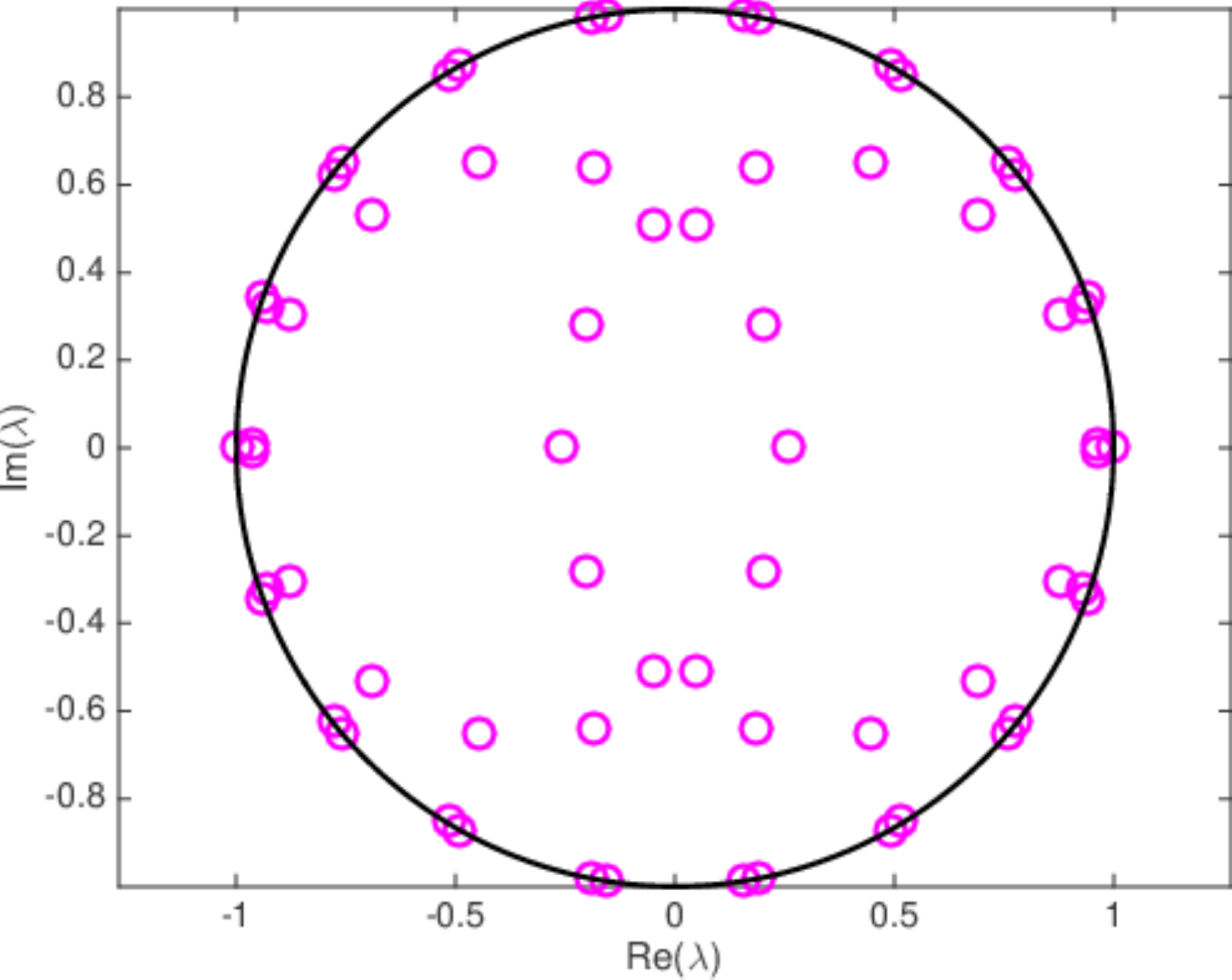}}

\subfloat[Upwind, $C = .1$]{\includegraphics[width=.3\textwidth]{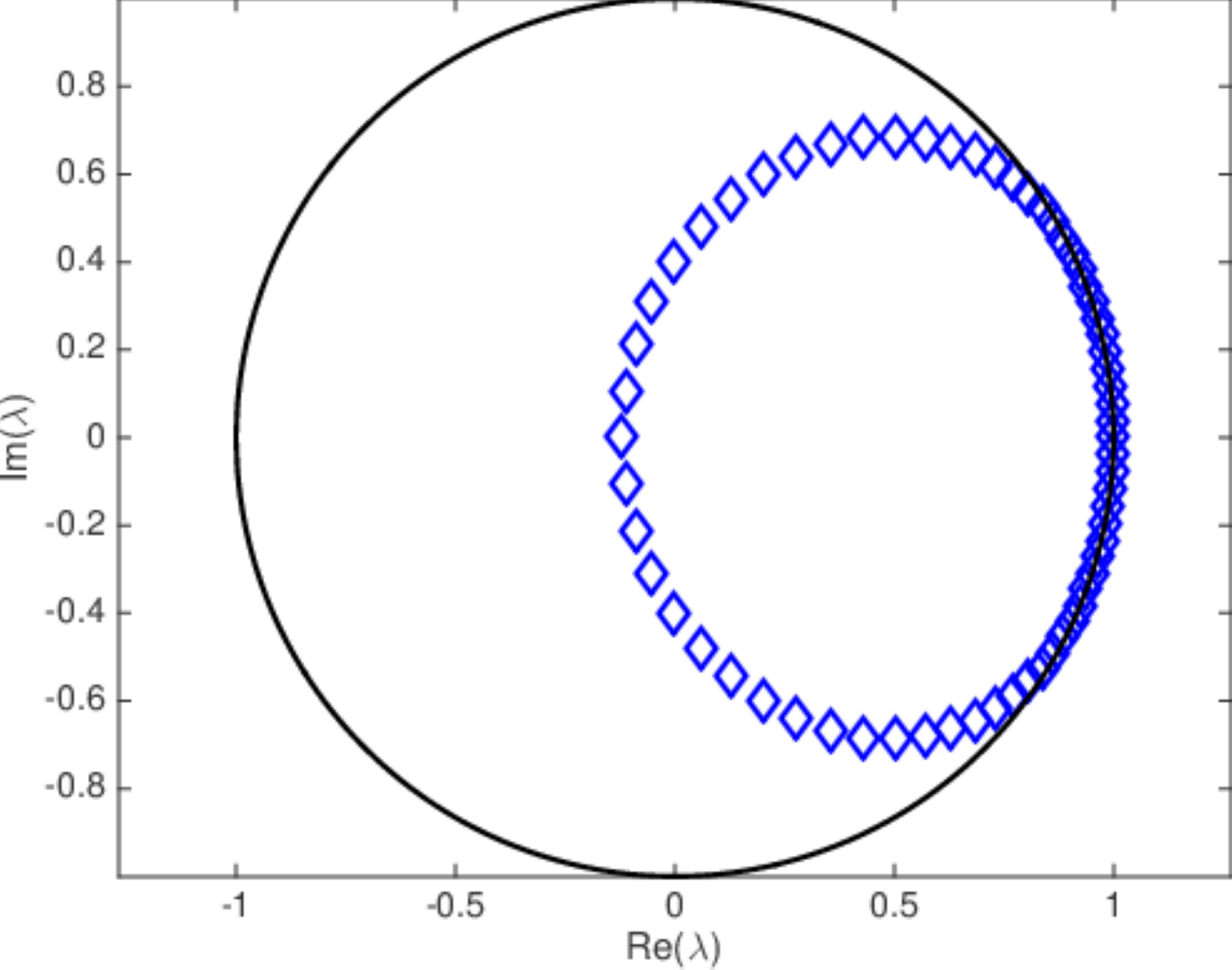}}
\hspace{1em}
\subfloat[Upwind, $C = .5$]{\includegraphics[width=.3\textwidth]{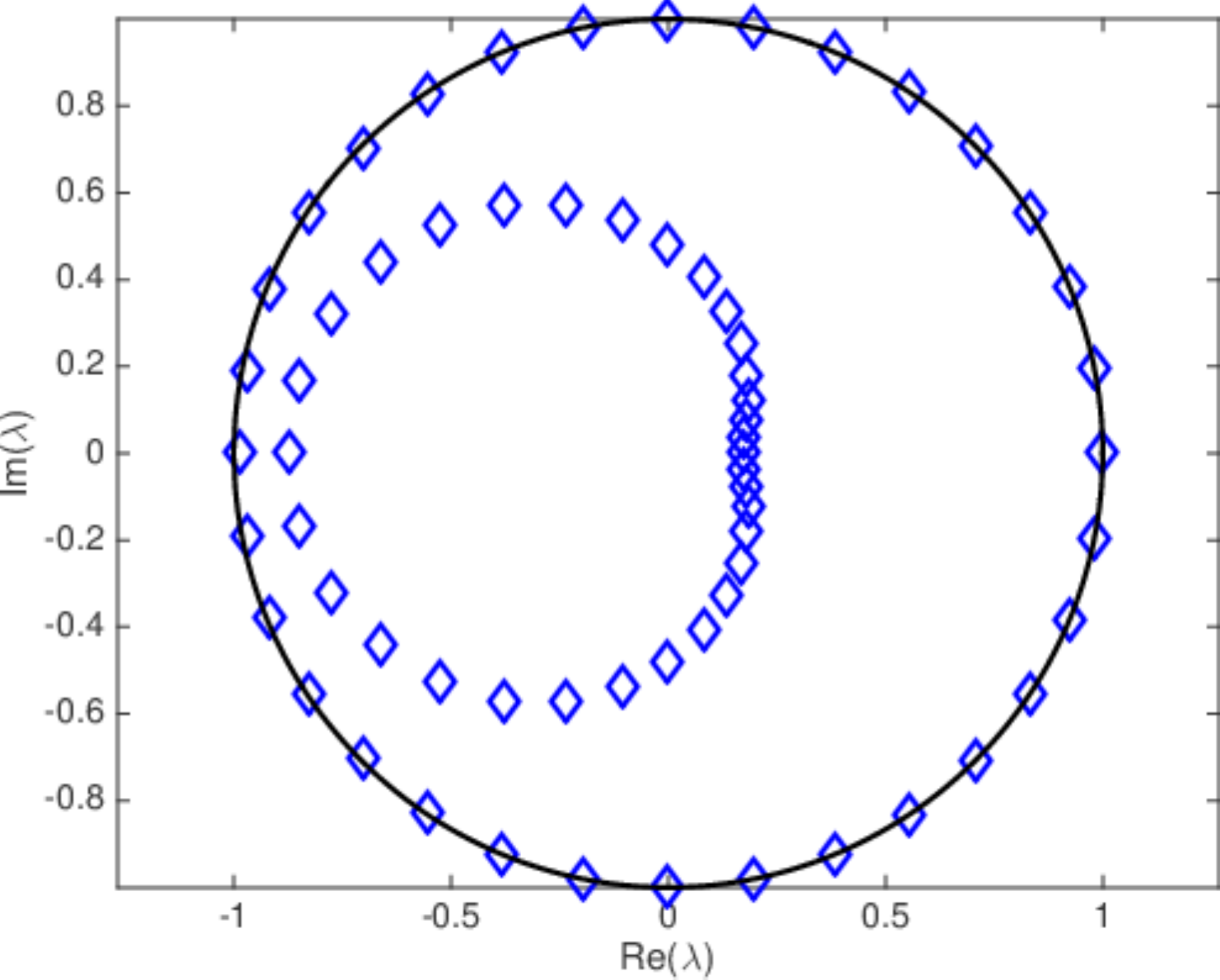}}
\hspace{1em}
\subfloat[Upwind, $C = .9$]{\includegraphics[width=.3\textwidth]{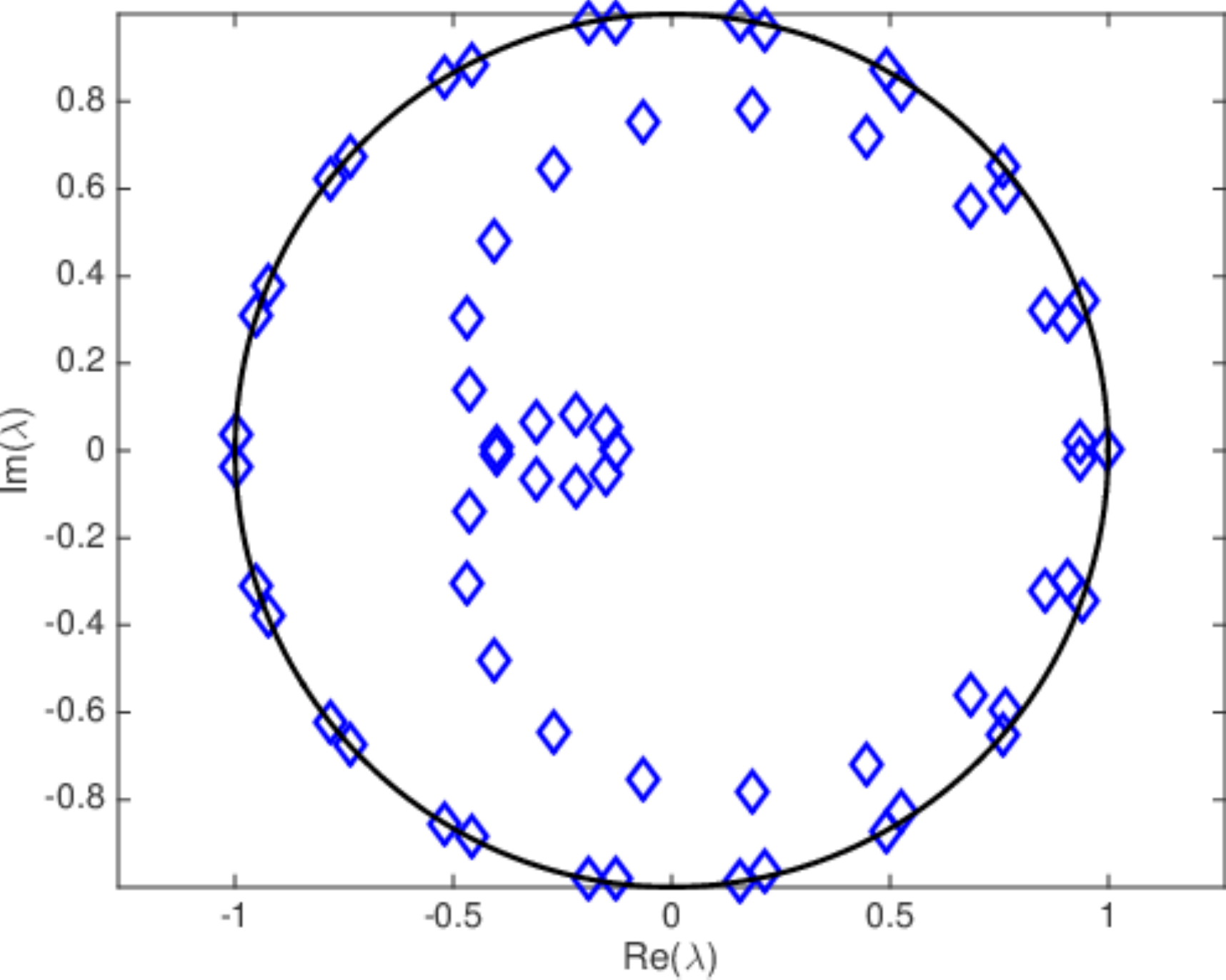}}
\caption{Spectra of the update matrix for each Hermite scheme at various CFL constants $C$.  The order of approximation and grid size are fixed to be $N = 3$ and $K = 16$, respectively.}
\label{fig:spectra}
\end{figure}

For the Dual, Central, and Upwind Hermite methods, setting $C=1$ exactly (with respect to machine precision) results in exact evolution of the solution, though this is unique to constant coefficient equations.  As a consequence, the update operator $\mb{A}$ becomes exactly equal to a circulant shift matrix.  As a result, the order of convergence of each of these methods with $C=1$ increases to $O(h^{2N+2})$, and coincides with Hermite interpolation estimates given in \cite[Lemma 3.1]{goodrich2006hermite}.  This same exact evolution property does not hold for the Virtual Hermite method.  

The spectra in Figure~\ref{fig:spectra} suggest that the Virtual and Dual Hermite methods should behave similarly, as the eigenvalues of $\mb{A}$ are distributed similarly for both methods.  Eigenvalues which lie on the unit circle are typically of the form $e^{i \omega}$, and are related to the non-dissipative propagation of modes of the form $e^{i \omega (x-ct)}$.  For example, for $C = .1$, the eigenvalues fall closest to the unit circle around the point $(1,0)$, corresponding to the non-dissipative propagation of modes with small $\omega$ (low frequency modes).  The remaining spectra lie within the unit circle, indicating dissipation of under-resolved modes.  In constrast, the spectra for the Central Hermite method clusters not only around $(1,0)$ but also around $(-1,0)$, suggesting that under-resolved high frequency modes may be propagated without dissipation.  These spurious modes may explain the behavior of the Central Hermite method observed in Figure~\ref{fig:gauss1}, where propagation of a Gaussian on a coarse grid resulted in ``spurious'' oscillatory behavior which remained over several periods of advection.  

%Finally, numerical experiments with the Upwind stencil indicated low numerical dissipation and optimal convergence rates for the advection equation in one space dimension.  Additionally, the magnitude of the error is approximately the same as the magnitude of the error for the Virtual and Dual Hermite methods, and Figure~\ref{fig:upwind_spectra} indicates that the spectra of the update matrix for the Upwind Hermite method is well-behaved and does not contain spurious modes.  

A study of the dispersion and dissipation error for the Dual Hermite method was reported in \cite{jang2012analysis} using a modified equation and Bloch wave analysis in one dimension, which we adapt and apply to the Hermite methods introduced in this work.  Dispersion and dissipation properties of Hermite methods depend mainly on the properties of one-dimensional Hermite interpolation matrix $\mb{H}$ 
\[
\mb{H} = \LRs{\mb{H}_L, \mb{H}_C, \mb{H}_R}
\]
where $\mb{H}_C,\mb{H}_L,\mb{H}_R$ act on Hermite data associated with a given node and it's left/right neighbors to produce a reconstruction.  

%The evolution matrix $\mb{S} \in \mbb{R}^{K(N+1)\times K(N+1)}$ is the system which advances the solution degrees of freedom $\mb{u}^n$ forward a single timestep to $\mb{u}^{n+1}$. 
For a periodic grid, $\mb{S} \in \mbb{R}^{K(N+1)\times K(N+1)}$  is a block tridiagonal matrix
\[
\mb{S} = \begin{bmatrix}
\mb{S}_L & \mb{S}_R & &\ldots & \mb{S}_L\\
\mb{S}_L & \mb{S}_C & \mb{S}_R  \\
& \mb{S}_L & \mb{S}_C & \ddots \\
& & \ddots & \ddots &  \\
\mb{S}_R & & & \mb{S}_L & \mb{S}_C
\end{bmatrix}
\]
where $\mb{S}_L = \mb{T}\mb{H}_L$, and similarly for $\mb{S}_C, \mb{S}_R$. For the Central and Upwind Hermite methods, $\mb{S}_C$ and $\mb{S}_R$ are zero, respectively.   

We perform a fully discrete Bloch analysis to examine dispersive and dissipative properties of each Hermite method.  This is done by representing the wave solution $e^{ik(x-ct)}$ in the Hermite basis, and noting that the solution is shifted in both space and time by scaling with a complex exponential %The exact evolution of the solution from timestep $t$ to timestep $t+dt$ is
\[
u(x,t+dt) = u(x,t)e^{-ikcdt}, \qquad u(x+h,t) = u(x,t)e^{ikh}.  
\]
Assuming a uniform grid spacing $h$, the discrete evolution of the interpolated exact solution at a node $x_m$ from time $t_n$ to $t_n + dt$ is then given by 
\[
%\LRp{e^{-ikh}\mb{S}_L + \mb{S}_C + e^{ikh}\mb{S}_R}\mb{u}^n_m = e^{-i \omega_h dt}\mb{u}^{n+1}_m.  
\LRp{e^{-ikh}\mb{S}_L + \mb{S}_C + e^{ikh}\mb{S}_R}\mb{u}^n_m = \lambda_h\mb{u}^{n+1}_m.  
\]
Since the timestep restrictions for the Dual and Virtual Hermite methods are $dt = Ch/(2c)$ as opposed to $dt = Ch/c$, the dispersion relations are measured over two timesteps.  For the Dual Hermite method, this implies that $\mb{S}$ captures the evolution of the solution from the primal to dual grid, then back to the primal grid.  For the Virtual Hermite method, this requires taking two timesteps and substituting for $\mb{S}$ the matrix $\mb{S}^2$, which is  block pentadiagonal.  

An eigenvalue problem may be solved for the discrete Floquet multiplier $\lambda_h$, whose real and imaginary parts correspond to numerical dispersion and dissipation, respectively.  For each Hermite method, we measure the relative error between the discrete and the true Floquet multiplier 
\[
E_{kh} = \frac{\LRb{\lambda_h - e^{-ikcdt}}}{\LRb{e^{-ikcdt}}}
\]
over a single timestep as a function of $kh$, the order of approximation $N$ and the CFL constant $C = .1, .5, .9$.  %Numerical experiments indicate that the real and imaginary parts of this error are of roughly equal magnitude.  

\begin{figure}
\centering
\subfloat[$C = .1$]{\includegraphics[width=.42\textwidth]{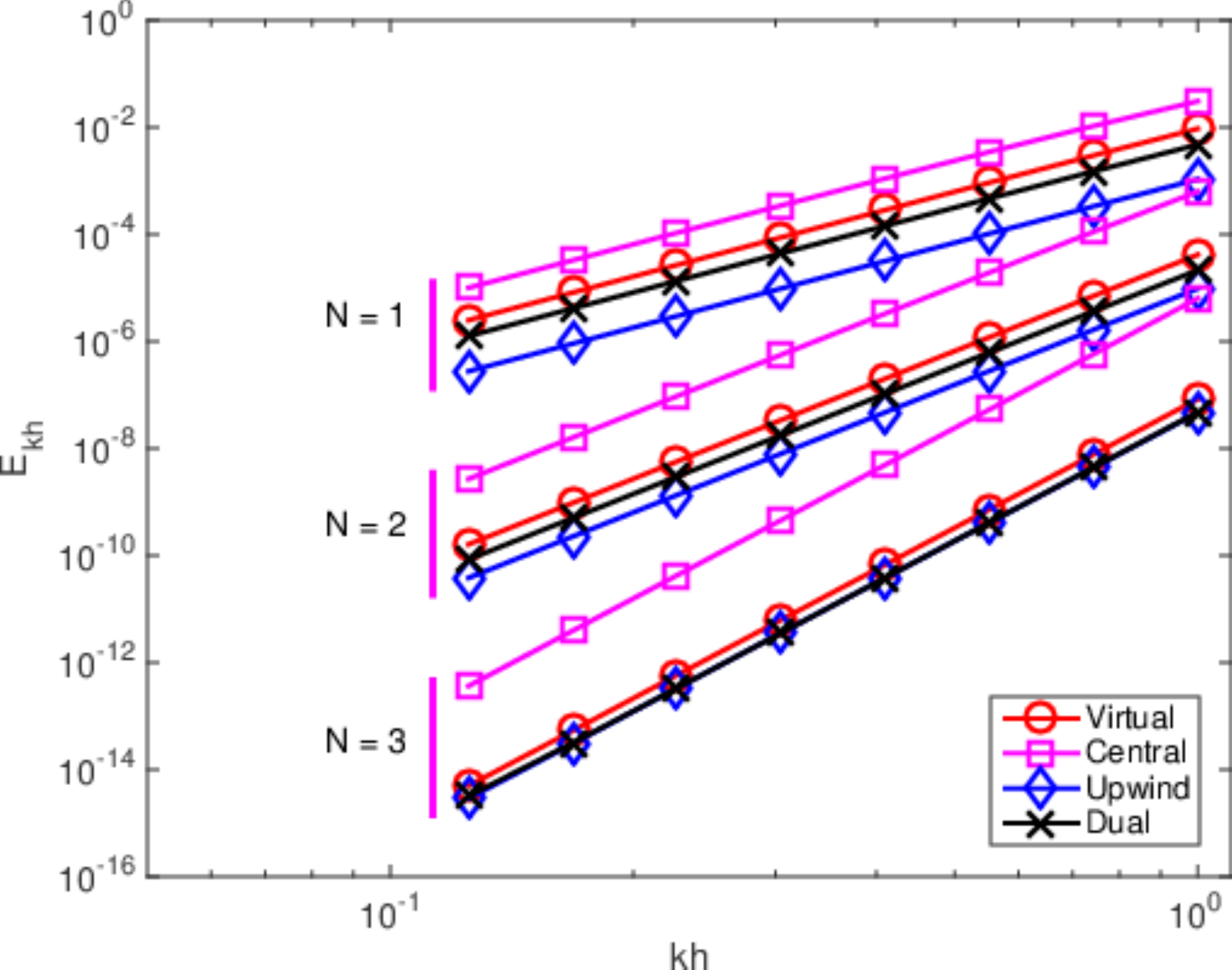}}
\hspace{1em}
\subfloat[$C = .9$]{\includegraphics[width=.42\textwidth]{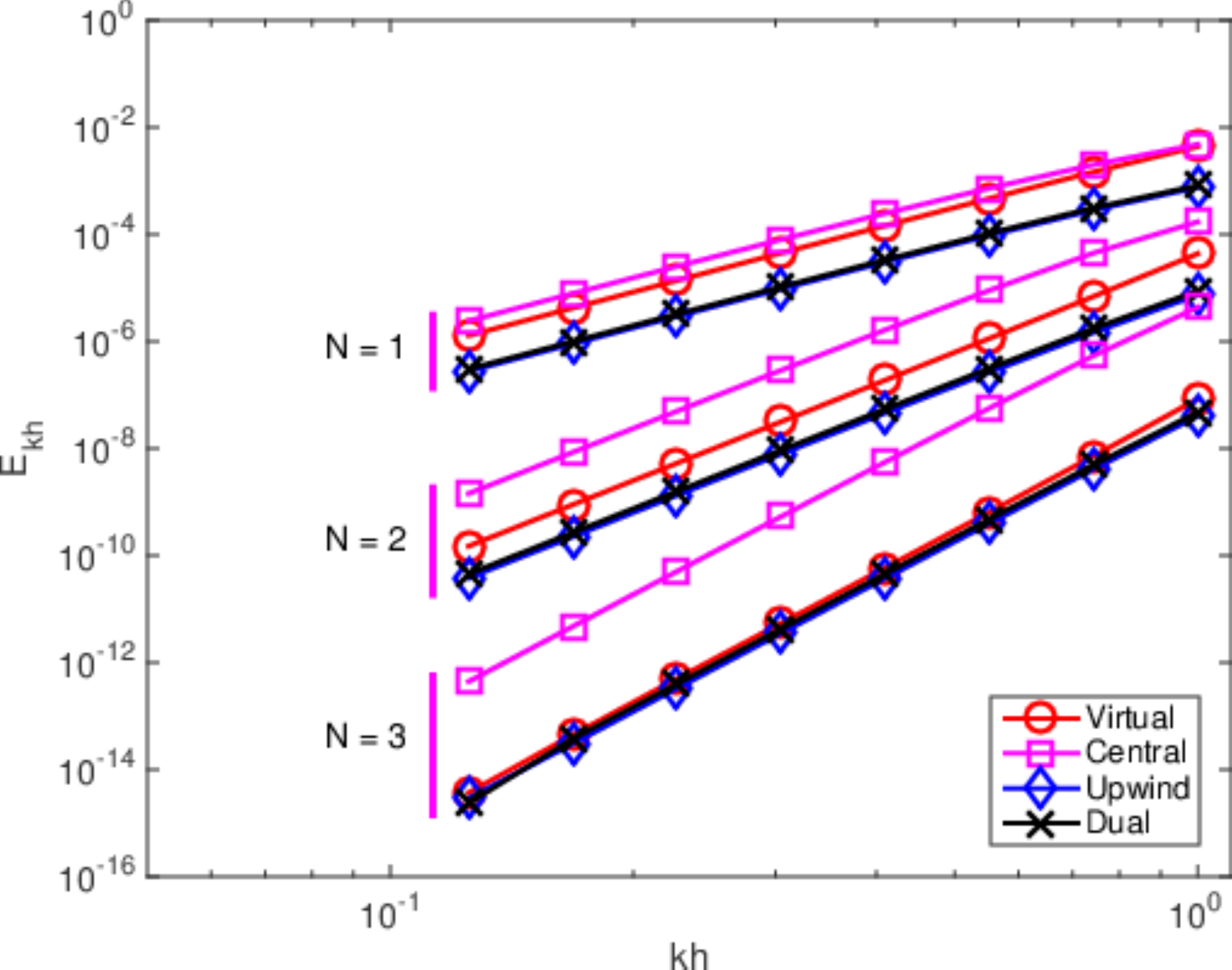}}
\caption{Dispersive and dissipative errors ${\LRb{\lambda_h - e^{-ikcdt}}}/{\LRb{e^{-ikcdt}}}$ for each Hermite scheme, with $N = 1,2, 3$.  The computed errors are observed to behave as $O\LRp{{kh}/{c}}^{2N+2}$.} %$\LRb{\omega - \omega_h}$.}
\label{fig:disp}
\end{figure}
Figure~\ref{fig:disp} shows the error $E_{kh}$ over a range of $kh$.  The Central Hermite method shows the largest errors, while Upwind Hermite shows the smallest errors.  The error for the Virtual and Dual methods lie in-between, with the Dual Hermite displaying smaller errors at low $N$.  Smaller values of the CFL constant $C$ increase the dispersion and dissipation error, though the effect is less noticable as $N$ increases.  

For each method, the error $E_{kh}$ is observed to follow
\begin{align*}
\frac{\LRb{{\lambda_h - e^{-ikcdt}}}}{\LRb{{e^{-ikcdt}}}} &\approx C_N \LRp{\frac{kh}{c}}^{2N+2}
\end{align*}
where $N$ is the number of degrees of freedom per node, and the underlying Hermite approximation space is of degree $2N+2$.\footnote{If we seek instead the error in the discrete and exact wavenumbers $\LRb{\omega - \omega_h}$, we recover convergence rates of $O(h^{2N+1})$ and $O(h^{2N+3})$ for the real and imaginary parts, respectively.  In comparison, DG with co-volume filtering achieves rates of $O(h^{2N+2})$ and $O(h^{2N+3})$ for the real and imaginary parts, respectively \cite{warburton2008taming}, while Galerkin methods result in rates of either $2N+1$ or $2N+3$ (for $N$ even or odd) under a degree $N$ approximation space \cite{ainsworth2014dispersive}.}.  

\subsubsection{Optimizing dispersive and dissipative errors}
Finally, motivated by Dispersion Relation Preserving (DRP) finite difference schemes \cite{tam1993dispersion}, dispersive and dissipative errors may be improved through optimization of entries of the interpolation matrix.  As mentioned in Section~\ref{sec:virtual}, the Virtual Hermite method produces a degree $2N+1$ reconstruction using degree $N$ data from three nodes, though there is sufficient data to define a higher $3N+2$ degree reconstruction.  We define the interpolation matrix $\tilde{\mb{H}}\in \mbb{R}^{(3N+3)\times (3N+3)}$
\[
\tilde{\mb{H}} = \begin{bmatrix}
\mb{H}\\
\mb{H}_2
\end{bmatrix},
\]
where $\mb{H}$ is the Virtual Hermite interpolation matrix.  DRP schemes enforce a fixed order of approximation for a given finite difference stencil, while using additional degrees of freedom to optimize the dispersion relation.  Similarly, fixing the first $2N+2$ rows of $\tilde{\mb{H}}$, the reconstruction implied by $\tilde{\mb{H}}$ is enforced to match that of the Virtual Hermite reconstruction for the first $2N+2$ coefficients, while entries of the matrix $\mb{H}_2$ (which determine higher order coefficients) are used to minimize dispersion and dissipation errors.  The entries of $\mb{H}$ depend on the ratio between $L_m-\tilde{x}$ (or $R_m-\tilde{x}$) and as the grid spacing $h$.  Since this ratio is constant as a function of $h$, $\mb{H}_2$ does not change drastically as a grid is refined.  However, the optimization does appear to be sensitive to the value of $C$.  

\begin{figure}
\subfloat[Spectra of $\mb{S}$]{\includegraphics[width=.335\textwidth]{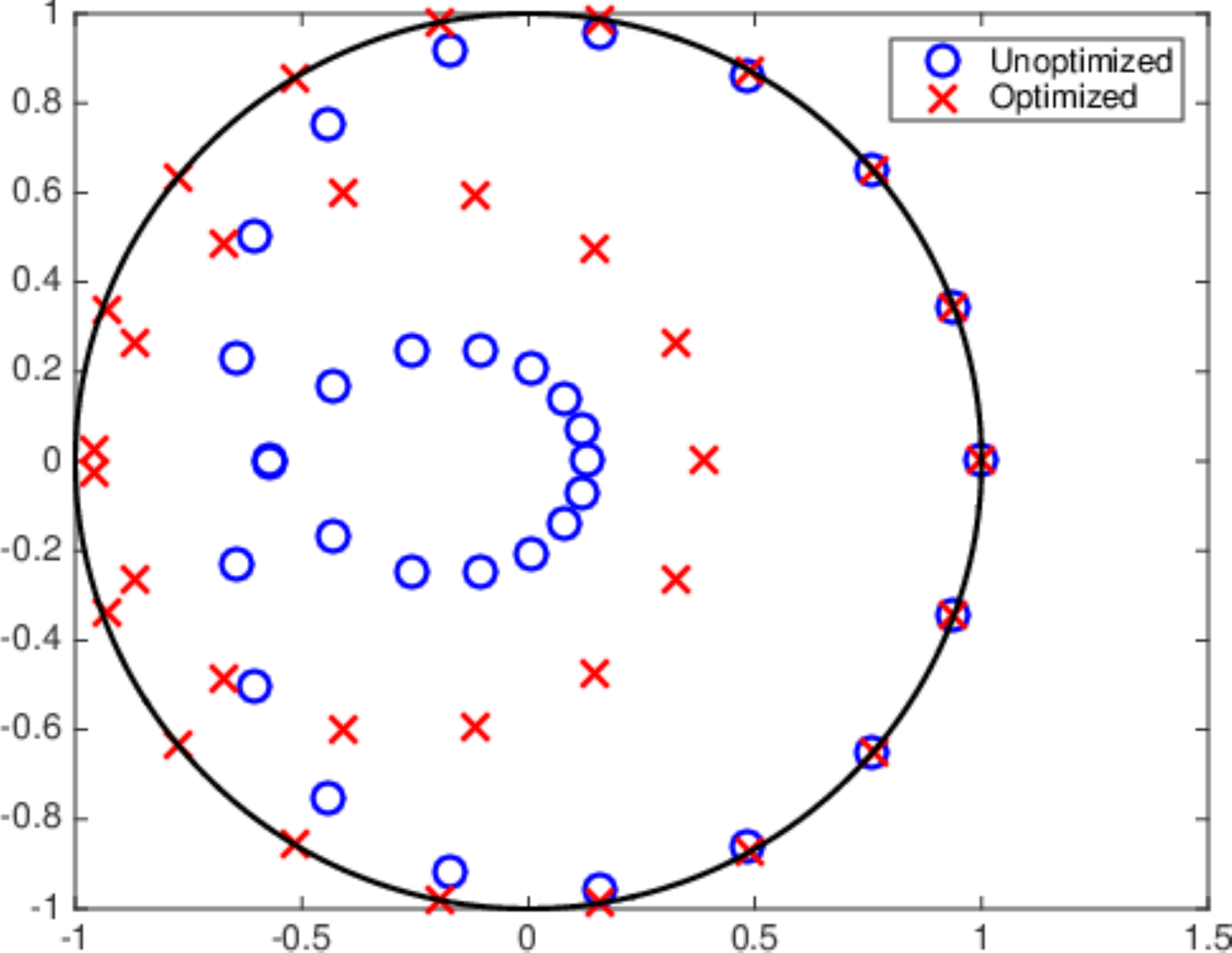}}
\subfloat[$E_{kh}$]{\includegraphics[width=.315\textwidth]{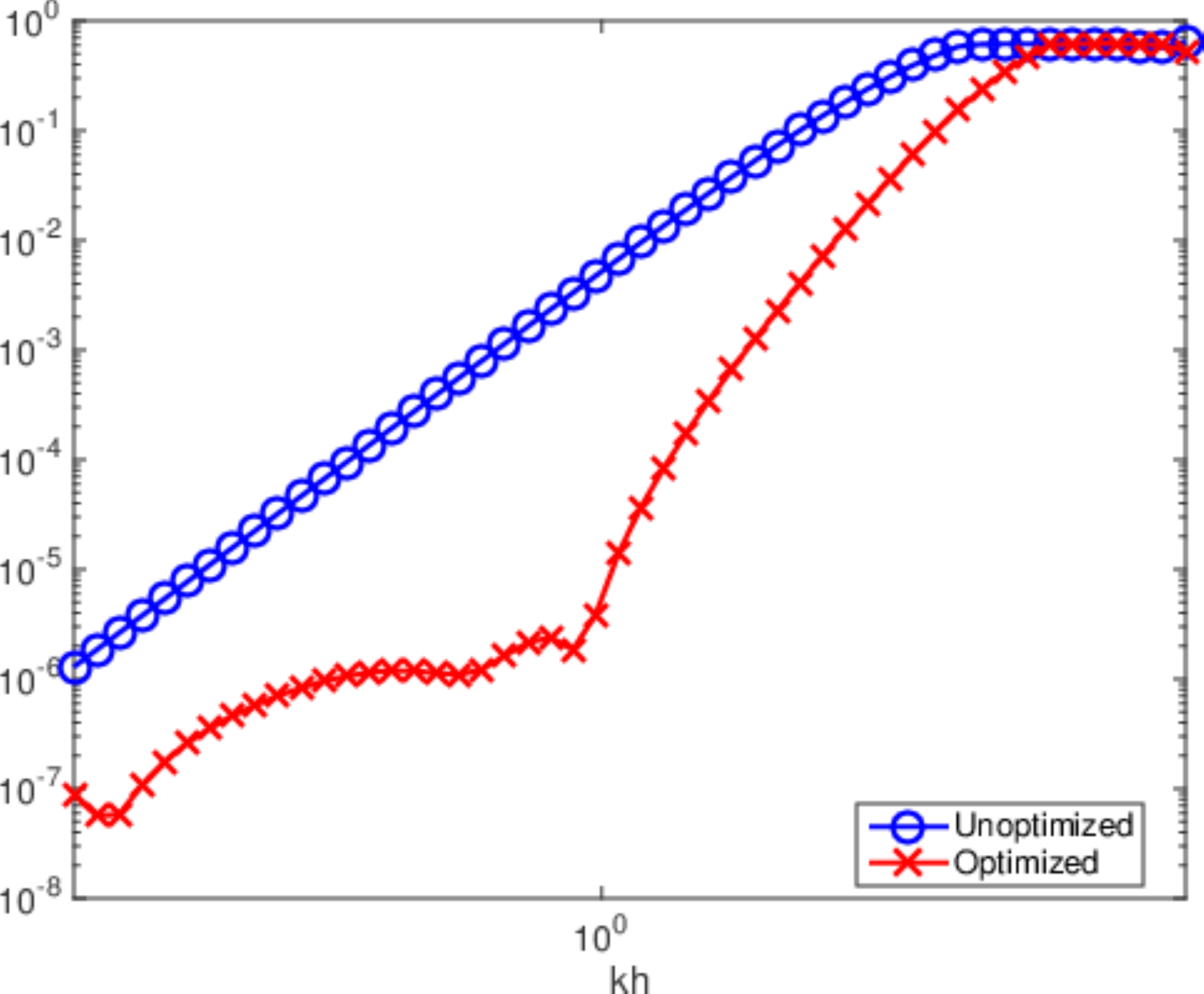}}
\subfloat[Convection over 5 periods]{\includegraphics[width=.335\textwidth]{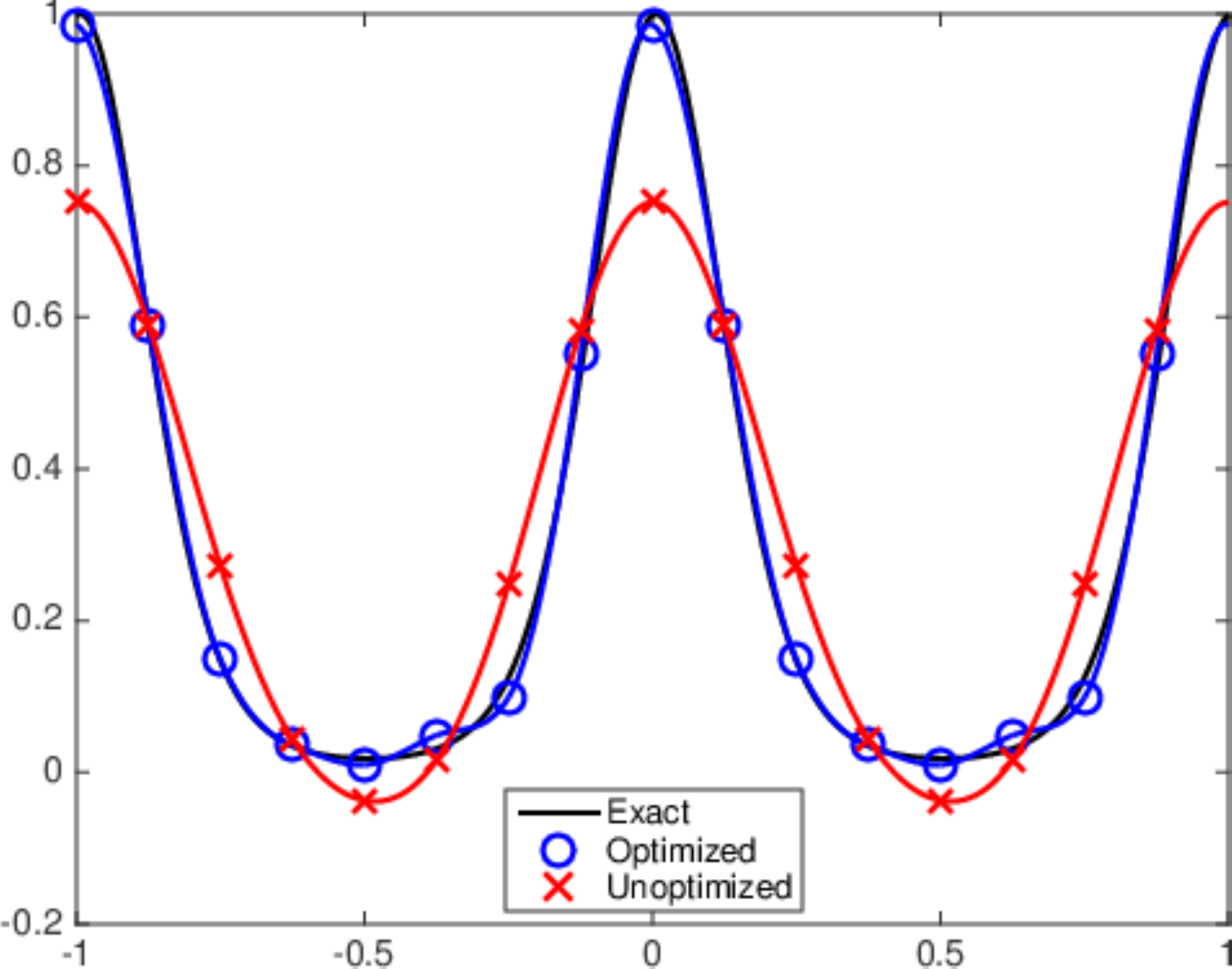}}
\caption{Spectra, dispersion/dissipation errors, and convection of the initial condition $e^{-4\sin(\pi x)^2}$ over 5 periods.  Results are shown for both standard and optimized Virtual Hermite methods with $N=1$, $K=16$, and $C=.9$.  Optimization is done on a coarse $K=8$ mesh.  }
\label{fig:vhopt}
\end{figure}

To demonstrate the effect of optimization, we compare the Virtual Hermite method to an optimized scheme for $N=1$ and CFL constant $C=.9$.  We produce the optimized submatrix $\mb{H}_2$ by minimizing the real and imaginary parts of the relative dispersion error for the advection equation 
\[
\LRp{\frac{{\rm Re}\LRp{{\lambda_h - e^{-ikcdt}}}}{{\rm Re}\LRp{{e^{-ikcdt}}}}}^2 + \LRp{\frac{{\rm Im}\LRp{{\lambda_h - e^{-ikcdt}}}}{{\rm Im}\LRp{{e^{-ikcdt}}}}}^2
\]
with $c=1$ and $K=8$ grid cells.  The same optimized submatrix $\mb{H}_2$ is then used on a finer $K=16$ grid, and computed solutions for the Virtual and optimized Hermite methods are compared for the initial condition $e^{-4\sin(\pi x)^2}$ in Figure~\ref{fig:vhopt}.  The spectra of the update matrix $\mb{S}$ and dispersion/disspation errors $E_{kh}$ are also compared in Figure~\ref{fig:vhopt}.  The dispersion error and spectra are shown to be significantly improved, and numerical results indicate that under-resolved features are convected with greater accuracy.  

Unfortunately, the benefits of such an approach appear to be limited to low orders of approximation.  At higher orders, optimization did not reduce the dispersion and dissipation error $E_{kh}$ significantly compared to the unoptimized scheme.  Additionally, the stability of such an approach is not guaranteed for Hermite methods (compared to DRP schemes, which optimized over symmetric stencils to guarantee stability).  For example, for $N=1$, the spectral radius of the update matrix for the unoptimized scheme was computed to be $\rho(\mb{S})=1$ to machine precision.  For the $N=1$ optimized scheme, $\rho(\mb{S}) \approx 1.0005$, and strict enforcement of $\rho(\mb{S}) \leq 1$ resulted in either non-convergence of the optimization problem or subpar dispersion and dissipation properties.  Further study is required to address these issues.

\section{Extension to two dimensions}
\label{sec:2D}
Each Hermite method may be extended to higher dimensions naturally through a tensor product construction.  In this work, we take the grid $\Omega$ to be the tensor product of one-dimensional grids.  Assuming grid spacings $h_x, h_y$ in the $x$ and $y$ directions, respectively, each point $(x_m, y_m) \in \Omega$ admits the tensor-product expansion 
\[
u_m(x,y) = \sum_{j=0}^N \sum_{k=0}^N \mb{u}_{jk} \LRp{\frac{x-x_m}{h_x}}^j\LRp{\frac{y-y_m}{h_y}}^k.  
\]
For linear autonomous equations, a temporal Taylor series may be used to evolve the solution in time.  Algorithm~\ref{alg:time2} describes this process for the two-dimensional scalar advection equation, using derivative matrices $\mb{D}_x, \mb{D}_y$ for the $x$ and $y$ coordinates, respectively.  Due to the tensor-product nature of the Hermite interpolants in higher dimensions, the Taylor series must be of order $d(2N+1)$ to be exact in $d$-dimensions \cite{goodrich2006hermite}.  Numerical experiments indicate that reducing the degree of the Taylor expansion in time results in a tighter timestep restriction; however, this only decreases the restriction by some constant factor, which is independent of the order of approximation.  In all experiments, the increase in the order of the Taylor expansion did not correspond with a significant decrease in error.  
\begin{algorithm}
\begin{algorithmic}[1]
\Procedure{Two-dimensional temporal Taylor series evaluation }{}
\State $\mb{w} = \tilde{\mb{u}}^n$
%\For{$\ell = 2N+1,\ldots, 2$}
\For{$\ell = \tilde{N},\tilde{N}-1,\ldots, 0$}
\State $\mb{w} = \mb{w} + \frac{dt}{1 + \ell} (-c \mb{D}_x - c\mb{D}_y){\tilde{\mb{u}}^n }$
\EndFor
\State $\tilde{\mb{u}}^{n+1} =  {\mb{w}}$
\EndProcedure
\end{algorithmic}
\caption{Time evolution procedure for 2D scalar advection. $\tilde{N}$ may be taken to be $d(2N+1)$ for exact time evolution.}
\label{alg:time2}
\end{algorithm}

\begin{figure}[!h]
\centering
\subfloat[$x$-reconstruction (Virtual)]{\includegraphics[width=.32\textwidth]{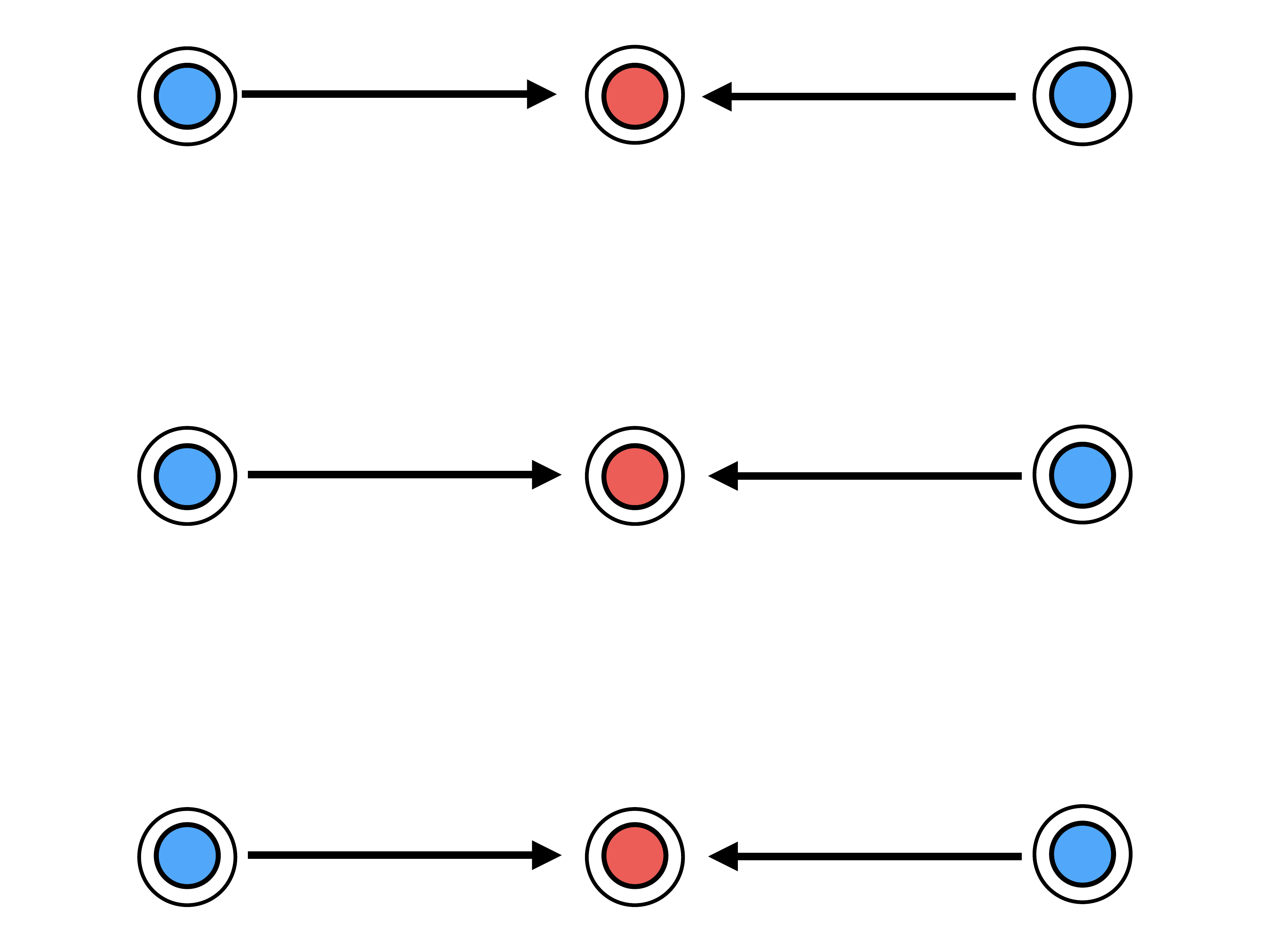}}
\subfloat[$x$-reconstruction (Central)]{\includegraphics[width=.32\textwidth]{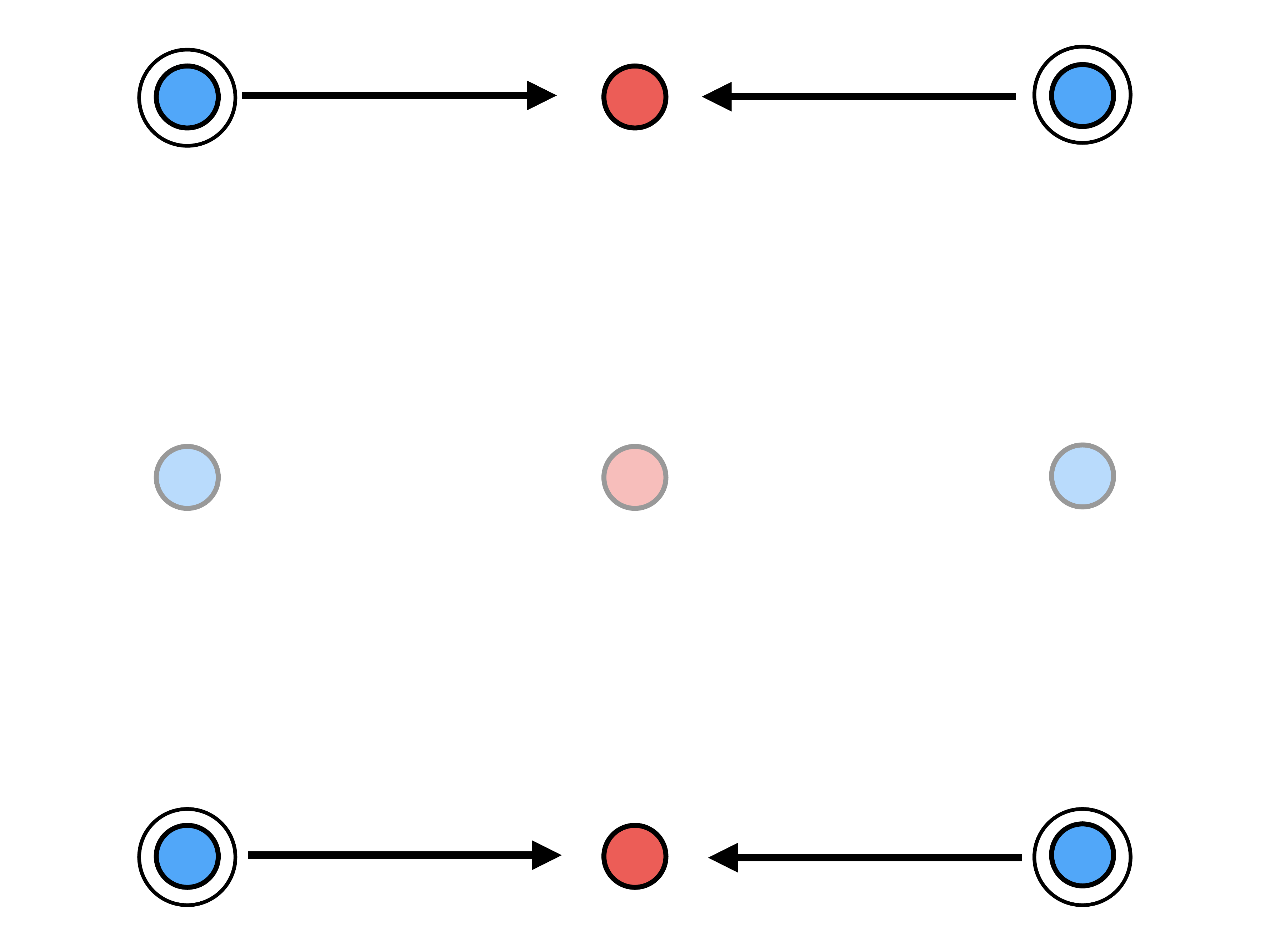}}
\subfloat[$x$-reconstruction (Upwind)]{\includegraphics[width=.32\textwidth]{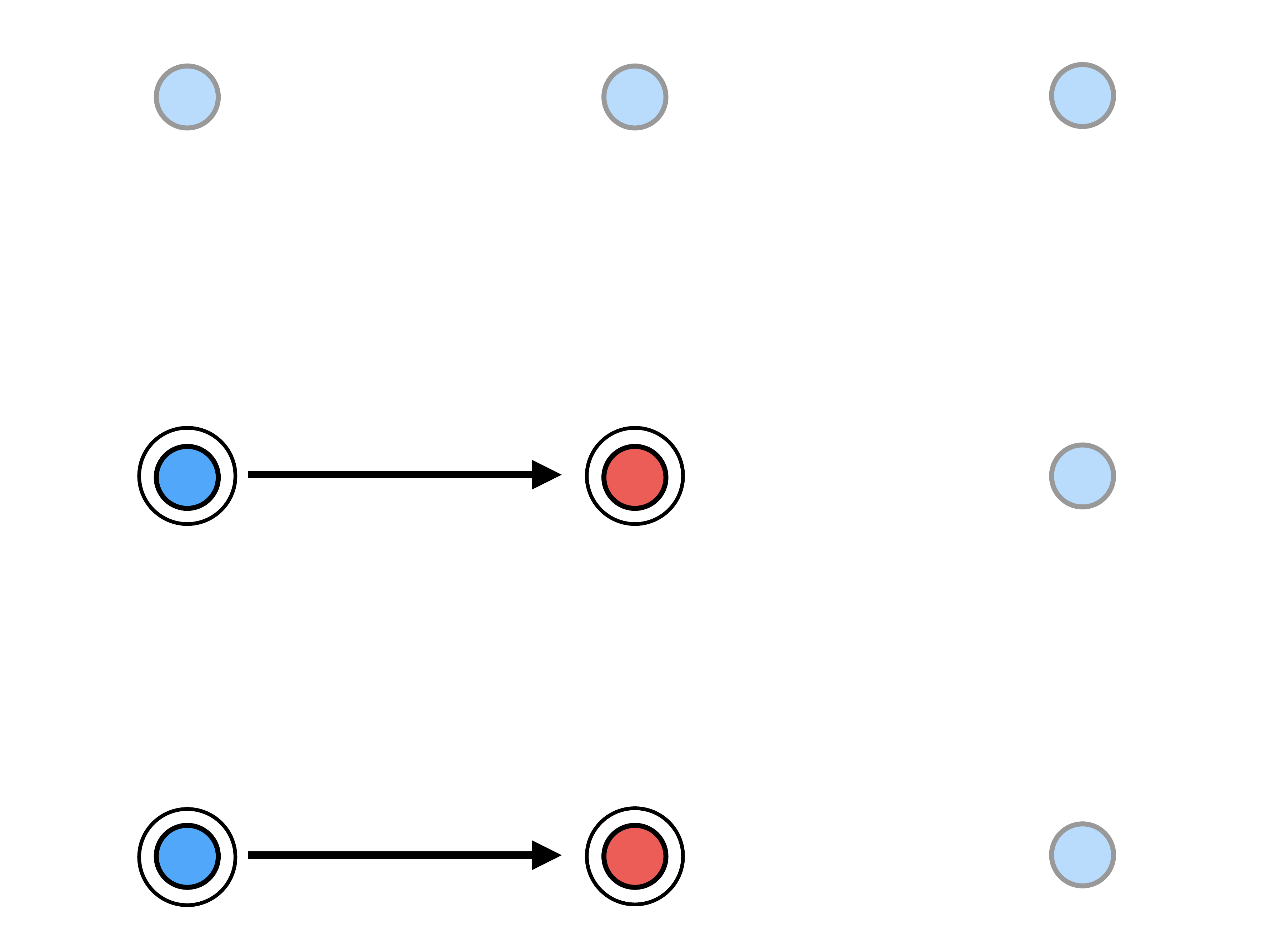}}\\
\subfloat[$y$-reconstruction (Virtual)]{\includegraphics[width=.32\textwidth]{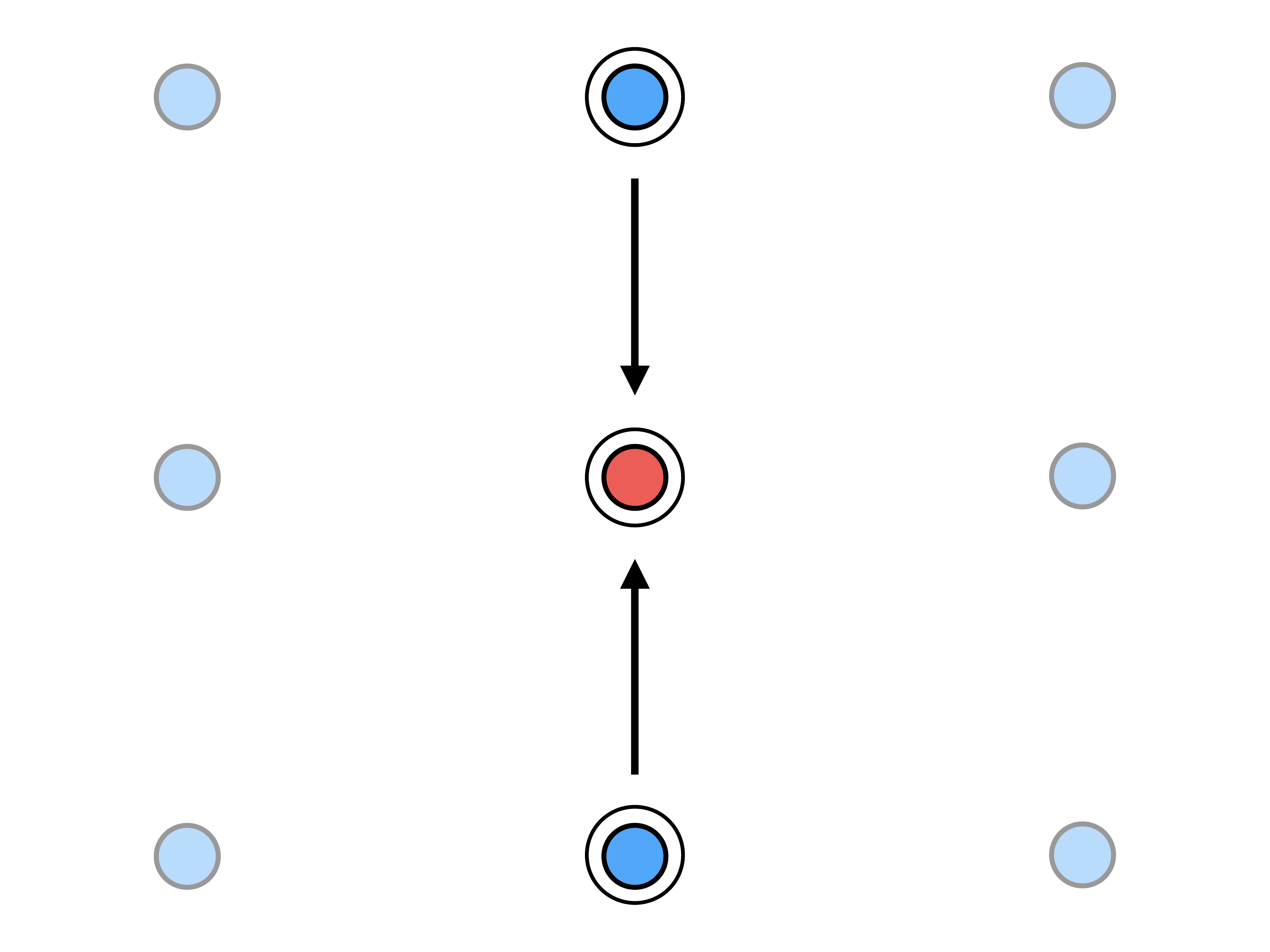}}
\subfloat[$y$-reconstruction (Central)]{\includegraphics[width=.32\textwidth]{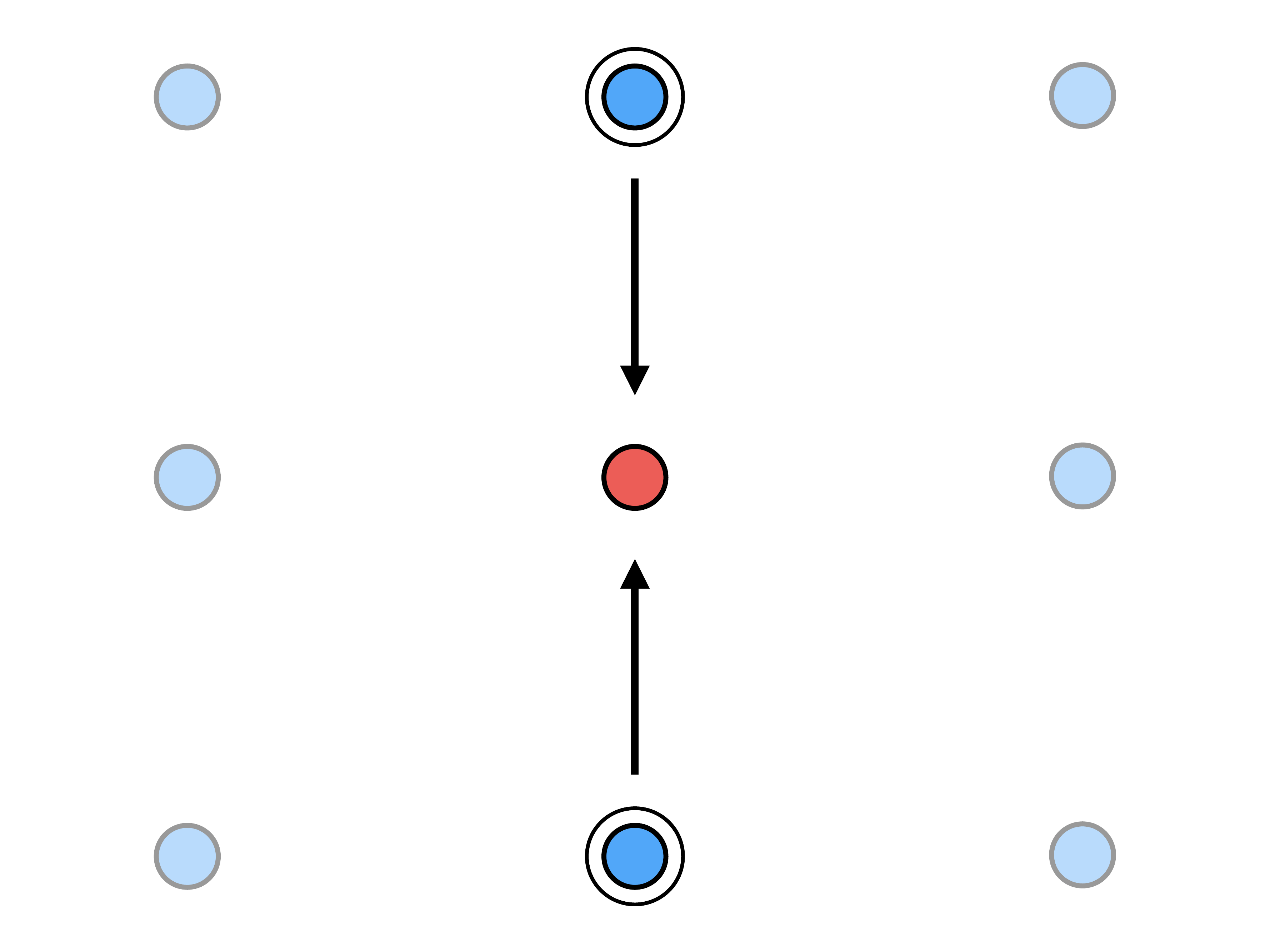}}
\subfloat[$y$-reconstruction (Upwind)]{\includegraphics[width=.32\textwidth]{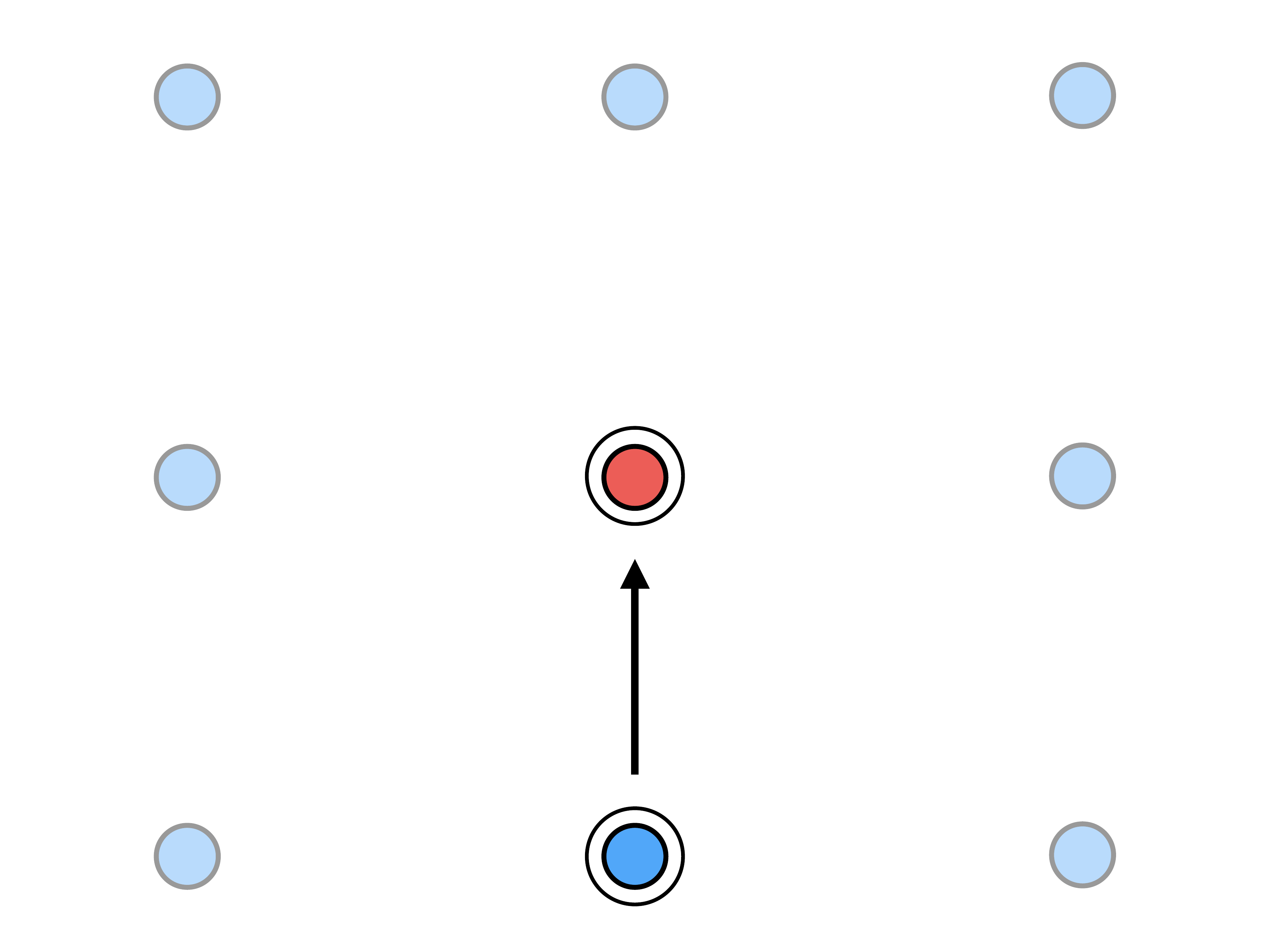}}
\caption{Two-dimensional Hermite reconstruction stencils.  The Upwind reconstruction assumes $c_x, c_y > 0$.  Nodes which contribute data to a reconstruction are circled.}
\label{fig:stencils2D}
\end{figure}

While time evolution is extended in a straightforward way regardless of spatial dimension, interpolation operators in higher dimensions are defined through applications of 1D interpolation operators along each coordinate direction.  The application of operators for the Virtual and Central Hermite methods is illustrated in Figure~\ref{fig:stencils2D}, and we refer the reader to \cite{goodrich2006hermite} for more details on the extension of the Dual Hermite method to multiple dimensions.  The Upwind Hermite reconstruction may be adapted to the advection equation in higher dimensions by considering the direction of advection along each coordinate, though the direction of the reconstruction will depend on the sign of $c_x$ at each point.  

\subsection{Numerical experiments in two dimensions}
We consider two model problems in two space dimensions: the periodic advection equation 
\[
\pd{u}{t} + c_x \pd{u}{x} + c_y \pd{u}{x} = 0.
\]
where $\mb{c} = \LRp{c_x,c_y}$ is a unit vector, and the isotropic wave equation in first order form
\begin{align*}
\frac{1}{c^2}\pd{p}{t} &= -\div\mb{u}, \qquad \pd{\mb{u}}{t} = -\grad{p},
\end{align*}
where $c$ is a specified wavespeed, $p$ is pressure and $\mb{u} = (u,v)$ is the velocity.  The CFL condition for the two-dimensional advection equation is
\[
\nor{\mb{c}} dt <  h, 
\]
while the wave equation depends on the maximum wavespeed in each coordinate direction.  For the non-dimensional isotropic wave equation above, this results in the CFL condition $cdt < h$. 

We note that behavior of the Upwind Hermite method is reported only for the advection equation. While the Upwind Hermite method may be extended to hyperbolic systems through a characteristic-based approach, numerical experiments with the isotropic wave equation indicated that the method resulted in a timestep restriction depending on the degree of approximation.  We intend to explore additional generalizations of the Upwind Hermite method to systems of equations in multiple dimensions in the future. 

\begin{figure}[!h]
\centering
\subfloat[$C = .1$]{\includegraphics[width=.475\textwidth]{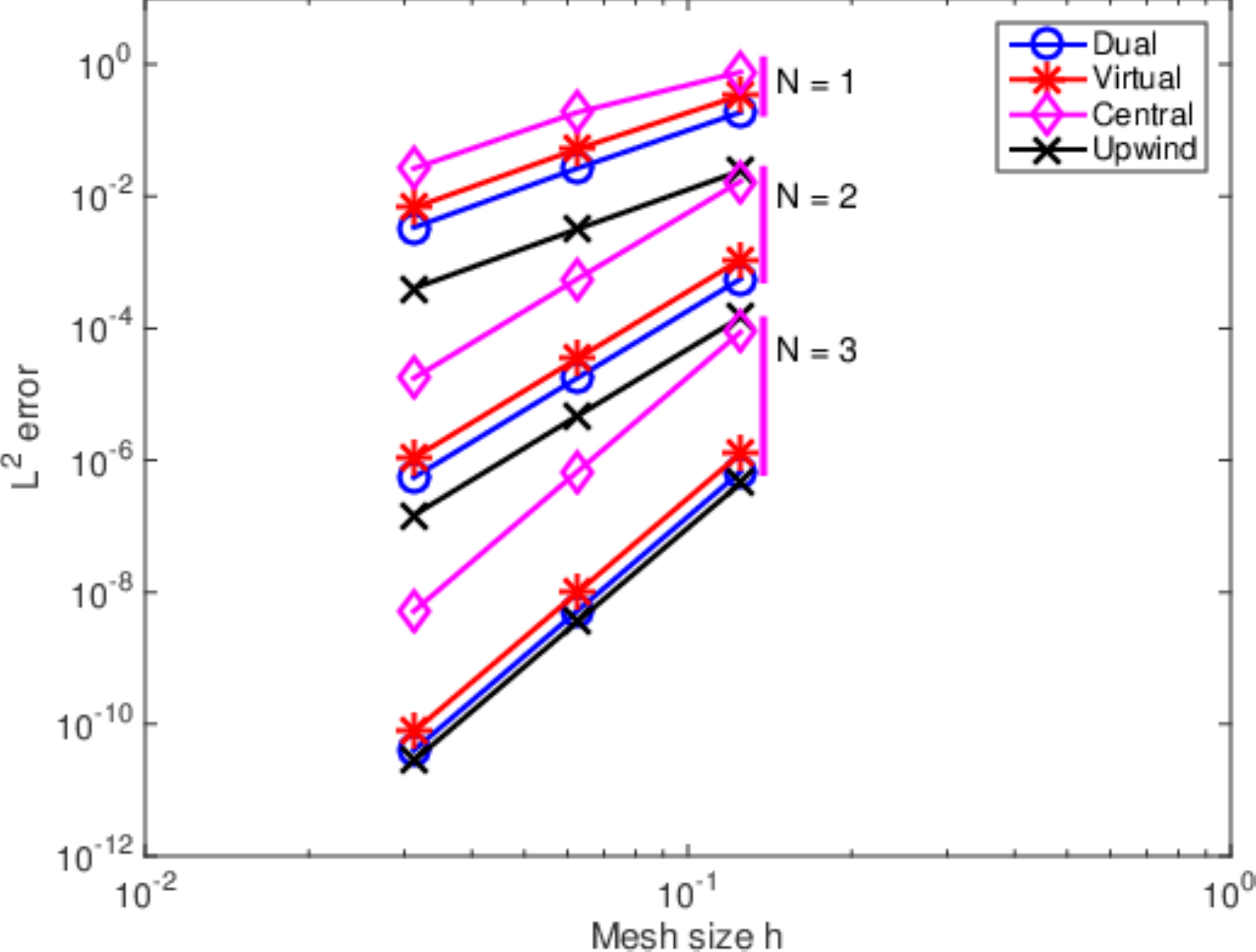}}
\subfloat[$C = .5$]{\includegraphics[width=.475\textwidth]{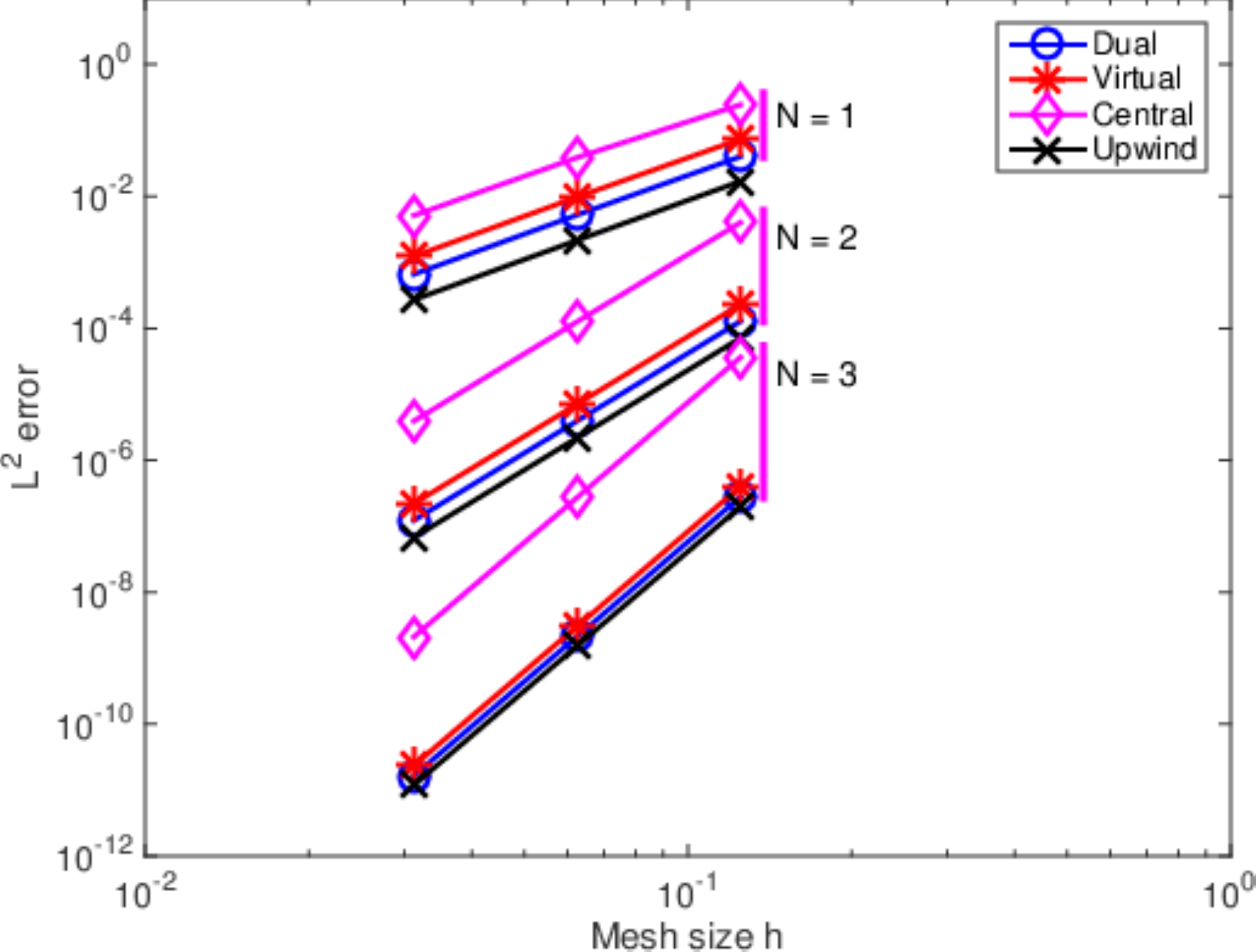}}\\
\subfloat[$C = .9$]{\includegraphics[width=.475\textwidth]{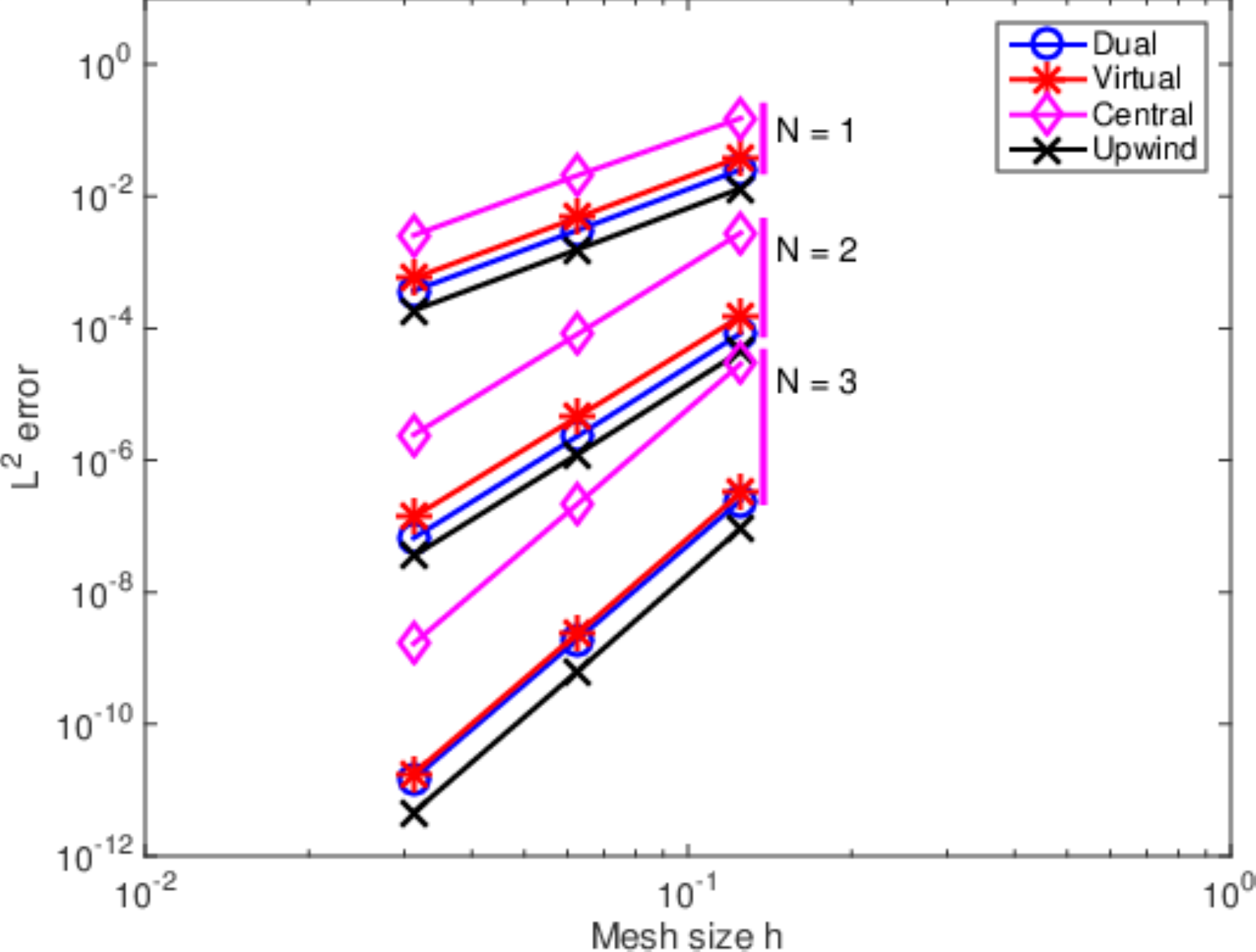}}
\caption{Convergence of $L^2$ errors for each Hermite scheme for the two-dimensional advection equation with a smooth sinusoidal solution.}
\label{fig:advec2D}
\end{figure}

As in one dimension, we introduce a CFL constant $C$ such that $dt = C h$, and examine the behavior of each Hermite method at various values $C$.  Figure~\ref{fig:advec2D} shows the $L^2$ convergence of the Virtual and Central Hermite methods for the advection equation.  We take the advection speeds $c_x, c_y$ and exact solution to be
\[
c_x = \cos(\pi/3), \quad c_y = \sin(\pi/3), \qquad u(x,y,t) = \sin(\pi(x-c_x t))\sin(\pi(y-c_y t)).
\]
and compute $L^2$ errors for isotropic grids of $8\times 8, 16\times 16$, and $32\times 32$ nodes.  
%We compute also the maximum pointwise error in the solution over all nodes in the primary grid (which we refer to this as the discrete $\ell_{\inf}$ error)
%\[
%\max_{m} \LRb{u(x_m,y_m,T) - u_h(x_m,y_m,T)},
%\]
%where $u_h(x_m,y_m,T)$ is the value of the Hermite solution at the $m$th node and final time.  

Unlike Hermite methods in one space dimension, a temporal Taylor series of degree $d(2N+1)$ is required for exact time evolution in $d$ dimensions \cite{goodrich2006hermite}.  All experiments use exact time evolution; however, decreasing the order of the temporal Taylor series to $2N+1$ did not result in a significant decrease in error, though the stable timestep restriction decreases by a factor of $.5$.  $L^2$ rates of convergence are reported in Table~\ref{table:advec2D} for $T=1$.  As in the one-dimensional case, the $L^2$ error is observed to converge at a rate close to $h^{2m+1}$.  %The convergence rates of the Dual and Virtual Hermite methods differ slightly, unlike the one-dimensional case; however, this difference is negligible.  
\begin{table}[!h]
\centering                                                          
\begin{tabular}{|c||c|c|c|c|c|c|c|c|c|}
\hline
& \multicolumn{3}{c|}{$C=.1$} &  \multicolumn{3}{c|}{$C=.5$}  &  \multicolumn{3}{c|}{$C=.9$} \\ 
\hline
$N$ & $1$ & $2$ & $3$ & $1$ & $2$ & $3$  & $1$ & $2$ & $3$ \\ 
\hhline{|=|=|=|=|=|=|=|=|=|=|}
Dual  & 2.90 & 5.00 & 7.01  & 2.97 &  5.02 & 7.02 &  3.05   & 5.15  & 7.00 \\ 
\hline                                                              
Virtual & 2.83 & 5.00 & 7.01 & 2.95  & 5.02  & 7.03 & 3.04 & 5.04  & 7.06 \\        
\hline                                                              
Central & 2.43 & 4.96 & 7.03  & 2.79 &  5.00 & 7.04 &  2.95   & 5.08  & 7.04 \\ 
\hline
Upwind &  2.96 &    4.98 &    6.99 &   2.98 &    5.02&    7.01 &    3.07 &    5.14&    7.14\\
\hline                                                             
\end{tabular}                                                       
\caption{$L^2$ rates of convergence of the Dual, Virtual, and Central Hermite methods for the advection equation in two dimensions.}
\label{table:advec2D}
\end{table}              

\begin{figure}[!h]
\centering
\subfloat[$C = .1$]{\includegraphics[width=.475\textwidth]{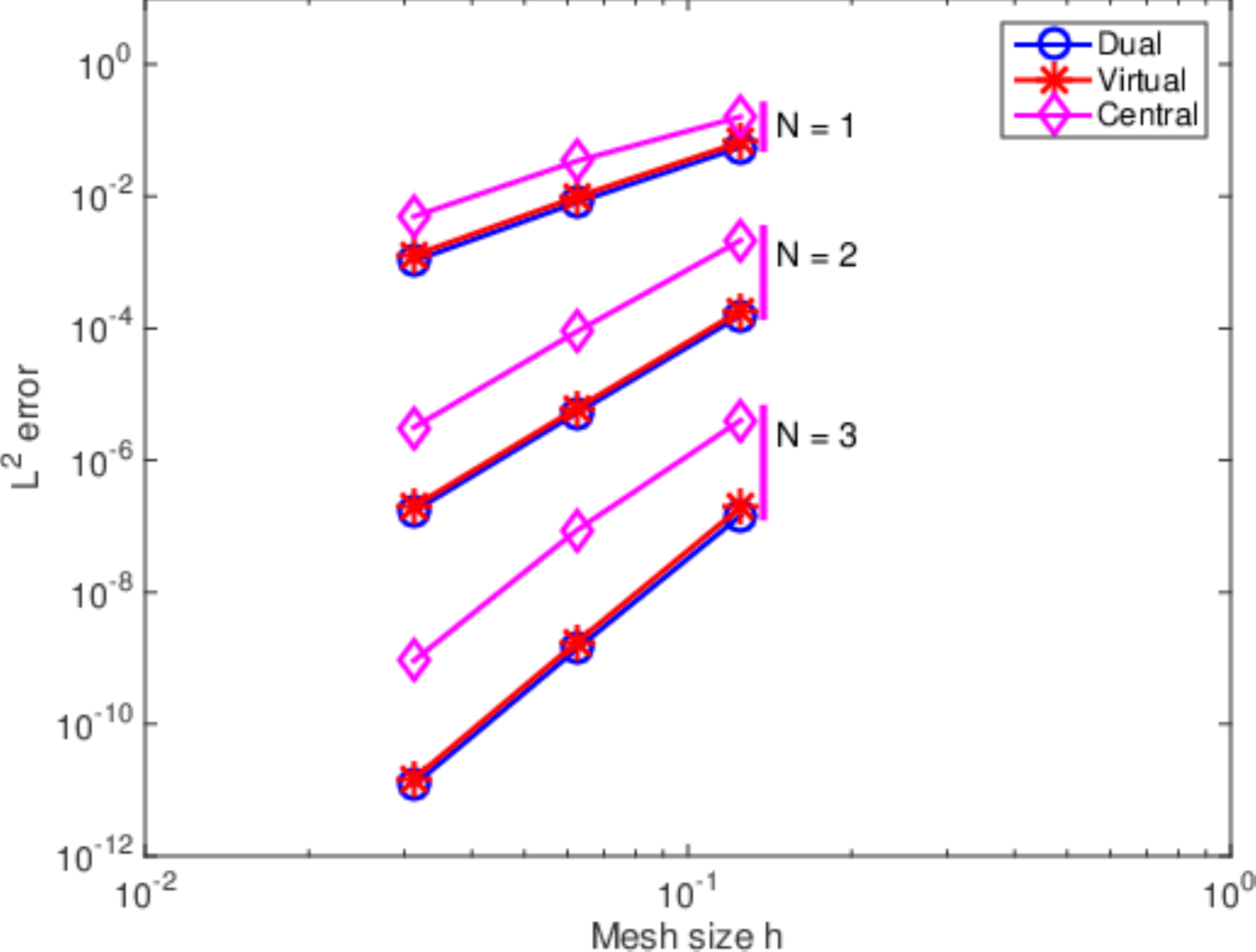}}
\subfloat[$C = .5$]{\includegraphics[width=.475\textwidth]{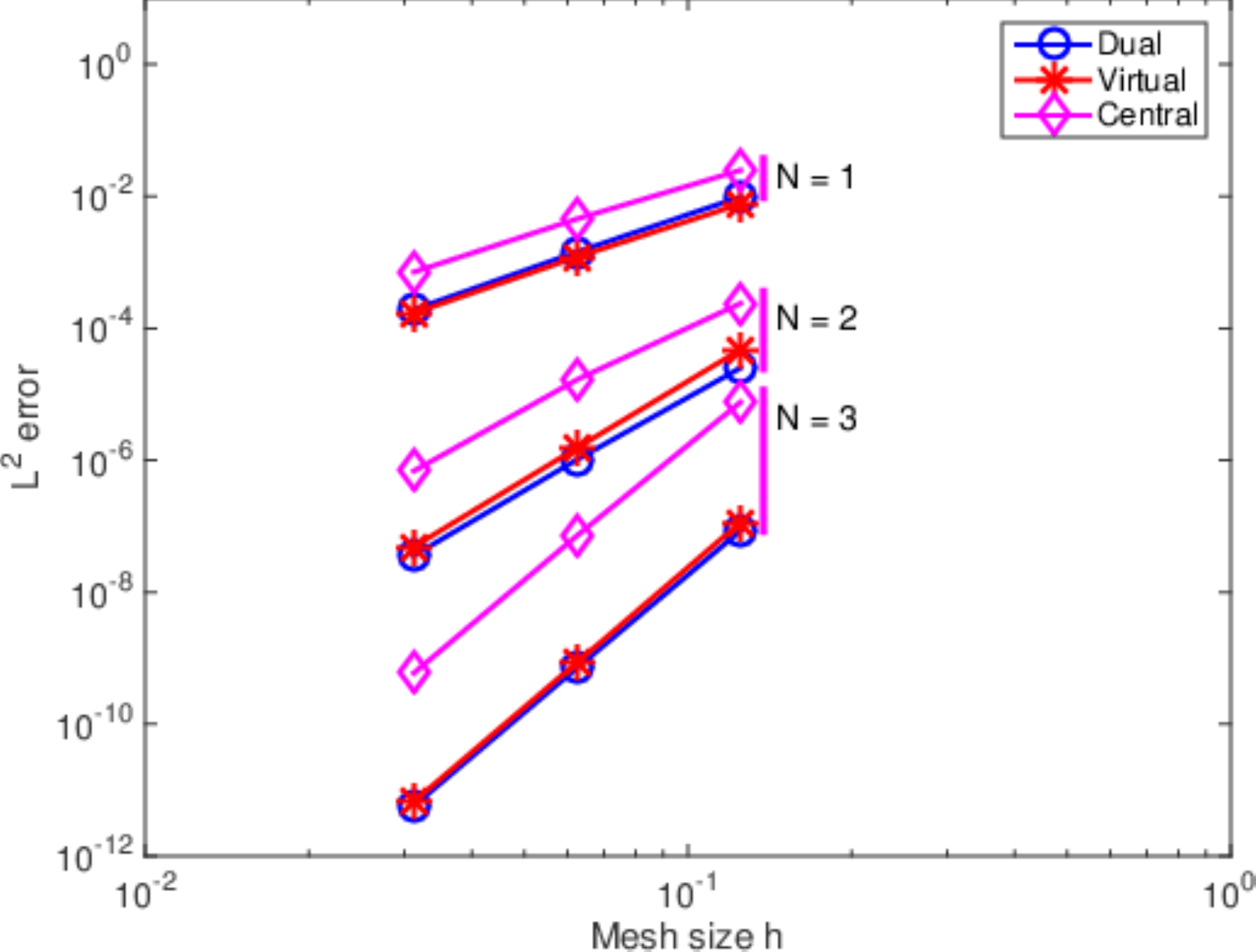}}\\
\subfloat[$C = .9$]{\includegraphics[width=.475\textwidth]{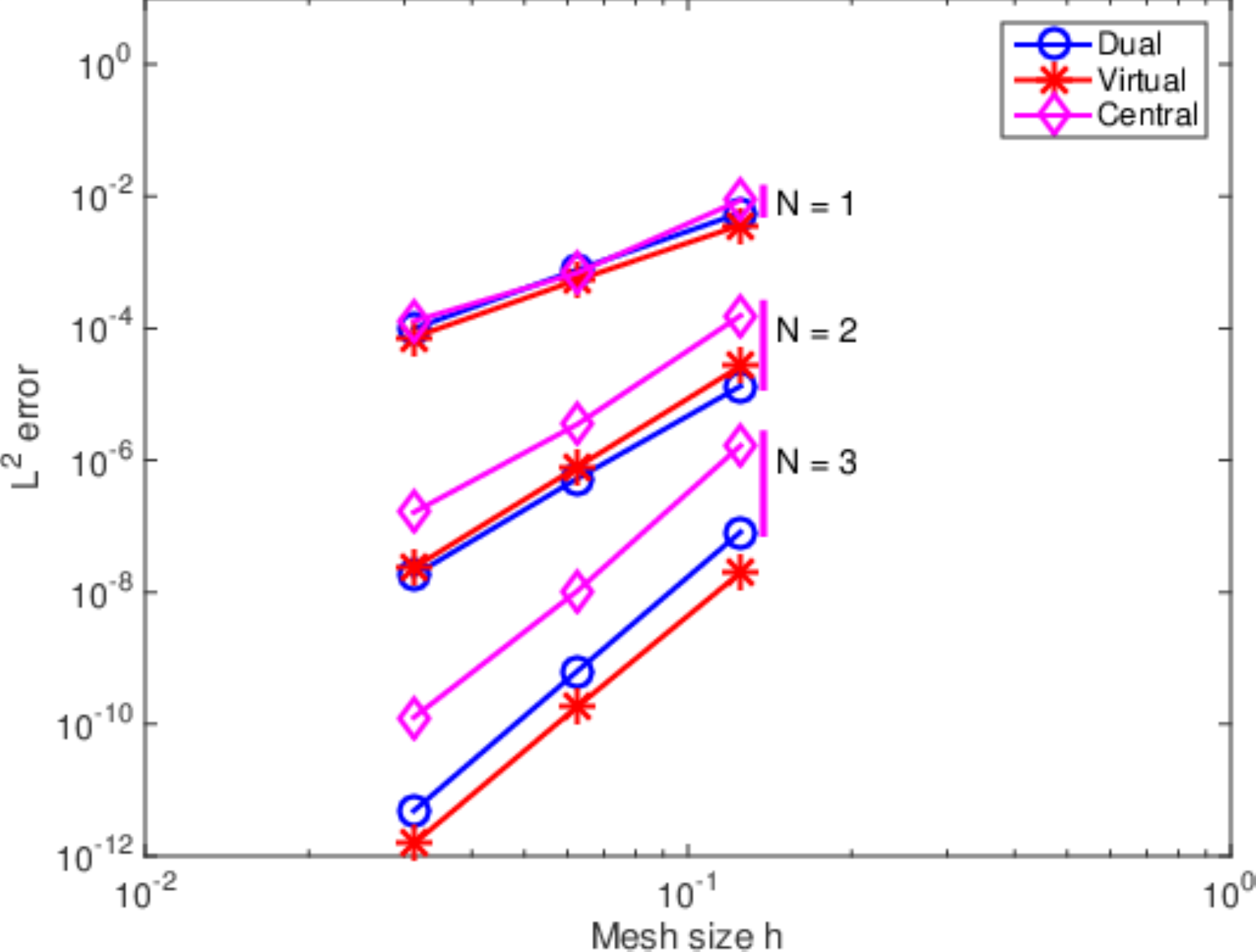}}
\caption{Convergence of $L^2$ errors for the Dual, Virtual, and Central Hermite schemes for the isotropic wave equation in two dimensions.}
\label{fig:wave2D}
\end{figure}

We repeat the Dual, Virtual, and Central Hermite convergence experiments for the periodic wave equation in 2D, using the exact solution
\[
p(x,y,t) = \sin(\pi x)\sin(\pi y)\cos(\sqrt{2}\pi t).  
\]
Upwind Hermite results are not reported, since a straightforward application of the two-dimensional upwind Hermite reconstruction does not yield a stable procedure for the wave equation.  Convergence experiments are repeated for the set of grids used for advection, and Figure~\ref{fig:wave2D} plots the $L^2$ errors in $p$ at time $T=1$ for the Dual, Virtual, and Central Hermite methods.  The Dual and Virtual Hermite methods produce errors of very similar magnitude, while the error for the Central Hermite method is larger by a factor of roughly $2^N$ as observed in 1D.  Surprisingly, at $C=.9$ and $N=3$, the error for the Virtual Hermite method is lower than that of the Dual Hermite method. 
The $L^2$ rates of convergence are reported in Table~\ref{table:wave2D}.  

\begin{table}[!h]
\centering                                                          
\begin{tabular}{|c||c|c|c|c|c|c|c|c|c|}
\hline
& \multicolumn{3}{c|}{$C=.1$} &  \multicolumn{3}{c|}{$C=.5$}  &  \multicolumn{3}{c|}{$C=.9$} \\ 
\hline
$N$ & $1$ & $2$ & $3$ & $1$ & $2$ & $3$  & $1$ & $2$ & $3$ \\ 
\hhline{|=|=|=|=|=|=|=|=|=|=|}
Dual  & 2.86 & 4.93 & 6.79  & 2.84 &  4.73 & 6.92 &  2.91   & 4.77  & 7.02 \\ 
\hline                                                              
Virtual & 2.85 & 4.93 & 6.84 & 2.75  & 4.96  & 7.03 & 2.82 & 5.07  & 6.83 \\        
\hline                                                              
Central & 2.51 & 4.71 & 6.05  & 2.56 &  4.22 & 6.82 &  3.05   & 4.95  & 6.84 \\ 
\hline                                                             
\end{tabular}                                                                
\caption{$L^2$ rates of convergence for the isotropic wave equation in 2D.}
\label{table:wave2D}
\end{table}              

\subsection{Coupling with Discontinuous Galerkin methods}

Hermite methods may also be coupled to Discontinuous Galerkin methods in order to tackle more complicated geometries and boundary conditions.  In \cite{chen2014hybrid}, coupling conditions between DG and the Hermite method are constructed for both the primary and auxiliary grids using a least squares reconstruction and high order finite difference stencils.  Since the Hermite methods introduced in this work do not require staggered grids, the transfer of information is simplified.  We will refer to the order of approximation for the DG method as $m$.  
\begin{figure}
\centering
\subfloat[Hermite-DG coupling (4-stage RK with 2 substeps)]{\includegraphics[width=.425\textwidth]{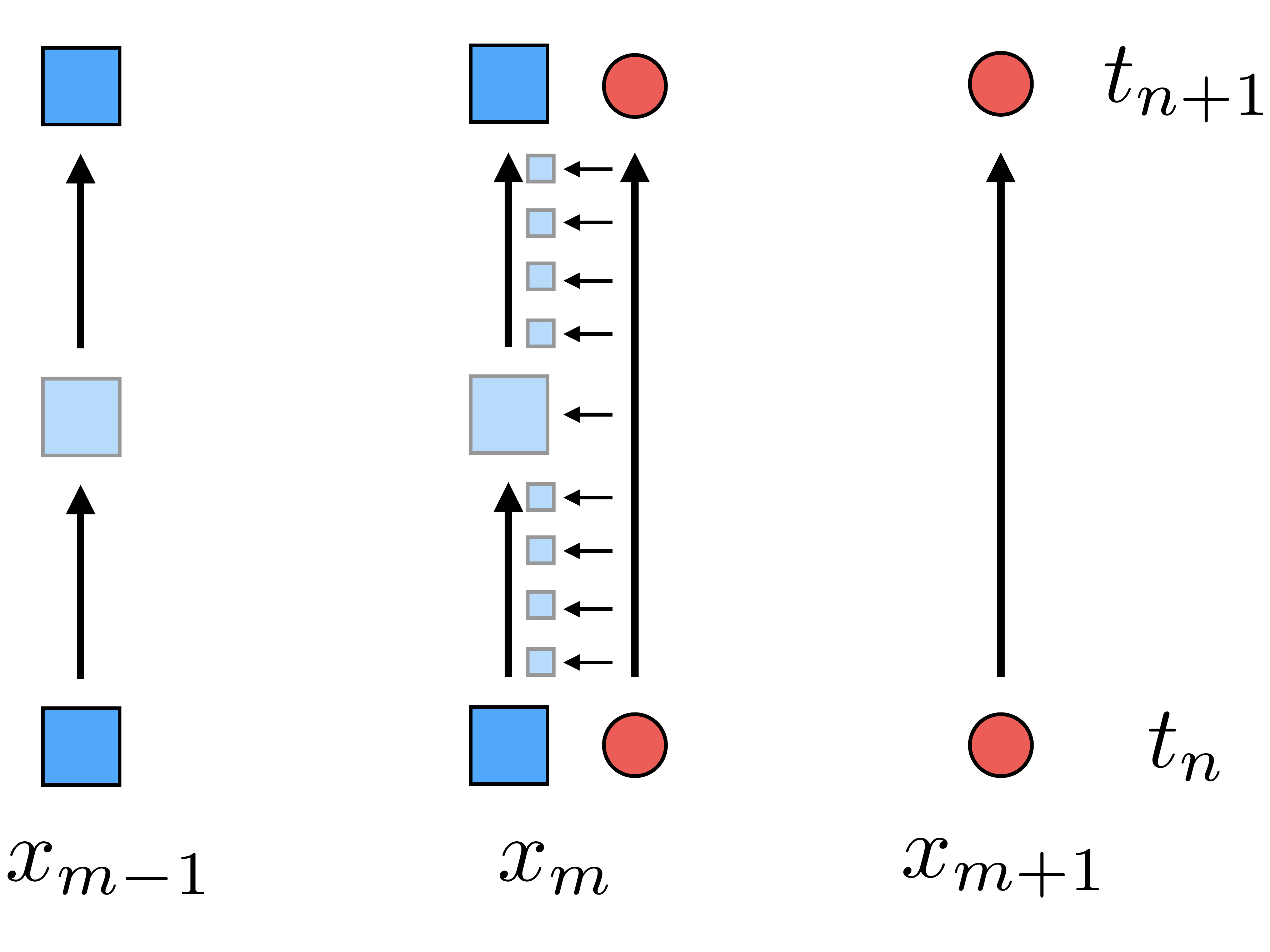}}
\hspace{3em}
\subfloat[DG-Hermite coupling (patch recovery)]{\includegraphics[width=.425\textwidth]{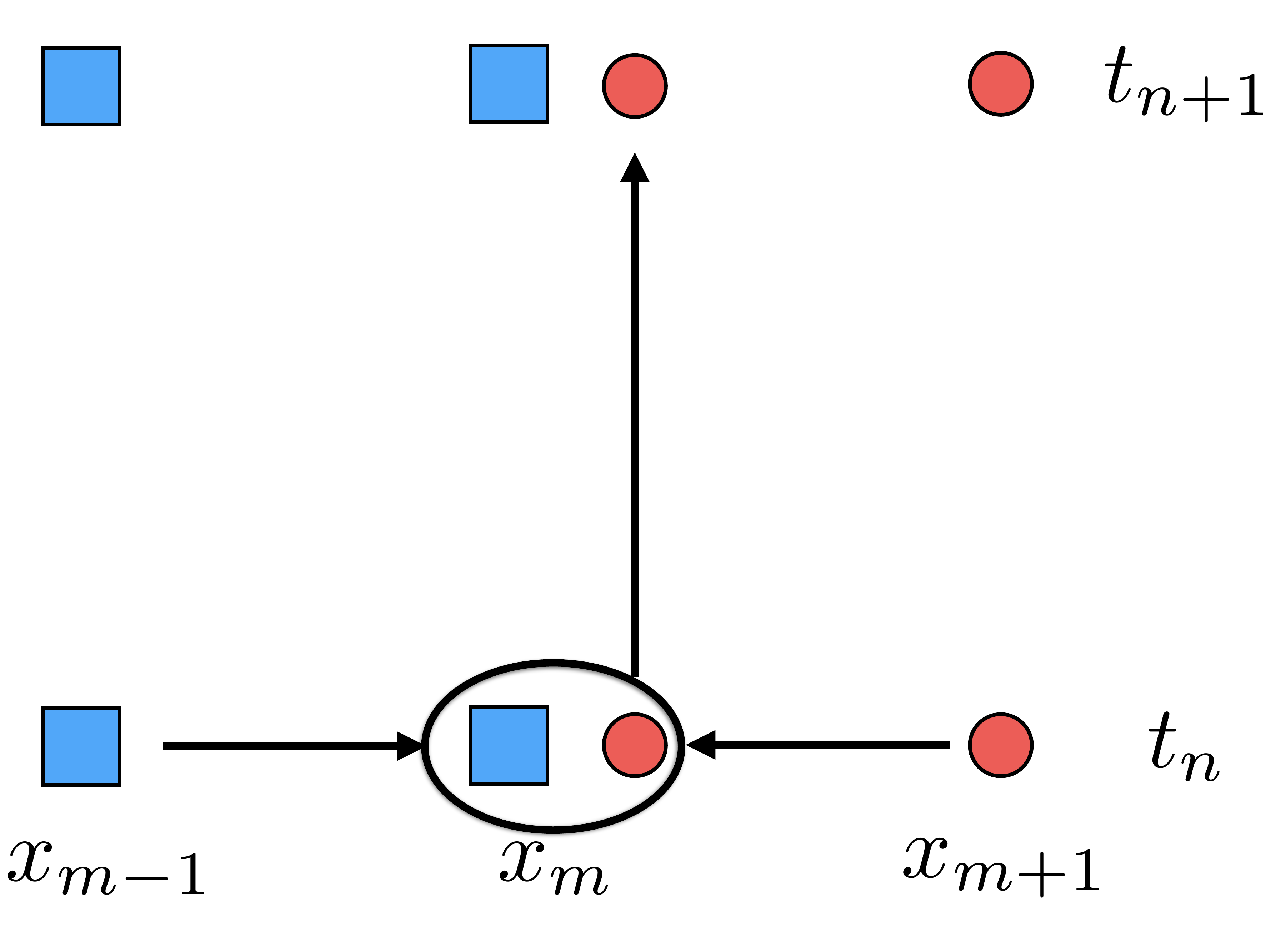}}
\caption{Coupling between Hermite and DG methods without staggered grids (DG nodes are squares, while Hermite nodes are circles).  }
\label{fig:coupling}
\end{figure}
Hermite methods may transfer information to DG methods through the numerical flux.  Due to a timestep restriction of $O(h/m^2)$ for DG compared to the $O(h)$ timestep restriction for Hermite methods, multiple DG substeps must be taken for each Hermite timestep.   The following numerical experiments use a 4th-order Runge-Kutta scheme with 5 stages, and necessitates the evaluation of the numerical flux for each stage.  To maintain high order convergence, we compute flux contributions by evaluating the high order Hermite interpolant in time, as shown in Figure~\ref{fig:coupling}.  The coupling from DG to Hermite is more fragile, as high order derivative information must be determined from the DG solution.  In the following experiments, these derivatives are provided via a patch reconstruction at Hermite nodes \cite{zienkiewicz1992superconvergent}. Projection onto the Hermite basis directly yields high order Hermite coefficients; alternatively, these coefficients may then be determined by taking derivatives of the reconstructed polynomial.  

We compute $L^2$ errors and convergence rates using the smooth solution 
\[
p(x,y,t) = \sin(2\pi x)\sin(2\pi y) \cos(2\pi t).
\]
for wavespeed $c = \sqrt{1/2}$.  On non-overlapping grids, a polynomial of degree $2N+1$ is constructed using both Hermite and DG data on neighboring elements.  Since the best possible convergence rate for DG is $O(h^{m+1})$,\footnote{Optimal convergence rates for upwind DG are typically observed in practice, and are provable on specific classes of meshes \cite{cockburn2008optimal}.  However, on general meshes, DG methods can expect at most $O(h^{m+1/2})$. } we take the order of approximation for DG to be $m=2N$ (where $N$ is the degree of the Hermite method) in order to preserve the $O(2N+1)$ convergence rate.  
\begin{figure}[!h]
\centering
\subfloat[Coupled mesh]{\includegraphics[width=.44\textwidth]{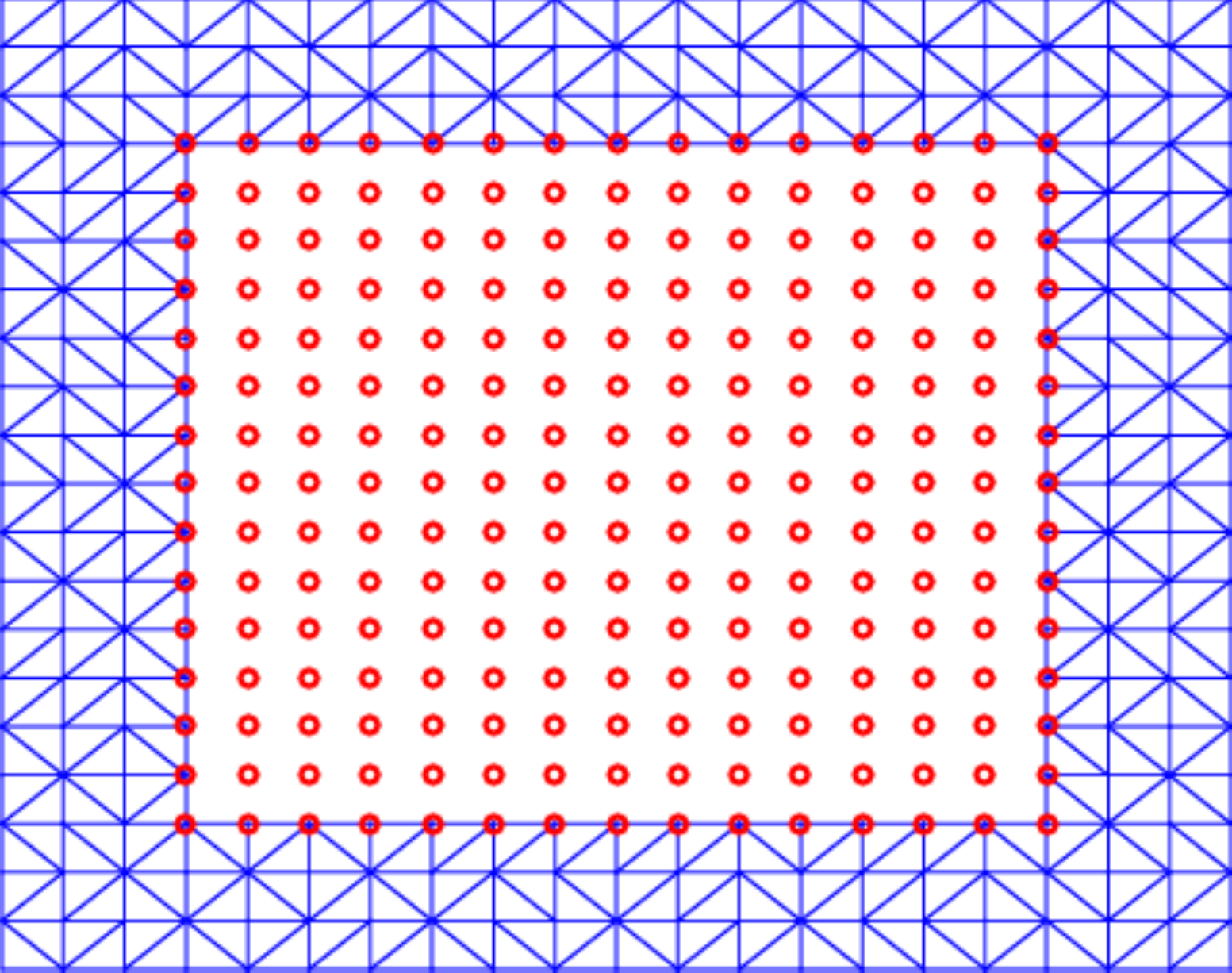}}
\hspace{2em}
\subfloat[Wave equation $L^2$ errors]{\includegraphics[width=.425\textwidth]{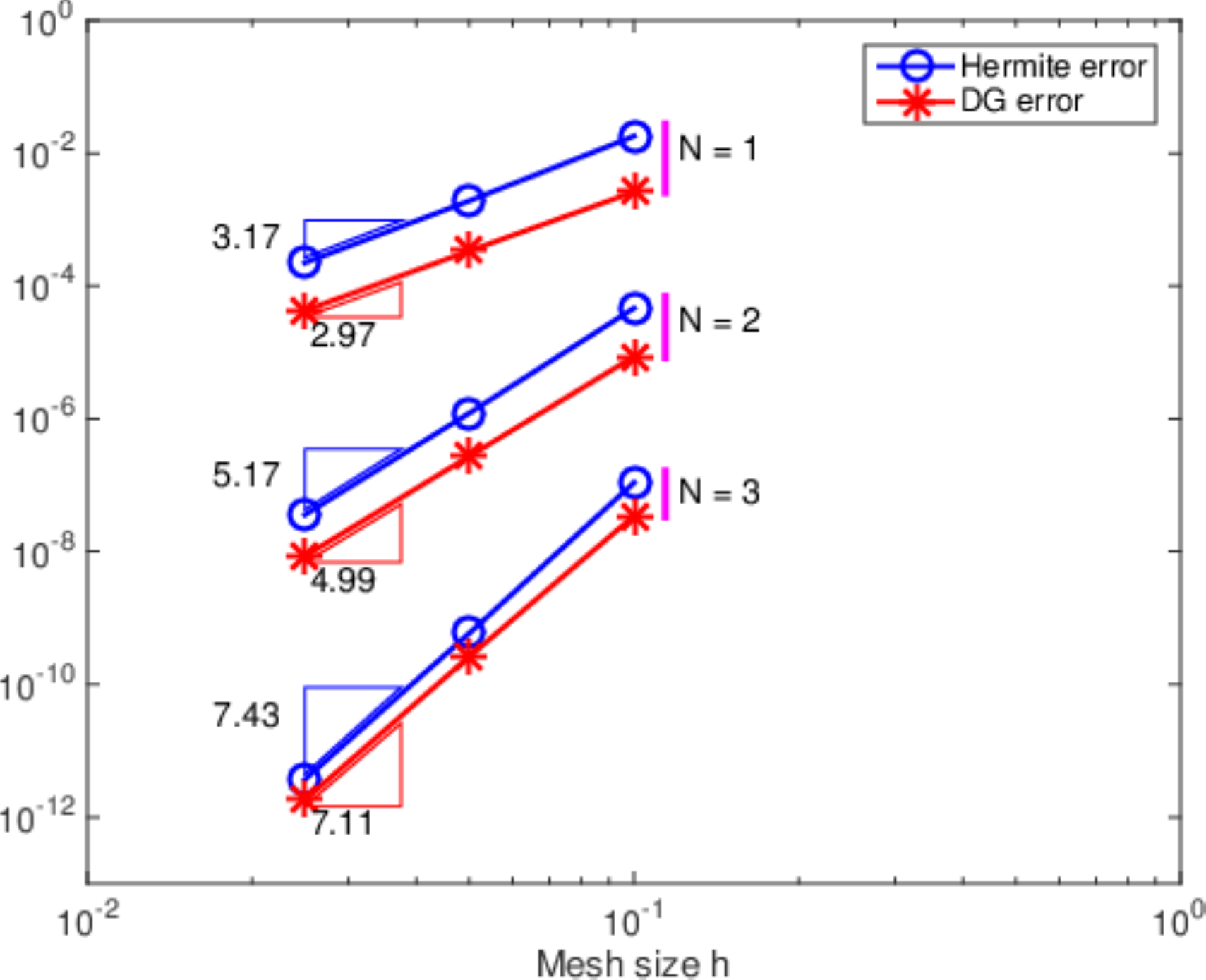}}
\caption{(Left) Coarsest coupled mesh used for coupling Virtual Hermite and DG (Hermite nodes are circled). (Right) Convergence of $L^2$ errors in 2D for the isotropic wave equation.  }
\label{fig:HDG}
\end{figure}

Results are shown in Figure~\ref{fig:HDG} for $N=1,2,3$ for the coupled Virtual Hermite-DG method, using $C=.8$ for the Hermite timestep.  For $N>2$, convergence was limited by the order of the DG Runge-Kutta scheme, and a smaller timestep must be taken to recover optimal convergence rates.  Similar observations were made when coupling Dual Hermite and DG methods \cite{chen2014hybrid}.  Similar behavior is observed when coupling Central Hermite and DG, though the Hermite error increases slightly.  

Unfortunately, the approximation Hermite coefficients using DG becomes less accurate at higher orders, as roughly an order of convergence is lost per derivative with patch recovery methods.  A salient alternative to patch recovery is Smoothness-Increasing Accuracy-Conserving (SIAC) postprocessing \cite{cockburn2003enhanced, ryan2005extension}, which produces smooth reconstructions of the solution which converge with rate $O(h^{2m+1})$.  Under such a method, optimal convergence rates could be preserved using DG and Hermite methods with degrees $m=N$.  The postprocessing of higher order derivatives may also yield additional accuracy in Hermite coefficients \cite{ryan2009local}.  %It may also be beneficial to more directly couple DG and Hermite methods by using the same time-evolution scheme.  For example, for linear autonomous systems, this might be done by postprocessing the DG solution to obtain derivative information \cite{Li2015} for Algorithm~\ref{alg:time} or ADER schemes \cite{titarev2002ader}.  

%\textcolor{red}{Mention simplified coupling using one-step Hermite methods + overlapping grids.   Dual-grid Hermite uses least squares to determine boundary valued DG derivatives.  Add various approaches to improving coupling (smooth reconstructions \cite{cockburn2003enhanced, ryan2005extension} and derivative reconstructions).  Mention possibility of more direct coupling based on smooth reconstructed DG data \cite{titarev2002ader,gassner2013analysis} }. 

\section{Conclusions and future work}

We have presented a generalization of Hermite methods for periodic problems, and have investigated two new methods within this framework and compared their performance to the original Hermite method in the literature.  The original Dual Hermite method results in a two-node stencil in one space dimension, and requires time integration on both primal and staggered (dual) grids.  The Virtual Hermite method increases the stencil to three nodes, but avoids explicit storage of staggered grid degrees of freedom by fusing operations on the auxiliary and primary grid together.  The Central Hermite method modifies the interpolation procedure in order to avoid a staggered grid, and in doing so, maintains a two-node stencil and doubles the timestep restriction.  However, to achieve a specific error resolution, the Central Hermite method requires almost as many degrees as the original Dual Hermite method.  Additionally, the Central Hermite method may suffer from the propagation of spurious modes.  The Upwind Hermite method achieves a resolution close to that of the Dual Hermite method, while maintaining the same two-node stencil and timestep restriction as the Central Hermite method.  However, the stability of Upwind Hermite schemes does not appear to generalize in a straightforward manner to systems of equations in higher dimensions.  Finally, since the Virtual, Central, and Upwind Hermite methods do not require dual grids, the coupling between Hermite and DG methods is simplified.  

Future work will address variable coefficient problems and explore stable extensions of upwind Hermite reconstructions to multi-dimensional wave problems, and to use these simplified methods to produce efficient many-core parallel implementations in two and three space dimensions.  

\section{Acknowledgments}

The authors wish to acknowledge the Matlab codes of the Hermite training library CHIDES (\href{http://www.chides.org}{http://www.chides.org}), as well as helpful discussions with Daniel Appello.  

Arturo Vargas is supported by an NSF graduate fellowship.  Jesse Chan and T. Warburton are supported by NSF (award number DMS-1216674).  Thomas Hagstrom is supported by NSF (award number DMS-1418871).  

\bibliographystyle{plain}
\bibliography{hermite}

\end{document}